\documentclass[a4paper,12pt,amsart,frenchb]{article}

\usepackage{amsmath,amsbsy,amsfonts,amssymb}
\usepackage[french]{babel}
\newcommand{\ass}[2]{\vskip0.3cm\noindent
{\bf {#1}}. { \sl {#2}}\vskip0.3cm\noindent
}

\begin{document}

  \title{  Caract\`eres de repr\'esentations de niveau $0$}
\author{J.-L. Waldspurger}
\date{1er f\'evrier 2016}
\maketitle

{\bf Introduction}

Soient $F$ un corps local non-archim\'edien et $G$ un groupe alg\'ebrique connexe d\'efini sur $F$. Notons $G_{AD}$ le groupe adjoint. Bruhat et Tits ont d\'efini l'immeuble $Imm(G_{AD})$ de ce groupe. Il est muni d'une d\'ecomposition en facettes et d'une action de $G(F)$ qui respecte cette d\'ecomposition. A toute facette ${\cal F}$ sont associ\'es plusieurs sous-groupes ouverts de $G(F)$: le groupe $ K_{{\cal F}}^{\dag}$ des \'el\'ements de $G(F)$ qui conservent ${\cal F}$;   le sous-groupe parahorique $K_{{\cal F}}^0$ d\'efini par Bruhat et Tits, qui est un sous-groupe ouvert compact de $K_{{\cal F}}$;  un certain sous-groupe  pro-$p$-unipotent $K_{{\cal F}}^+$ de $K_{{\cal F}}^{0}$ qui est distingu\'e dans $K_{{\cal F}}^{\dag}$. Le quotient ${\cal N}({\cal F})=K_{{\cal F}}^{\dag}/K_{{\cal F}}^0$ est ab\'elien et l'application de Kottwitz l'identifie \`a un sous-groupe d'un certain groupe ab\'elien  ${\cal N}$ ind\'ependant de ${\cal F}$. Pour $\nu\in {\cal N}({\cal F})$, notons $K_{{\cal F}}^{\nu}$ l'image r\'eciproque de $\nu$ dans $K_{{\cal F}}$. L'action sur l'immeuble d'un \'el\'ement $k\in K_{{\cal F}}^{0}$ conserve ${\cal F}$ point par point. Si $\nu\not=0$, l'action sur l'immeuble d'un \'el\'ement $k\in K_{{\cal F}}^{\nu}$ conserve ${\cal F}$ mais pas point par point en g\'en\'eral. On sait que ${\cal F}$ s'identifie \`a un ouvert d'un espace vectoriel r\'eel euclidien et qu'il existe une isom\'etrie $\sigma_{{\cal F},\nu}$ de cet espace, conservant ${\cal F}$ et telle que tout \'el\'ement $k\in K_{{\cal F}}^{\nu}$ agisse dans ${\cal F}$ par cette isom\'etrie. On note ${\cal F}^{\nu}$ le sous-ensemble des points fixes de $\sigma_{{\cal F},\nu}$ dans ${\cal F}$. On note $Fac^*_{max}(G)$ l'ensemble des couples $({\cal F},\nu)$ form\'es d'une facette ${\cal F}$   de $Imm(G_{AD})$ et d'un \'el\'ement $\nu\in {\cal N}({\cal F})$, tels que ${\cal F}^{\nu}$ est r\'eduit \`a un point.

Soit $\pi$ une repr\'esentation lisse de longueur finie de $G(F)$ dans un espace complexe $V$. Pour toute facette ${\cal F}$, notons $V^{K_{{\cal F}}^+}$ le sous-espace des invariants par $K_{{\cal F}}^+$. Il est de dimension finie. On dit que $\pi$ est de niveau $0$ si $V$ est engendr\'e  par les sous-espaces $V^{K_{{\cal F}}^+}$ quand ${\cal F}$ d\'ecrit les facettes de l'immeuble. Supposons $\pi$ de niveau $0$. Pour toute facette ${\cal F}$, le groupe $K_{{\cal F}}^{\dag}$ conserve $V^{K_{{\cal F}}^+}$ et cette action se quotiente en une repr\'esentation $\pi_{{\cal F}}$ de dimension finie  du groupe  $K_{{\cal F}}^{\dag}/K_{{\cal F}}^+$. Notons $trace\,\pi_{{\cal F}}$ le caract\`ere de cette repr\'esentation. On sait qu'il existe un groupe r\'eductif ${\bf G}_{{\cal F  }}$ sur ${\mathbb F}_{q}$ (le corps r\'esiduel de $F$) tel que $K_{{\cal F}}^0/K_{{\cal F}}^+\simeq {\bf G}_{{\cal F  }}({\mathbb F}_{q})$. Pour $\nu\in {\cal N}({\cal F})$,  il existe un "espace tordu" ${\bf G}_{{\cal F}}^{\nu}$ sous ${\bf G}_{{\cal F}}$ de sorte que ${\bf G}_{{\cal F  }}^{\nu}({\mathbb F}_{q})\simeq K_{{\cal F}}^{\nu}/K_{{\cal F}}^+$.  En utilisant la th\'eorie des repr\'esentations des  espaces tordus, on d\'efinit pour tout $\nu$ la projection cuspidale $\phi_{\pi,{\cal F},\nu,cusp
}$ de la restriction du caract\`ere $trace\,\pi_{{\cal F}}$ \`a $K_{{\cal F}}^{\nu}/K_{{\cal F}}^+$. On peut consid\'erer que c'est une fonction sur $G(F)$, \`a support dans $K_{{\cal F}}^{\nu}$ et biinvariante par $K_{{\cal F}}^+$. Elle est invariante par conjugaison par $K_{{\cal F}}^{\dag}$. Fixons une mesure de Haar sur $G(F)$. Pour toute fonction $f\in C_{c}^{\infty}(G(F))$, on pose
$$\Theta^{G}_{\pi,cusp}(f)=\sum_{({\cal F},\nu)\in Fac^*_{max}(G)}\int_{G(F)}f(g)\phi_{\pi,{\cal F},\nu,cusp}(g)\,dg.$$
Cette d\'efinition a un sens car on montre qu'il n'y a qu'un nombre fini de $({\cal F},\nu)\in Fac^*_{max}(G)$ tels que l'int\'egrale ci-dessus ne soit pas nulle. Cela d\'efinit une distribution $\Theta^{G}_{\pi,cusp}$ sur $G(F)$, qui est invariante par conjugaison. Pour tout sous-groupe de Levi $M$ de $G$ (d\'efini sur $F$), notons ${\cal L}(M)$ l'ensemble des sous-groupes de Levi de $G$ contenant $M$. Fixons un sous-groupe de Levi minimal $M_{min}$ de $G$   et, pour tout  $M\in {\cal L}(M_{min})$, notons $W^M$ le groupe de Weyl de $M$ relatif \`a $M_{min}$. Fixons un sous-groupe parabolique $P$ de composante de Levi $M$ et des mesures de Haar sur $M(F)$ et sur $U_{P}(F)$, o\`u $U_{P}$ est le radical unipotent de $P$. On d\'efinit la repr\'esentation $\pi_{P}$ de $M(F)$ dans le module de Jacquet de $\pi$ relatif \`a $P$. Elle est de niveau $0$.  On d\'efinit  la distribution $\Theta^{M}_{\pi_{P},cusp}$ sur $M(F)$. On l'induit \`a $G(F)$ de la fa\c{c}on habituelle. C'est-\`a-dire, pour toute $f\in C_{c}^{\infty}(G(F))$ et tout  $g\in G(F)$, notons $(^gf)_{[P]}$ la fonction $m\mapsto \delta_{P}(m)^{1/2}\int_{U_{P}(F)}f(g^{-1}mug)\,du$ sur $M(F)$, o\`u $\delta_{P}$ est le module usuel. Des mesures fix\'ees se d\'eduit une (pseudo-)mesure sur $P(F)\backslash G(F)$ (pseudo car elle ne s'applique pas \`a des fonctions sur ce quotient mais \`a des sections d'un certain fibr\'e). On pose
$$ind_{M}^G(\Theta^{M}_{\pi_{P},cusp})(f)=\int_{P(F)\backslash G(F)}\Theta^{M}_{\pi_{P},cusp}((^gf)_{[P]})\,dg.$$
Les distributions $\Theta^{M}_{\pi_{P},cusp}$ et $ind_{M}^G(\Theta^{M}_{\pi_{P},cusp})$ d\'ependent du choix du sous-groupe parabolique $P$. Mais disons qu'un \'el\'ement $g\in G(F)$ est compact modulo $Z(G)$ si l'image de $g$ dans $G_{AD}(F)$ appartient \`a un sous-groupe compact de ce groupe. Alors la restriction de $\Theta^{M}_{\pi_{P},cusp}$ aux \'el\'ements de $M(F)$ qui sont compacts modulo $Z(G)$ ne d\'epend pas de $P$ et la restriction de $ind_{M}^G(\Theta^{M}_{\pi_{P},cusp})$  aux \'el\'ements de $G(F)$ qui sont compacts modulo $Z(G)$ n'en d\'epend pas non plus.

On d\'efinit le caract\`ere-distribution $\Theta_{\pi}$ par $\Theta_{\pi}(f)=trace\,\pi(f)$ pour toute $f\in C_{c}^{\infty}(G(F))$.   Notre r\'esultat principal est le suivant.

\ass{Th\'eor\`eme}{ Soit $\pi$ une repr\'esentation lisse de $G(F)$ de longueur finie et de niveau $0$. Soit $f\in C_{c}^{\infty}(G(F))$, supposons que le support de $f$ est form\'e d'\'el\'ements de $G(F)$ qui sont compacts modulo $Z(G)$. Alors on a l'\'egalit\'e
$$\Theta_{\pi}(f)=\sum_{M\in {\cal L}(M_{min})}\vert W^M\vert \vert W^G\vert ^{-1}ind_{M}^G(\Theta^{M}_{\pi_{P(M)},cusp})(f),$$
o\`u, pour tout $M$, on a fix\'e un sous-groupe parabolique $P(M)$ de composante de Levi $M$.}

Pour $\pi$ comme ci-dessus, on sait gr\^ace \`a Harish-Chandra que $\Theta_{\pi}$ est une distribution localement int\'egrable, associ\'ee \`a une fonction localement int\'egrable $\theta_{\pi}$ sur $G(F)$, localement constante sur les \'el\'ements semi-simples fortement r\'eguliers. Soit $M\in {\cal L}(M_{min})$, notons $M(F)_{ell}$ le sous-ensemble des \'el\'ements de $M(F)$ qui sont semi-simples, fortement r\'eguliers dans $G(F)$ et elliptiques dans $M(F)$. On note $A_{M}$ le plus grand sous-tore du centre de $M$ qui soit d\'eploy\'e sur $F$. Pour $({\cal F},\nu)\in Fac^*_{max}(G)$ et $m\in M(F)_{ell}$, on d\'efinit en suivant Arthur l'int\'egrale orbitale pond\'er\'ee $J_{M}^G(m,\phi_{\pi,{\cal F},\nu,cusp}) $, modulo certains choix de mesures. Fixons    un ensemble de repr\'esentants $\underline{Fac}_{max}^*(G)$ des classes de conjugaison par $G(F)$ dans $Fac_{max}^*(G)$. Le th\'eor\`eme ci-dessus se traduit plus concr\`etement par l'\'enonc\'e suivant, o\`u $D^G$ est l'habituel discriminant de Weyl.

 \ass{Th\'eor\`eme}{ Soit $\pi$ une repr\'esentation lisse de $G(F)$ de longueur finie et de niveau $0$. Soit $M\in {\cal L}(M_{min})$ et soit $m\in M(F)_{ell}$. Supposons que $m$ soit compact modulo $Z(G)$. Alors on a l'\'egalit\'e
$$\theta_{\pi}(m)= D^G(m)^{-1/2}\sum_{L\in {\cal L}(M)} (-1)^{dim(A_{M})-dim(A_{L})}$$
$$\sum_{({\cal F}_{L},\nu)\in \underline{Fac}^*_{max}(L)}  mes(A_{L}(F)\backslash K_{{\cal F}_{L}}^{\dag})^{-1}J_{M}^L(m,\phi_{\pi_{P(L)}{\cal F}_{L},\nu,cusp}).$$}

En utilisant un r\'esultat bien connu de Casselman, on obtient une expression plus g\'en\'erale pour $\theta_{\pi}(m)$ pour tout $m\in M(F)_{ell}$ (c'est-\`a-dire  sans supposer de  condition de compacit\'e modulo $Z(G)$). Cf. paragraphe 18 pour l'\'enonc\'e pr\'ecis. 

Si on se limite aux \'el\'ements de $G(F)_{ell}$, le r\'esultat ci-dessus co\"{\i}ncide avec le th\'eor\`eme  II.4.16 de \cite{SS}. Dans le cas o\`u $\pi$ est une repr\'esentation cuspidale, Arthur a montr\'e que le caract\`ere-fonction  $\theta_{\pi}$ \'etait \'egale \`a l'int\'egrale orbitale pond\'er\'ee d'un certain coefficient de $\pi$, cf. \cite{A1}. Il a obtenu une expression qui vaut plus g\'en\'eralement pour $\pi$ de la s\'erie discr\`ete (ou m\^eme elliptique en son sens), cf. \cite{A2} th\'eor\`eme  5.1. Mais il utilise alors non pas des int\'egrales orbitales pond\'er\'ees mais leurs versions invariantes, qui sont beaucoup plus difficilement calculables. Ici, on ne suppose pas que $\pi$ est de la s\'erie discr\`ete mais on impose une condition assez s\'ev\`ere: $\pi$ est de niveau $0$. 

Notre r\'esultat s'inspire beaucoup de l'article de Court\`es \cite{Co}.  Disons qu'une fonction $f\in C_{c}^{\infty}(G(F))$ est de niveau $0$ si elle est combinaison lin\'eaire de fonctions $f'$ pour lesquelles il existe une facette ${\cal F}'$ de sorte que $f'$ soit biinvariante par $K_{{\cal F}'}^+$. Court\`es a \'etabli un th\'eor\`eme similaire \`a notre premier th\'eor\`eme ci-dessus pour des distributions \`a support dans les \'el\'ements compacts modulo $Z(G)$, que l'on restreint \`a des fonctions de niveau $0$. En un sens, notre r\'esultat est dual de celui de Court\`es: ce sont les distributions que l'on suppose "de niveau $0$" et on les restreint \`a des fonctions \`a support compact modulo $Z(G)$. Si on travaillait sur l'alg\`ebre de Lie et non pas sur le groupe, il est probable que notre r\'esultat se d\'eduirait de celui de Court\`es par transformation de Fourier. En fait, nos preuves reprennent assez souvent celles de Court\`es. On ajoute de plus un ingr\'edient essentiel: la r\'esolution de $\pi$ introduite par Schneider et Stuhler, cf. \cite{SS},  dans une version am\'elior\'ee par Meyer et Solleveld, cf. \cite{MS}. Elle permet d'exprimer $\Theta_{\pi}$ \`a l'aide des fonctions $trace\,\pi_{{\cal F}}$. Le reste de la preuve n'est rien d'autre qu'un calcul combinatoire.

Je remercie beaucoup B. Lemaire pour diverses explications et r\'ef\'erences concernant la th\'eorie des immeubles. Je remercie \'egalement J.-F. Dat pour  m'avoir indiqu\'e une r\'eference tr\`es utile.

\bigskip

\section{Notations}

 Soit $G$ un groupe r\'eductif connexe d\'efini sur un corps $k$. Tous les sous-groupes alg\'ebriques de $G$ que nous consid\'ererons sont suppos\'es d\'efinis sur $k$.     On note $Z(G)$ le centre de $G$ et $A_{G}$ le plus grand sous-tore de $Z(G)$ qui est d\'eploy\'e sur $k$.   On appelle sous-groupe de Levi de $G$ une composante de Levi  d'un sous-groupe parabolique de $G$. Pour un tel sous-groupe $M$, on note ${\cal L}(M)$, resp. ${\cal P}(M)$, ${\cal F}(M)$, l'ensemble des sous-groupes de Levi contenant $M$, resp. celui des sous-groupes paraboliques de $G$ de composante de Levi $M$, resp. celui des sous-groupes paraboliques de $G$ contenant $M$.   
  Pour un sous-groupe parabolique $P$ de $G$, on note $U_{P}$ son radical unipotent.  On note $G_{AD}$ le groupe adjoint de $G$.  
  
   Les objets ci-dessus, ainsi que d'autres d\'efinis plus loin, d\'ependent du groupe $G$.   Le cas \'ech\'eant, on ajoutera un $G$ dans leur notation pour pr\'eciser ce groupe "ambiant": ${\cal L}^G(M)$ au lieu de ${\cal L}(M)$. Les objets analogues relatifs \`a un autre groupe $H$ seront alors not\'es en rempla\c{c}ant la lettre $G$ par $H$: par exemple, si $L$ est un sous-groupe de Levi de $G$ et $M$ est un sous-groupe de Levi de $L$, ${\cal L}^L(M)$ est l'ensemble des sous-groupes de Levi de $L$ contenant $M$.
  
  Un espace tordu sous $G$ est une vari\'et\'e alg\'ebrique $G^{\nu}$ sur $k$ munie de deux actions de $G$ \`a gauche et \`a droite telles que, pour chacune des actions, $G^{\nu}$ soit un espace principal homog\`ene. On impose que $G^{\nu}(k)\not=\emptyset$. Pour tout $x\in G^{\nu}$, il existe un unique automorphisme $ad_{x}$ de $G$ \'echangeant les actions de $G$ \`a gauche et \`a droite, c'est-\`a-dire tel que $ad_{x}(g)x=xg$ pour tout $g\in G$. La restriction de $ad_{x}$ \`a $Z(G)$ ne d\'epend pas de $x$ et on impose que cet automorphisme de $Z(G)$ soit d'ordre fini. 
   Soit $P$ un sous-groupe parabolique de $G$ de composante de Levi $M$. On note $P^{\nu}$ l'ensemble des $x\in G^{\nu}$ tels que $ad_{x}$ conserve $P$ et $M^{\nu}$ l'ensemble des $x\in P^{\nu}$ tels que $ad_{x}(M)=M$.    Les quatre conditions suivantes sont \'equivalentes: 
  $P^{\nu}\not=\emptyset$; $P^{\nu}(k)\not=\emptyset$; $M^{\nu}\not=\emptyset$; $M^{\nu}(k)\not=\emptyset$. Si  elles sont v\'erifi\'ees, on dit que $P^{\nu}$ est un sous-espace parabolique de $P^{\nu}$ et que $M^{\nu}$ est un sous-espace de Levi. Remarquons que $M^{\nu}$ n'est pas d\'etermin\'e par $M$, il d\'epend aussi de $P$. Par contre, $M^{\nu}$ d\'etermine $M$: $M$ est l'ensemble des $g\in G$ tels que la multiplication par $g$, \`a droite ou \`a gauche, conserve $M^{\nu}$. Dans le cas o\`u $P$ est un sous-groupe 
  parabolique minimal et $M$ est un sous-groupe de Levi minimal, les conditions pr\'ec\'edentes sont toujours v\'erifi\'ees.
     On g\'en\'eralise \`a ces espaces les notations introduites plus haut:  ${\cal L}(M^{\nu})$, ${\cal P}(M^{\nu})$ et ${\cal F}(M^{\nu})$.
  
  Si $A$ est un tore d\'efini et d\'eploy\'e sur $k$, on note $X_{*}(A)$ le groupe des homomorphismes alg\'ebriques  de $GL(1)$ dans $A$ et $X^*(A)$ le groupe des homomorphismes alg\'ebriques de $A$ dans $GL(1)$.  En particulier, on  dispose du groupe $X_{*}(A_{G})$. On pose ${\cal A}_{G}=X_{*}(A_{G})\otimes_{{\mathbb Z}}{\mathbb R}$. 
   La d\'efinition de $X^*(A)$ se g\'en\'eralise au cas o\`u $A$ est un sous-groupe ferm\'e pas forc\'ement connexe d'un tore comme ci-dessus. 
   
  \bigskip

   \section{Fonctions sur les groupes finis}

Soient ${\bf G}$ un groupe r\'eductif d\'efini sur un corps fini ${\mathbb F}_{q}$ et ${\bf G}^{\nu}$ un espace tordu sous ${\bf G}$.   On note $C^{inv}({\bf G}^{\nu})$ l'espace des fonctions de ${\bf G}^{\nu}({\mathbb F}_{q})$ dans ${\mathbb C}$ qui sont invariantes par conjugaison par ${\bf G}({\mathbb F}_{q})$. Soit ${\bf M}^{\nu}$ un sous-espace de Levi de ${\bf G}^{\nu}$. On d\'efinit les homomorphismes
$$res_{{\bf M}^{\nu}}^{{\bf G}^{\nu}}: C^{inv}({\bf G}^{\nu})\to C^{inv}({\bf M}^{\nu})$$
 $$ind_{{\bf M}^{\nu}}^{{\bf G}^{\nu}}: C^{inv}({\bf M}^{\nu})\to C^{inv}({\bf G}^{\nu})$$
 de la fa\c{c}on suivante. On fixe un sous-espace parabolique ${\bf P}^{\nu}\in {\cal P}({\bf M}^{\nu})$. Pour $f\in C^{inv}({\bf G}^{\nu})$ et $m\in {\bf M}^{\nu}({\mathbb F}_{q})$, on pose
 $$res_{{\bf M}^{\nu}}^{{\bf G}^{\nu}}(f)(m)=\vert {\bf U}_{{\bf P}}({\mathbb F}_{q})\vert ^{-1}\sum_{u\in {\bf U}_{{\bf P}}({\mathbb F}_{q})} f(mu).$$
 Pour $f\in C^{inv}({\bf M}^{\nu})$, on d\'efinit d'abord une fonction $f[{\bf P}^{\nu}]$ sur ${\bf G}^{\nu}({\mathbb F}_{q})$. Elle est \`a support dans ${\bf P}^{\nu}({\mathbb F}_{q})$. Pour $m\in {\bf M}^{\nu}({\mathbb F}_{q})$ et $u\in {\bf U}_{{\bf P}}({\mathbb F}_{q})$, on a $f[{\bf P}^{\nu}](mu)=f(m)$. 
 Pour $g\in {\bf G}^{\nu}({\mathbb F}_{q})$, on pose
 $$ind_{{\bf M}^{\nu}}^{{\bf G}^{\nu}}(f)(g)=\vert {\bf P}({\mathbb F}_{q})\vert ^{-1}\sum_{x\in {\bf G}({\mathbb F}_{q} )}f[{\bf P}^{\nu}] (x^{-1}gx).$$
 Il est connu que ces d\'efinitions ne d\'ependent pas du choix de ${\bf P}^{\nu}$. On dit qu'une fonction $f\in C^{inv}({\bf G}^{\nu})$ est cuspidale si et seulement si $res_{{\bf M}^{\nu}}^{{\bf G}^{\nu}}(f)=0$ pour tout espace de Levi propre ${\bf M}^{\nu}$. On note $C_{cusp}^{inv}({\bf G}^{\nu})$ le sous-espace des fonctions cuspidales. Il est connu que cet espace est un suppl\'ementaire dans $C^{inv}({\bf G}^{\nu})$ de la somme (non directe) sur tous les espaces de Levi propres ${\bf M}^{\nu}$ des images des homomorphismes $ind_{{\bf M}^{\nu}}^{{\bf G}^{\nu}}$, cf. \cite{Co} proposition 1.1. Cela d\'efinit une projection $C^{inv}({\bf G}^{\nu})\to C^{inv}_{cusp}({\bf G}^{\nu})$ que l'on appelle la projection cuspidale. Plus g\'en\'eralement, soit ${\bf M}^{\nu}$ un espace de Levi. On note $proj_{cusp,{\bf M}^{\nu}}$ la compos\'ee de $res_{{\bf M}^{\nu}}^{{\bf G}^{\nu}}$ et de la projection cuspidale $C^{inv}({\bf M}^{\nu})\to C^{inv}_{cusp}({\bf M}^{\nu})$. Si ${\bf M}_{1}^{\nu}$ et ${\bf M}_{2}^{\nu}$ sont deux espaces de Levi et si $g\in {\bf G}({\mathbb F}_{q})$ conjugue ${\bf M}_{1}^{\nu}$ en ${\bf M}_{2}^{\nu}$, la conjugaison par $g$ envoie naturellement $C^{inv}_{cusp}({\bf M}_{1}^{\nu})$ sur $C^{inv}_{cusp}({\bf M}_{2}^{\nu})$. Il est clair que  la compos\'ee de $proj_{cusp,{\bf M}_{1}^{\nu}}$ et de cet isomorphisme est \'egale \`a $proj_{cusp,{\bf M}_{2}^{\nu}}$. Pour tout espace de Levi ${\bf M}^{\nu}$, notons $n^{{\bf G}}({\bf M}^{\nu})$ le nombre d'\'el\'ements du quotient du groupe $Norm_{{\bf G}({\mathbb F}_{q})}({\bf M}^{\nu})/{\bf M}({\mathbb F}_{q})$, o\`u $Norm_{{\bf G}({\mathbb F}_{q})}({\bf M}^{\nu})$ est le normalisateur de ${\bf M}^{\nu}$ dans ${\bf G}({\mathbb F}_{q})$.  Fixons un ensemble de repr\'esentants $\underline{Levi}({\bf G}^{\nu})$ des classes de conjugaison par ${\bf G}({\mathbb F}_{q})$ de sous-espaces de Levi de ${\bf G}^{\nu}$.  On v\'erifie que
 
 (1) pour tout $f\in C^{inv}({\bf G}^{\nu})$, on a l'\'egalit\'e
 $$f=\sum_{{\bf M}^{\nu}\in \underline{Levi}({\bf G}^{\nu})}n^{{\bf G}}({\bf M}^{\nu})^{-1}ind_{{\bf M}^{\nu}}^{{\bf G}^{\nu}}(proj_{cusp,{\bf M}^{\nu}}(f)).$$ 
 
 Plus g\'en\'eralement, si ${\bf L}^{\nu}$ est un sous-espace de Levi, on a
 
 (2) pour tout $f\in C^{inv}({\bf G}^{\nu})$, on a l'\'egalit\'e
 $$res_{{\bf L}^{\nu}}^{{\bf G}^{\nu}}(f)=\sum_{{\bf M}^{\nu}\in \underline{Levi}({\bf L}^{\nu})}n^{{\bf L}}({\bf M}^{\nu})^{-1}ind_{{\bf M}^{\nu}}^{{\bf L}^{\nu}}(proj_{cusp,{\bf M}^{\nu}}(f)).$$ 
 
 On utilisera la formule suivante, qui est \'equivalente \`a (1). Notons $Par({\bf G}^{\nu})$ l'ensemble des sous-espaces paraboliques de ${\bf G}^{\nu}$ (et non pas un ensemble de repr\'esentants des classes de conjugaison). Pour tout espace parabolique ${\bf P}^{\nu}$, fixons un nombre $z({\bf P}^{\nu})\in {\mathbb R}$ de sorte que les conditions suivantes soient v\'erifi\'ees:
 
 (3) la fonction ${\bf P}^{\nu}\mapsto z({\bf P}^{\nu})$ est constante sur les classes de conjugaison par ${\bf G}({\mathbb F}_{q})$; 
 
 (4) pour tout espace de Levi ${\bf M}^{\nu}$, on a $\sum_{{\bf P}^{\nu}\in {\cal P}({\bf M}^{\nu})}z({\bf P}^{\nu})=1$.
 
 Alors 

 (5) pour tout $f\in C^{inv}({\bf G}^{\nu})$, on a l'\'egalit\'e
 $$f=\sum_{{\bf P}^{\nu}\in Par({\bf G}^{\nu})}z({\bf P}^{\nu})(proj_{cusp,{\bf M}_{{\bf P}}^{\nu}}(f))[{\bf P}^{\nu}],$$
 o\`u, pour chaque ${\bf P}^{\nu}$, on a not\'e ${\bf M}_{{\bf P}}^{\nu}$  une de ses composantes de Levi.  
 
 Preuve.  Consid\'erons la formule (1). Pour chaque ${\bf M}^{\nu}\in \underline{Levi}({\bf G}^{\nu})$, soit ${\bf Q}^{\nu}\in {\cal P}({\bf M}^{\nu})$. Par d\'efinition, on a
 $$ind_{{\bf M}^{\nu}}^{{\bf G}^{\nu}}(proj_{cusp,{\bf M}^{\nu}}(f))=\sum_{{\bf P}^{\nu}}(proj_{cusp,{\bf M}_{{\bf P}}^{\nu}}(f))[{\bf P}^{\nu}],$$
 o\`u on somme sur les ${\bf P}^{\nu}\in Par({\bf G}^{\nu})$ conjugu\'es \`a ${\bf Q}^{\nu}$. En vertu de (4), on peut sommer le membre de droite sur ${\bf Q}^{\nu}\in {\cal P}({\bf M}^{\nu})$, \`a condition de multiplier chaque terme par $z({\bf Q}^{\nu})$. Gr\^ace \`a (3), on obtient
 $$ind_{{\bf M}^{\nu}}^{{\bf G}^{\nu}}(proj_{cusp,{\bf M}^{\nu}}(f))=\sum_{{\bf P}^{\nu}\in Par({\bf G}^{\nu})}z({\bf P}^{\nu})c({\bf P}^{\nu},{\bf M}^{\nu})(proj_{cusp,{\bf M}_{{\bf P}}^{\nu}}(f))[{\bf P}^{\nu}],$$
 o\`u $c({\bf P}^{\nu},{\bf M}^{\nu})$ est le nombre d'\'el\'ements  ${\bf Q}^{\nu}\in {\cal P}({\bf M}^{\nu})$ qui sont conjugu\'es \`a ${\bf P}^{\nu}$. Mais, pour ${\bf P}^{\nu}$ fix\'e, il y a un unique ${\bf M}^{\nu}\in \underline{Levi}({\bf G}^{\nu})$ tel que $c({\bf P}^{\nu},{\bf M}^{\nu})$ soit non nul et, pour celui-ci, on a $c({\bf P}^{\nu},{\bf M}^{\nu})=n^{{\bf G}}({\bf M}^{\nu})$. Alors (5) se d\'eduit de (1). $\square$
 
 \bigskip
 
\section{Le groupe et son immeuble}
   Pour la suite de l'article, on fixe un corps local non-archim\'edien $F$ de caract\'eristique nulle et un groupe alg\'ebrique connexe $G$ d\'efini sur $F$. On note ${\mathbb F}_{q}$ le corps r\'esiduel de $F$, avec la notation usuelle: $q$ est le nombre d'\'el\'ements de ce corps. On note $p$ la caract\'eristique de ${\mathbb F}_{q}$. Introduisons le centre $Z(\hat{G})$ du dual de Langlands $\hat{G}$ de $G$. Ce sont des groupes complexes. Le groupe $Z(\hat{G})$ est muni d'une action du groupe de Galois $\Gamma=Gal(\bar{F}/F)$, o\`u $\bar{F}$ est une cl\^oture alg\'ebrique de $F$. On note $I\subset \Gamma$ le sous-groupe d'inertie et  $Z(\hat{G})^{I}$ le sous-groupe des points fixes par $I$ dans $Z(\hat{G})$.  Le groupe $\Gamma/I$ agit sur  $X^*(Z(\hat{G})^{I})$. On note ${\cal N}= X^*(Z(\hat{G})^{I})^{\Gamma/I}$ le sous-groupe des invariants et $w_{G}:G(F)\to  {\cal N}$ l'homomorphisme d\'efini par Kottwitz. Cet homomorphisme est surjectif, cf. \cite{K} 7.7. Le sous-groupe de torsion ${\cal N}_{tors}$ est fini. 
   
   On note $Imm(G_{AD})$ l'immeuble de $G_{AD}$ sur $F$. Cet ensemble se d\'ecompose en r\'eunion disjointe de facettes, on note $Fac(G)$ l'ensemble de ces facettes. Le groupe $G(F)$ agit sur $Imm(G_{AD})$. Pour ${\cal F}\in Fac(G)$, notons $K_{{\cal F}}^{\dag}$ le stabilisateur de ${\cal F}$ dans $G(F)$.
   
   {\bf Remarque 1.}  Soulignons qu'un \'el\'ement de ce groupe conserve la facette  ${\cal F}$ mais ne fixe pas forc\'ement tout point de cette facette. La litt\'erature immobili\`ere consid\`ere  souvent des fixateurs plut\^ot que des stabilisateurs. Par contre, elle consid\`ere des fixateurs de sous-ensembles compacts de l'immeubles qui ne sont pas forc\'ement des facettes. Notre groupe $K_{{\cal F}}^{\dag}$ est un tel fixateur. En effet, comme on le sait et comme on le rappellera au paragraphe 4, ${\cal F}$ s'identifie naturellement \`a un ouvert d'un espace affine   r\'eel muni d'une distance euclidienne. Alors $K_{{\cal F}}^{\dag}$ est le fixateur du barycentre de ${\cal F}$. 
   
   \bigskip

   L'image $K_{{\cal F}}^{\dag}/Z(G)(F)$ de $K_{{\cal F}}^{\dag}$
    dans $G_{AD}(F)$ est un sous-groupe compact de $G_{AD}(F)$.  Pour tout $\nu\in {\cal N}$, on note $K_{{\cal F}}^{\nu}$ l'ensemble des $x\in K_{{\cal F}}^{\dag}$ tels que $w_{G}(x)=\nu$. On note ${\cal N}({\cal F})$ le sous-groupe des $\nu\in {\cal N}$ tels que $K_{{\cal F}}^{\nu}\not=\emptyset$.  Le plus grand sous-groupe compact de $K_{{\cal F}}^{\dag}$ est la r\'eunion des $K_{{\cal F}}^{\nu}$ pour $\nu\in {\cal N}({\cal F})\cap {\cal N}_{tors}$. L'ensemble $K_{{\cal F}}^0$ est un groupe, c'est le groupe parahorique associ\'e \`a ${\cal F}$.
   
   {\bf Remarque 2.} Cette description des groupes parahoriques est tir\'ee de \cite{L} 2.16 (malheureusement non publi\'e) et de \cite{HR} proposition 3 et remarque 4. Cette derni\`ere r\'ef\'erence s'applique en vertu de la remarque 1 ci-dessus. 
   \bigskip

     Bruhat et Tits ont d\'efini un sous-groupe  $K_{{\cal F}}^+$ de $K_{{\cal F}}^0$ qui est  distingu\'e dans $K_{{\cal F}}^{\dag}$ et pro$-p$-unipotent et un groupe r\'eductif connexe ${\bf G}_{{\cal F  }}$ sur ${\mathbb F}_{q}$ de sorte que    $K_{{\cal F}}^0/K_{{\cal F}}^+$ soit isomorphe \`a ${\bf G}_{{\cal F  }}({\mathbb F}_{q})$. Pour $\nu\in {\cal N}({\cal F})$, il existe un espace tordu ${\bf G}_{{\cal F  }}^{\nu}$ sous ${\bf G}_{{\cal F  }}$ de sorte que $K^{\nu}_{{\cal F}}/K^+_{{\cal F}}$ s'identifie \`a $G^{\nu}_{{\cal F}}({\mathbb F}_{q})$, cf. \cite{L} proposition 2.10.3. Dans le cas o\`u $\nu=0$, on a ${\bf G}_{{\cal F  }}^0={\bf G}_{{\cal F  }}$.

    Le groupe $K^0_{{\cal F}}$ fixe tout point de ${\cal F}$. Pour $\nu\in {\cal N}({\cal F})$, il existe une permutation $\sigma_{{\cal F},\nu}$ de ${\cal F}$ telle que tout \'el\'ement de $K^{\nu}_{{\cal F}}$ agisse sur ${\cal F}$ par cette permutation. Comme on l'a rappel\'e ci-dessus, ${\cal F}$ s'identifie naturellement \`a un ouvert d'un espace affine   r\'eel muni d'une distance euclidienne. La permutation $\sigma_{{\cal F},\nu}$ est une isom\'etrie.  On note ${\cal F}^{\nu}$ le sous-ensemble des \'el\'ements de ${\cal F}$ qui sont fix\'es par $\sigma_{{\cal F},\nu}$. C'est donc l'intersection de ${\cal F}$ avec un sous-espace affine. On note $Fac^*(G)$ l'ensemble des couples $({\cal F},\nu)$ o\`u ${\cal F}\in Fac(G)$ et $\nu\in {\cal N}({\cal F})$. 
    
    On dit qu'un \'el\'ement  de $g\in  G(F)$ est compact modulo $Z(G)$ si son image dans $G_{AD}(F)$ est contenue dans un sous-groupe compact de ce groupe. Cette condition est \'equivalente \`a ce qu'il existe $({\cal F},\nu)\in Fac^*(G)$ tel que $g\in K_{{\cal F}}^{\nu}$.  

  \bigskip

\section{Les facettes d'un appartement}
 
 Fixons un sous-groupe de Levi minimal $M_{min}$ de $G$. Posons $A=A_{M_{min}}$ et ${\cal A}={\cal A}_{M_{Min}}=X_{*}(A)\otimes_{{\mathbb Z}}{\mathbb R}$.  
   Notons $Norm_{G(F)}(A)$  le normalisateur   de $A$ dans $G(F)$. Le quotient $W=Norm_{G(F)}(A)/ M_{min}(F)$ est le groupe de Weyl $W$ de $G$ relatif \`a $A$.    Le groupe $W$ agit sur l'espace ${\cal A}$. On   munit celui-ci d'une norme euclidienne invariante par cette action. On note $\Sigma$   l'ensemble des racines  {\bf r\'eduites} de $A$ dans $G$. A tout $\alpha\in \Sigma$ est associ\'e un sous-groupe radiciel $U_{\alpha}$ de $G$.

 Au tore $A$  est associ\'e un appartement $App(A)\subset Imm(G_{AD})$. C'est une r\'eunion de facettes, on note $Fac(G;A)$ l'ensemble des facettes contenues dans $App(A)$ et $Fac^*(G;A)$ l'ensemble des $({\cal F},\nu)\in Fac^*(G)$ tels que ${\cal F}$ soit contenue dans $App(A)$. 
 L'appartement $App(A)$ s'identifie \`a $ {\cal A}/{\cal A}_{G}$. 

{\bf Remarque 1.} Nous distinguons dans la notation $App(A)$ et ${\cal A}/{\cal A}_{G}$ car $App(A)$ appara\^{\i}t  comme un espace affine sur l'espace vectoriel ${\cal A}/{\cal A}_{G}$. Les actions naturelles sur $App(A)$ sont affines tandis que celles sur ${\cal A}/{\cal A}_{G}$ sont lin\'eaires.  Pour deux \'el\'ements $x,y\in App(A)$,   nous consid\'ererons $x-y$ comme un \'el\'ement de ${\cal A}/{\cal A}_{G}$.

{\bf Remarque 2.} La th\'eorie g\'en\'erale des immeubles fait intervenir toutes les racines de $A$ dans $G$, pas seulement les racines r\'eduites. C'est n\'ecessaire pour d\'ecrire pr\'ecis\'ement les groupes $K_{{\cal F}}^0$. Mais nous n'aurons pas besoin d'un telle pr\'ecision, la description des facettes nous suffit et, pour cela, les racines affines n'interviennent que par  les hyperplans affines sur lesquels elles s'annulent. Si le syst\`eme de racines contient des racines $\alpha$ tels que $2\alpha$ est encore une racine, on peut donc remplacer toute racine affine de la forme $2\alpha+c$ par la fonction affine $\alpha+c/2$:  les hyperplans associ\'es \`a ces fonctions sont les m\^emes. Cela nous permet de ne consid\'erer que les racines r\'eduites, comme nous le faisons. 

\bigskip

L'action sur l'immeuble du sous-groupe  $Norm_{G(F)}(A)$ de $G(F)$ conserve $App(A)$ (inversement, un \'el\'ement de $G(F)$ dont l'action sur l'immeuble conserve $App(A)$ appartient \`a $Norm_{G(F)}(A)$).  Cette action est compatible avec l'action lin\'eaire de $Norm_{G(F)}(A)$ sur ${\cal A}$ via son quotient $W$ dans le sens suivant: soient $n\in Norm_{G(F)}(A)$, $x\in App(A)$ et $e\in {\cal A}/{\cal A}_{G}$; notons $w$ l'image de $n$ dans $W$; alors $n(x+e)=n(x)+w(e)$. 
A tout $\alpha\in \Sigma$ est associ\'e un certain sous-ensemble $\Gamma_{\alpha}$ de ${\mathbb Q}$. C'est l'image r\'eciproque dans ${\mathbb Q}$ d'un sous-ensemble fini de ${\mathbb Q}/{\mathbb Z}$ et on a ${\mathbb Z}\subset \Gamma_{\alpha}$. On a $\Gamma_{-\alpha}=-\Gamma_{\alpha}$. Pour $c\in \Gamma_{\alpha}$, on note $c^+$ le plus petit \'el\'ement de $\Gamma_{\alpha}$ strictement sup\'erieur \`a $c$ et on note $c^-$ le plus grand \'el\'ement de $\Gamma_{\alpha}$ strictement inf\'erieur \`a $c$. On note  $H_{\alpha,c}$ l'hyperplan  affine de  $ App(A)$  d\'efini par l'\'equation $\alpha(x)=c$.   Alors la d\'ecomposition en facettes de $App(A)$ est d\'efinie par la famille d'hyperplans  $(H_{\alpha,c})_{\alpha\in \Sigma, c\in \Gamma_{\alpha}}$.  Ainsi, \`a toute facette ${\cal F}\in Fac(G;A)$ est associ\'e un sous-ensemble  $\Sigma_{{\cal F}}\subset \Sigma(A)$ et, pour tout $\sigma\in \Sigma(A)$, un \'el\'ement $c_{\alpha,{\cal F}}\in  \Gamma_{\alpha}$ de sorte que ${\cal F}$ soit le sous-ensemble des \'el\'ements $x\in App(A)$ v\'erifiant les relations
 
(1)  $\alpha(x)=c_{\alpha,{\cal F}}$ pour tout $\alpha\in \Sigma_{{\cal F}}$;
 
(2)  $c_{\alpha,{\cal F}}<\alpha(x)<c_{\alpha,{\cal F}}^+$ pour tout $\alpha\in \Sigma-\Sigma_{{\cal F}}$.

{\bf Remarque.} Dans le cas o\`u $G$ est quasi-d\'eploy\'e,  $\Gamma_{\alpha}$ est un groupe dont la description est donn\'ee en  \cite{L} 2.5 ou  \cite{BT2} 4.2.21. On en  d\'eduit la description ci-dessus de $\Gamma_{\alpha}$ quand $G$ n'est pas quasi-d\'eploy\'e par descente \'etale, cf. \cite{BT2} 5.1.19.

\bigskip

Introduisons le sous-tore $A_{{\cal F}}\subset A$ tel que $X_{*}(A_{{\cal F}})$ soit le sous-ensemble des \'el\'ements de $X_{*}(A)$ annul\'es par tous les \'el\'ements de $\Sigma_{{\cal F}}$.  
Soit $M_{{\cal F}}$ le commutant de $A_{{\cal F}}$ dans $G$. C'est un sous-groupe de Levi de $G$. Puisque $A_{{\cal F}}$ est contenu dans $A$, $M_{{\cal F}}$ contient $M_{min}$. Notons $\Sigma^{M_{{\cal F}}}\subset \Sigma$ le sous-ensemble des racines  r\'eduites de $A$ dans $M_{{\cal F}}$ autrement dit des  \'el\'ements de $\Sigma$ qui sont triviaux sur $A_{M_{{\cal F}}}$. On a

(3) $A_{{\cal F}}=A_{M_{{\cal F}}}$; $\Sigma_{{\cal F}}\subset \Sigma^{M_{{\cal F}}}$ et ces deux ensembles engendrent le m\^eme ${\mathbb Q}$-sous-espace vectoriel de $X_{*}(A)\otimes_{{\mathbb Z}}{\mathbb Q}$.

Preuve. Pour $\alpha\in \Sigma_{{\cal F}}$, $\alpha$ s'annule sur $X_{*}(A_{{\cal F}})$ donc le groupe radiciel $U_{\alpha}$ commute \`a $A_{{\cal F}}$ et est inclus dans $M_{{\cal F}}$. Cela d\'emontre que $\Sigma_{{\cal F}}\subset \Sigma^{M_{{\cal F}}}$. Par d\'efinition de $M_{{\cal F}}$, $A_{{\cal F}}$ est inclus dans le centre de ce groupe, donc aussi dans son plus grand tore central d\'eploy\'e $A_{M_{{\cal F}}}$. Inversement, $X_{*}(A_{M_{{\cal F}}})$ est le sous-ensemble des \'el\'ements de $X_{*}(A)$ annul\'es par tous les \'el\'ements de $\Sigma^{M_{{\cal F}}}$.  Cet ensemble contenant $\Sigma_{{\cal F}} $, on a $X_{*}(A_{M_{{\cal F}}})\subset X_{*}(A_{{\cal F}})$, d'o\`u $A_{M_{{\cal F}}}\subset A_{{\cal F}}$. Enfin, la derni\`ere assertion r\'esulte de ce que les deux ensembles ont le m\^eme annulateur dans $X_{*}(A)$. $\square$ 

 {\bf Remarque.} Par contre, $\Sigma^{M_{{\cal F}}}$ n'est en g\'en\'eral pas engendr\'e sur ${\mathbb Z}$ par $\Sigma_{{\cal F}}$. 
 
 \bigskip
 
  Posons ${\cal A}_{{\cal F}}= X_{*}(A_{{\cal F}})\otimes_{{\mathbb Z}}{\mathbb R}$, autrement dit, ${\cal A}_{{\cal F}}={\cal A}_{M_{{\cal F}}}$ avec la notation introduite en 1. Il r\'esulte de (1) et (2) que ${\cal A}_{{\cal F}}/{\cal A}_{G}$ est le sous-espace de ${\cal A}/{\cal A}_{G}$ engendr\'e par les $x-y$ pour $x,y\in {\cal F}$. Cela entra\^{\i}ne  que, pour $n\in K_{{\cal F}}^0\cap Norm_{G(F)}(A)$, l'image $w$ de $n$ dans $W$ conserve ${\cal A}_{{\cal F}}$ et y agit trivialement. Autrement dit, $w$ appartient au sous-groupe $W^{M_{{\cal F}}}$ qui est le groupe de Weyl de $M_{{\cal F}}$ relatif \`a $A$. Une racine $\alpha\in \Sigma$ est constante sur ${\cal F}$ si et seulement si elle annule l'espace engendr\'e par les $x-y$ pour $x,y\in {\cal F}$, c'est-\`a-dire ${\cal A}_{{\cal F}}$. Donc 
  $\Sigma^{M_{{\cal F}}}$ est  l'ensemble  des \'el\'ements de $\Sigma$ qui sont constants sur  ${\cal F}$. Pour $\alpha\in \Sigma^{M_{{\cal F}}}$, la valeur constante de $\alpha$ sur ${\cal F}$ appartient \`a $\Gamma_{\alpha}$ si et seulement si $\alpha\in \Sigma_{{\cal F}}$.
La facette ${\cal F}$ est un ouvert du sous-espace affine de $App(A)$ d\'efini par la relation (1), dont l'espace vectoriel associ\'e est ${\cal A}_{{\cal F}}$.  

Les constructions de Bruhat-Tits associent \`a tout $\alpha\in \Sigma$ et \`a tout  $c\in \Gamma_{\alpha}$ un sous-groupe ouvert compact $U_{\alpha,c}$ de $U_{\alpha}(F)$. Les propri\'et\'es suivantes sont bien connues (cf. \cite{L} 2.5 et 2.13, \cite{BT1} 6.4.9 et 7.4.4):
 
 (4)  l'application $c\mapsto U_{\alpha,c}$ est strictement  croissante (si $c<c'$, $U_{\alpha,c}\subsetneq U_{\alpha,c'}$);
 
 (5)  $\cup_{c\in \Gamma_{\alpha}}U_{\alpha,c}=U_{\alpha}(F)$;
 
 (6) quelle que soit ${\cal F}\in Fac(G;A)$, on a $U_{\alpha}(F)\cap K_{{\cal F}}^0=U_{\alpha,c_{\alpha,{\cal F}}}$;  
 
(7)  quelle  que soit ${\cal F}\in Fac(G;A)$, on a $  U_{\alpha}(F)\cap K_{{\cal F}}^+=U_{\alpha,c_{\alpha,{\cal F}}^-}$ si $\alpha\in \Sigma_{{\cal F}}$ et  et $U_{\alpha}(F)\cap K_{{\cal F}}^+=U_{\alpha}(F)\cap K_{{\cal F}}^0=U_{\alpha,c_{\alpha,{\cal F}}}$ si $\alpha\not\in \Sigma_{{\cal F}}$. 

 Soit ${\cal F}\in Fac(G;A)$.  Il est bien connu que l'ensemble des sous-groupes paraboliques  de ${\bf G}_{{\cal F  }}$ est en bijection avec celui des facettes ${\cal F}'$ telles que ${\cal F}$ soit contenue dans l'adh\'erence $\overline{{\cal F}'}$ de ${\cal F}'$. Pr\'ecis\'ement, pour une telle facette ${\cal F}'$, on a $K_{{\cal F}'}^0\subset K_{{\cal F}}^0$ et le parabolique $ {\bf P}_{{\cal F}'}$ associ\'e \`a ${\cal F}'$ est celui pour lequel l'image de $K_{{\cal F}'}^0$ dans ${\bf G}_{{\cal F  }}({\mathbb F}_{q})$ est \'egale \`a ${\bf P}_{{\cal F}'}({\mathbb F}_{q})$.   Il est utile de donner une interpr\'etation plus alg\'ebrique de ces  sous-groupes paraboliques.  Le groupe ${\bf G}_{{\cal F}}$ est la fibre sp\'eciale d'un sch\'ema en groupes sur l'anneau des entiers de $F$ qui contient un mod\`ele entier du tore $A$, cf. \cite{L} 2.11. La fibre sp\'eciale de celui-ci est un sous-tore d\'eploy\'e maximal ${\bf A}$ de ${\bf G}_{{\cal F}}$,   qui v\'erifie 
  les propri\'et\'es suivantes: 

(8) $X_{*}({\bf A}  )\simeq X_{*}(A)$;

(9) soit $n\in Norm_{G(F)}(A)\cap K_{{\cal F}}^{\dag}$; la conjugaison $ad_{n}$ par $n$ conserve $K_{{\cal F}}^0$ et $K_{{\cal F}}^+$ donc se descend en un automorphisme $ad_{n,{\cal F}}$ de ${\bf G}_{{\cal F  }}({\mathbb F}_{q})$; alors $ad_{n,{\cal F}}$ conserve ${\bf A}  $
 et les actions d\'eduites de $ad_{n}$ sur $X_{*}(A)$ et de $ad_{n,{\cal F}}$ sur $X_{*}({\bf A}  )$ co\"{\i}ncident via l'isomorphisme (8).

Soit ${\cal F}'$ une facette de $App(A)$ dont l'adh\'erence contient ${\cal F}$. Via l'isomorphisme (8), le sous-tore $A_{{\cal F}'}$ de $A$ correspond \`a un sous-tore ${\bf A}  _{{\cal F}'}$ de ${\bf A}  $. Il est connu que le commutant ${\bf M}_{{\cal F}'}$ de ${\bf A}_{{\cal F}'}$ dans ${\bf G}_{{\cal F  }}$ est une composante de Levi du sous-groupe parabolique ${\bf P}_{{\cal F}'}$ et que ${\bf A}  _{{\cal F}'}$ est
le plus grand sous-tore d\'eploy\'e dans le centre de ${\bf M}_{{\cal F}'}$. Le sous-groupe ${\bf P}_{{\cal F}'}$ d\'etermine une chambre ouverte dans l'espace ${\cal A}_{{\cal F}'}$ (puisque cet espace s'identifie \`a $X_{*}({\bf A}  _{{\cal F}'})\otimes_{{\mathbb Z}}{\mathbb R}$). En utilisant (5) et (6) on voit que cette chambre est d\'efinie par les relations $\alpha(x)>0$ pour les $\alpha\in \Sigma_{{\cal F}}-\Sigma_{{\cal F}'}$ telles que  $c_{\alpha,{\cal F}}=c_{\alpha,{\cal F}'}$. On voit aussi que ${\bf M}_{{\cal F}'}({\mathbb F}_{q})$ s'identifie \`a ${\bf G}_{{\cal F  }'}({\mathbb F}_{q})$. On montre que cet isomorphisme est alg\'ebrique c'est-\`a-dire qu'il provient d'un isomorphisme ${\bf M}_{{\cal F}'}\simeq {\bf G}_{{\cal F  }'}$. 

Les conditions (1) et (2) sont en g\'en\'eral surabondantes pour d\'ecrire une facette.  Donnons-en une autre description.  D\'ecomposons $G_{AD}$ en produit de composantes simples $G_{1,AD}\times...\times G_{k,AD}$.  L'image $A_{ad}$ de $A$  dans $G_{AD}$ se d\'ecompose conform\'ement en produit  $A_{ad,1}\times...\times A_{ad,k}$. L'ensemble $\Sigma$ se d\'ecompose de m\^eme en union disjointe de composantes irr\'eductibles $\Sigma_{1}\sqcup...\sqcup \Sigma_{k}$. Alors $App(A)$  s'identifie au produit $App(A_{ad,1})\times ..;\times App(A_{ad,k})$ et toute facette ${\cal F}$ se d\'ecompose en produit de facettes ${\cal F}_{1}\times...\times {\cal F}_{k}$. Cela nous ram\`ene au cas o\`u $G_{AD}$ est simple, ce que nous supposons ci-apr\`es. 
 Consid\'erons une facette ouverte ${\cal F}_{min}\in Fac(G;A)$. Alors $\Sigma_{{\cal F}_{min}}=\emptyset$. Il existe  un sous-ensemble $\{\alpha_{0},\alpha_{1},...,\alpha_{n}\}$  de $\Sigma$  tel que:

(10) $\{\alpha_{1},...,\alpha_{n}\}$ est une base de $\Sigma$ (au sens des syst\`emes de racines);

(11)  $\alpha_{0}$ est l'oppos\'e de la racine r\'eduite associ\'ee \`a la plus grande racine relativement \`a cette base; 

(12) ${\cal F}_{min}$ est l'ensemble des $x\in App(A)$ tels que $c_{\alpha_{i},{\cal F}_{min}}<\alpha_{i}(x)$ pour tout $i=0,...,n$. 

 D'apr\`es (10) et (11), il y a
 une relation $\sum_{i=0,...,n}h_{i}\alpha_{i}=0$ avec des entiers $h_{i}>0$ et cette relation est, \`a homoth\'etie pr\`es, la seule relation entre les $\alpha_{i}$. Les sommets de l'adh\'erence de ${\cal F}_{min}$ sont les points $s_{i}$, $i=0,...,n$ d\'efinis par les relations $\alpha_{j}(s_{i})=c_{\alpha_{i},{\cal F}_{min}}$ pour $j\not=i$.  A tout sous-ensemble non vide $I\subset \{0,...,n\}$, associons l'ensemble ${\cal F}_{I}$ des $x\in App(A) $ tels que $c_{\alpha_{i},{\cal F}_{min}}(x)=\alpha_{i}(x)$ pour $i\not\in I$ et $c_{\alpha_{i},{\cal F}_{min}}(x)<\alpha_{i}(x)$ pour $i\in I$. L'application $I\mapsto {\cal F}_{I}$ est une bijection entre l'ensemble des sous-ensembles non vides de $\{0,...,n\}$ et l'ensemble des facettes 
  contenues dans l'adh\'erence de ${\cal F}_{min}$.  Les sommets de l'adh\'erence de la facette ${\cal F}_{I}$ sont les $s_{i}$ pour $i\in I$. La facette est l'int\'erieur relatif de l'enveloppe convexe de ses sommets (int\'erieur relatif au sens o\`u l'on consid\`ere cette enveloppe convexe comme un sous-ensemble du plus petit sous-espace affine de $App(A)$ qui la contient).   
  
   Pour toute facette ${\cal F}\in Fac(G;A)$, il existe $g\in Norm_{G(F)}(A)$ tel que $g({\cal F})$ soit contenue dans l'adh\'erence de ${\cal F}_{min}$.

  \bigskip
   
   \section{Description des ensembles ${\cal F}^{\nu}$} 
  
  Soit $({\cal F},\nu)\in Fac^*(G;A)$. La facette  ${\cal F}$ est d\'efinie par les relations (1)  et (2) du paragraphe pr\'ec\'edent.  On sait  (cf. \cite{BT1} 7.4.4, \cite{L} 2.13) que 
  
  (1) $K_{{\cal F}}^{\nu}$ est engendr\'e par $K_{{\cal F}}^0$ et un \'el\'ement quelconque de $Norm_{G(F)}(A)\cap K_{{\cal F}}^{\nu}$. 
  
  Fixons un \'el\'ement $n$ de cette intersection, notons $w$ son image dans $W$. La permutation $\sigma_{{\cal F},\nu}$ (cf. paragraphe 3)  est la restriction \`a ${\cal F}$ de l'action de $n$ sur $App(A)$. Celle-ci est une action affine et sa composante lin\'eaire est l'action  de $w$ sur ${\cal A}/{\cal A}_{G}$. De plus, l'action de $n$ fixe le barycentre de ${\cal F}$. 
 Cela entra\^{\i}ne que $w$ conserve l'ensemble $\Sigma_{{\cal F}}$, qu'il normalise $W^{M_{{\cal F}}}$ et que  l'action lin\'eaire de $w$ sur ${\cal A}$ conserve ${\cal A}_{{\cal F}}$. Notons ${\cal A}_{{\cal F}}^{\nu}$ le sous-espace des points fixes. On peut aussi bien d\'efinir  un tore $A_{{\cal F}}^{\nu}$ comme la composante neutre de l'ensemble des points fixes par $w$ dans $A_{{\cal F}}$ et alors ${\cal A}_{{\cal F}}^{\nu}=X_{*}(A_{{\cal F}}^{\nu})\otimes_{{\mathbb Z}}{\mathbb R}$.   Comme on le sait, parmi les sous-groupes de Levi $L$ de $G$, contenant $M_{{\cal F}}$ et tels que $w$ appartienne \`a $W^{L}$, il en existe un unique qui soit minimal.  On le note $M_{{\cal F},\nu}$. Il est caract\'eris\'e par l'\'egalit\'e ${\cal A}_{M_{{\cal F},\nu}}={\cal A}_{{\cal F}}^{\nu}$. On voit que ${\cal A}_{{\cal F}}^{\nu}/{\cal A}_{G}$ est  aussi l'espace engendr\'e par les $x-y$ pour $x,y\in {\cal F}^{\nu}$. L'ensemble de racines $\Sigma^{M_{{\cal F},\nu}}$ \'etant celui des \'el\'ements de $\Sigma$ qui annulent $A_{{\cal F}}^{\nu}$, c'est aussi celui des \'el\'ements de $\Sigma$ qui sont constants sur ${\cal F}^{\nu}$. Une cons\'equence de ce qui pr\'ec\`ede, on a l'inclusion
$$(2) \qquad Norm_{G(F)}(A)\cap K_{{\cal F}}^{\nu}\subset Norm_{M_{{\cal F},\nu}(F)}(A).$$

L'ensemble ${\cal F}^{\nu}$ est un polysimplexe. En effet,  on peut supposer que ${\cal F}$ est contenue dans l'adh\'erence d'une facette ${\cal F}_{min}$ comme au paragraphe 4. Puisque la permutation $\sigma_{{\cal F},\nu}$ provient de l'action d'un \'el\'ement de $Norm_{G(F)}(A)$, on se ram\`ene facilement au cas o\`u $G_{AD}$ est simple, ce que l'on suppose ci-apr\`es. Il existe un unique sous-ensemble non vide $I\subset \{0,...,n\}$ de sorte que ${\cal F}={\cal F}_{I}$. L'adh\'erence  $\bar{{\cal F}}$ de la facette ${\cal F}$ est l'enveloppe convexe des points $s_{i}$ pour $i\in I$. Fixons un point $s\in {\cal F}^{\nu}$. Alors ${\cal F}$ est l'ensemble des $x\in App(A)$ pour lesquels il existe une famille $(a_{i})_{i\in I}$ de r\'eels telle que

$a_{i}>0$ pour tout $i\in I$ et $\sum_{i\in I}a_{i}=1$;

$x-s=\sum_{i\in I}a_{i}(s_{i}-s)$.

 L'espace affine engendr\'e par les $s_{i}$ pour $i\in I$ est de dimension $\vert I\vert -1$. On v\'erifie que cela entra\^{\i}ne que la famille $(a_{i})_{i\in I}$  ci-dessus est unique. Puisque la permutation $\sigma_{{\cal F},\nu}$ conserve ${\cal F}$, elle permute les sommets de $\bar{{\cal F}}$  et d\'efinit une permutation de $I$. Notons $I^{\nu}$ l'ensemble des orbites. Un point $x$ \'ecrit comme ci-dessus est fix\'e par $\sigma_{{\cal F},\nu}$ si et seulement si, pour toute orbite $i^{\nu}$, la fonction $i\mapsto a_{i}$ est constante sur $i^{\nu}$. On note $a_{i^{\nu}}$ cette valeur constante multipli\'ee par $\vert i^{\nu}\vert $ et on note $s_{i^{\nu}}$ le point d\'efini par l'\'egalit\'e $s_{i^{\nu}}-s=\vert i^{\nu}\vert ^{-1}\sum_{i\in i^{\nu}}s_{i}-s$. Alors $\sum_{i^{\nu}\in I^{\nu}}a_{i^{\nu}}=1$  et
$$x-s=\sum_{i^{\nu}\in I^{\nu}}a_{i^{\nu}}(s_{i^{\nu}}-s).$$
Autrement dit, $x$ appartient \`a l'int\'erieur relatif de l'enveloppe convexe des points $s_{i^{\nu}}$. La r\'eciproque est imm\'ediate, donc ${\cal F}^{\nu}$ est \'egale \`a l'int\'erieur relatif de cette enveloppe convexe. On v\'erifie que l'espace affine engendr\'e par les $s_{i^{\nu}}$ pour $i^{\nu}\in I^{\nu}$ est de dimension $\vert I^{\nu}\vert -1$. Donc cette enveloppe convexe est un simplexe. Cette description montre de plus que les facettes de l'adh\'erence du simplexe $\bar{{\cal F}}^{\nu}$ sont les ${\cal F}_{J}\cap \bar{{\cal F}}^{\nu}$ pour les ensembles non vides $J\subset I$ qui sont invariants par $\sigma_{{\cal F},\nu}$. 

  On a

(3) pour ${\cal F}_{1}\in Fac(G;A)$, les conditions suivantes sont \'equivalentes:
 
 \noindent (a) $\nu\in {\cal N}({\cal F}_{1})$ et ${\cal F}_{1}\subset \bar{{\cal F}}$;
 
\noindent  (b) ${\cal F}_{1}\cap \bar{{\cal F}}^{\nu}\not=\emptyset$;
 
 \noindent  (c) $\nu\in {\cal N}({\cal F}_{1})$ et ${\cal F}_{1}^{\nu}={\cal F}_{1}\cap \bar{{\cal F}}^{\nu}$;
 
 (4) l'application ${\cal F}_{1}\mapsto {\cal F}_{1}\cap  \bar{{\cal F}}^{\nu}$ est une bijection entre l'ensemble des facettes de $App(A)$ v\'erifiant les conditions de (3) et celui des facettes de  l'adh\'erence du polysimplexe ${\cal F}^{\nu}$. 
 
 Preuve de (3). Modifions un instant nos hypoth\`eses: on ne suppose plus que $\nu\in {\cal F}$ mais on suppose  que ${\cal F}$ est une facette ouverte. Soit ${\cal F}_{1}$  v\'erifiant (a). Fixons $n_{1}\in Norm_{G(F)}(A)\cap K_{{\cal F}_{1}}^{\nu}$.  L'action de $n_{1}$ envoie ${\cal F}$ sur une facette ouverte ${\cal F}'$ dont l'adh\'erence contient encore ${\cal F}_{1}$. Il correspond aux facettes ${\cal F}$ et ${\cal F}'$ des sous-groupes de Borel ${\bf B}_{{\cal F}}$ et ${\bf B}_{{\cal F}'}$ de ${\bf G}_{{\cal F  }_{1}}$. On sait qu'il existe un \'el\'ement de ${\bf G}_{{\cal F  }_{1}}({\mathbb F}_{q})$ qui conjugue ${\bf B}_{{\cal F}'}$ en ${\bf B}_{{\cal F}}$. En relevant cet \'el\'ement en un \'el\'ement $k\in K_{{\cal F}_{1}}^0$, cet \'el\'ement $k$ envoie ${\cal F}'$ sur ${\cal F}$. Alors $g=kn_{1}$ est un \'el\'ement de $K_{{\cal F}_{1}}^{\nu}$ qui conserve ${\cal F}$. Cela entra\^{\i}ne que $\nu\in {\cal N}({\cal F})$ et que  $g\in K_{{\cal F}_{1}}^{\nu}\cap K_{{\cal F}}^{\nu}$. Alors ${\cal F}_{1}^{\nu}$, resp. ${\cal F}^{\nu}$, est l'intersection de ${\cal F}_{1}$, resp. ${\cal F}$, avec l'ensemble des points fixes de $g$. Ce dernier ensemble est convexe. Fixons $y\in {\cal F}^{\nu}$.  Pour $x\in {\cal F}_{1}^{\nu}$, le segment $[x,y]$ joignant $x$ \`a $y$ est contenu dans cet ensemble de points fixes. Mais le segment $]x,y]$ ouvert en $x$  est aussi contenu dans ${\cal F}$. Donc $]x,y]$ est inclus dans ${\cal F}^{\nu}$ et $x$ appartient \`a l'adh\'erence de cet ensemble. Cela prouve (c). 
 
 {\bf Remarque.} On vient de prouver que, pour  toute facette ouverte ${\cal F}$ dont l'adh\'erence contient une facette ${\cal F}_{1}$ telle que $\nu\in {\cal N}({\cal F}_{1})$, on a $\nu\in {\cal N}({\cal F})$. Puisque toutes les facettes ouvertes sont conjugu\'ees par l'action de $G$, on obtient
 
 (5) soient deux facettes ${\cal F}$ et ${\cal F}_{1}$; supposons ${\cal F}$ ouverte; alors ${\cal N}({\cal F}_{1})\subset {\cal N}({\cal F})$.
 
 \bigskip
 
 Revenons maintenant \`a nos hypoth\`eses: $n\in {\cal N}({\cal F})$ mais ${\cal F}$ n'est plus suppos\'ee ouverte. Soit ${\cal N}_{1}$ v\'erifiant (a). On fixe une facette ouverte ${\cal F}_{min}$ dont l'adh\'erence contient ${\cal F}$, donc aussi ${\cal F}_{1}$. D'apr\`es ce que l'on vient de prouver, on a $\nu\in {\cal N}({\cal F}_{min})$ et 
${\cal F}_{1}^{\nu}={\cal F}_{1}\cap \bar{{\cal F}}^{\nu}_{min}$.  On peut appliquer le m\^eme r\'esultat pour ${\cal F}$:  ${\cal F}^{\nu}={\cal F}\cap \bar{{\cal F}}^{\nu}_{min}$. Fixons $y\in {\cal F}^{\nu}$ et soit $x\in {\cal F}_{1}^{\nu}$. De nouveau, le segment $[x,y]$ est inclus dans $\bar{{\cal F}}^{\nu}_{min}$. Le segment $]x,y]$ est inclus dans ${\cal F}$, donc dans ${\cal F}^{\nu}$. Mais alors $x$ appartient \`a l'adh\'erence de ${\cal F}^{\nu}$, ce qui d\'emontre (c). Il est clair que (c) implique (b). Enfin, si (b) est v\'erifi\'ee, l'\'el\'ement $n$ introduit plus haut fixe tout point de $ \bar{{\cal F}}^{\nu}$ donc fixe un point de ${\cal F}_{1}$, donc conserve cette facette. Puisque $w_{G}(n)=\nu$, cela entra\^{\i}ne $\nu\in {\cal N}({\cal F}_{1})$. Enfin, ${\cal F}_{1}$ coupe $\bar{{\cal F}}$ et est donc incluse dans cette adh\'erence. Cela prouve (a) et ach\`eve la d\'emonstration de (3). 

Preuve de (4). D'apr\`es la description de $\bar{{\cal F}}^{\nu}$ donn\'ee plus haut, il existe un sous-ensemble fini ${\cal X}$ de $Fac(G;A)$ tel que l'application ${\cal F}_{1}\mapsto {\cal F}_{1}\cap \bar{{\cal F}}^{\nu}$ soit une bijection de ${\cal X}$ sur l'ensemble des facettes du polysimplexe $\bar{{\cal F}}^{\nu}$. Puisque deux facettes distinctes sont disjointes, ${\cal X}$ ne peut \^etre que l'ensemble des ${\cal F}_{1}\in Fac(G;A)$ telles que ${\cal F}_{1}\cap \bar{{\cal F}}^{\nu}$ ne soit pas vide. C'est-\`a-dire celui des facettes v\'erifiant le (b) de (3). Cela d\'emontre (4). $\square$

Montrons que

(6)  l'ensemble des sous-espaces paraboliques de ${\bf G}_{{\cal F  }}^{\nu}$ est en bijection avec celui des facettes ${\cal F}'$ telles que ${\cal F}$ soit contenue dans $\bar{{\cal F}}'$ et que $\nu$ appartienne \`a ${\cal N}({\cal F}')$; pour une telle facette ${\cal F}'$,   on a $K_{{\cal F}'}^{\nu}\subset K_{{\cal F}}^{\nu}$ et l'espace parabolique ${\bf P}_{{\cal F  }'}^{\nu}$ associ\'e \`a ${\cal F}'$ est tel que ${\bf P}_{{\cal F  }'}^{\nu}({\mathbb F}_{q})$ soit l'image de $K_{{\cal F}'}^{\nu}$ dans ${\bf G}_{{\cal F  }}^{\nu}({\mathbb F}_{q})$. 

Preuve. Soit ${\bf P}^{\nu}$ un sous-espace parabolique de ${\bf G}_{{\cal F  }}^{\nu}$.  Il lui est associ\'e un sous-groupe parabolique ${\bf P}$ de ${\bf G}_{{\cal F  }}$, donc aussi   une facette ${\cal F}'$ telle que ${\cal F}\subset \bar{{\cal F}}'$ et ${\bf P}={\bf P}_{{\cal F  }'}$.    Soit $g\in K_{{\cal F}}^{\nu}$ dont l'image ${\bf g}$ dans ${\bf G}_{{\cal F  }}^{\nu}({\mathbb F}_{q})$ appartienne \`a  ${\bf P}^{\nu}({\mathbb F}_{q})$. L'action de $g$ sur l'immeuble envoie ${\cal F}'$ sur une  facette ${\cal F}''$ dont l'adh\'erence contient encore ${\cal F}$. Il lui est associ\'e un sous-groupe parabolique $ {\bf P}_{{\cal F  }''}$ de ${\bf G}_{{\cal F  }}$. Alors l'action de ${\bf g}$ par  conjugaison sur ${\bf G}_{{\cal F  }}$  envoie ${\bf P}_{{\cal F  }'}$ sur ${\bf P}_{{\cal F  }''}$. Puisque ${\bf g}$ appartient \`a  ${\bf P}^{\nu}({\mathbb F}_{q})$, l'action de ${\bf g}$ conserve ${\bf P}_{{\cal F  }'}$.  Cela entra\^{\i}ne ${\bf P}_{{\cal F  }''}={\bf P}_{{\cal F  }'}$, d'o\`u ${\cal F}''={\cal F}'$. Donc  $g$ appartient \`a  $K_{{\cal F}'}^{\dag}$. Puisque $w_{G}(g)=\nu$, cela entra\^{\i}ne $g\in K_{{\cal F}'}^{\nu}$. A fortiori, cet ensemble n'est pas vide, donc $\nu\in {\cal N}({\cal F}')$.

R\'eciproquement, soit ${\cal F}'$ une facette de l'immeuble telle que ${\cal F}$ soit contenue dans $\bar{{\cal F}}'$ et que $\nu$ appartienne \`a ${\cal N}({\cal F}')$. Puisque deux facettes sont toujours contenues dans un appartement, on ne perd rien \`a supposer que ${\cal F}'$ est contenue dans $App(A)$. D'apr\`es (3), ${\cal F}^{\nu}$ est contenu dans $\bar{{\cal F}}^{_{'}\nu}$. Soit $n'$ un  \'el\'ement de $Norm_{G(F)}(A)\cap K_{{\cal F}'}^{\nu}$. Alors l'action de $n'$ fixe tout point  de $\bar{{\cal F}}^{_{'}\nu}$, donc fixe tout point de ${\cal F}^{\nu}$, donc conserve ${\cal F}$, donc appartient \`a $K_{{\cal F}}^{\nu}$. Notons ${\bf n}'$ 'image  de $n'$ dans ${\bf G}_{{\cal F  }}^{\nu}({\mathbb F}_{q})$. L'action de ${\bf n}'$ dans ${\bf G}_{{\cal F  }}$ conserve le sous-groupe parabolique ${\bf P}_{{\cal F  }'}$. Donc ${\bf P}_{{\cal F  }'}^{\nu}={\bf n}'{\bf P}_{{\cal F  }'}$ est un sous-espace parabolique de ${\bf G}_{{\cal F  }}^{\nu}$. Puisqu'on sait d'apr\`es (1) que  $K_{{\cal F}'}^{\nu}=n'K_{{\cal F}'}^0$, les derni\`eres assertions r\'esultent de la d\'efinition ci-dessus de ${\bf P}_{{\cal F}'}^{\nu}$ et du fait que ${\bf P}_{{\cal F}'}({\mathbb F}_{q})$ est l'image naturelle de $K_{{\cal F}'}^0$ dans ${\bf G}_{{\cal F  }}({\mathbb F}_{q})$. $\square$

Introduisons le sous-tore d\'eploy\'e maximal ${\bf A}  $ de ${\bf G}_{{\cal F  }}$ v\'erifiant les conditions (8) et (9) du paragraphe 4. Soit ${\cal F}'$ une facette de $App(A)$ dont l'adh\'erence contient ${\cal F}$ et telle que $\nu\in {\cal N}({\cal F}')$. Via l'isomorphisme $X_{*}({\bf A}  )\simeq X_{*}(A)$, le sous-tore $A_{{\cal F}'}$ de $A$ correspond \`a un sous-tore ${\bf A}  _{{\cal F}'}$ de ${\bf A}  $. On a d\'ej\`a dit que le commutant ${\bf M}_{{\cal F}'}$ de ${\bf A}_{{\cal F}'}$ dans ${\bf G}_{{\cal F  }}$ \'etait une composante de Levi du sous-groupe parabolique ${\bf P}_{{\cal F}'}$ et que ${\bf A}  _{{\cal F}'}$ \'etait
le plus grand sous-tore d\'eploy\'e dans le centre de ${\bf M}_{{\cal F}'}$.   Le tore $A_{{\cal F}}^{\nu}$ correspond de m\^eme \`a un sous-tore ${\bf A}  _{{\cal F}}^{\nu}$ de ${\bf A}  $.  En appliquant la condition 4(9)  \`a un \'el\'ement de $Norm_{G(F)}(A)\cap K_{{\cal F}'}^{\nu}$, on voit que ${\bf A}  _{{\cal F}}^{\nu}$ est le plus grand sous-tore d\'eploy\'e contenu dans le commutant de ${\bf M}_{{\cal F}'}^{\nu}$ dans ${\bf M}_{{\cal F}'}$. La th\'eorie g\'en\'erale des  espaces tordus nous dit que ${\bf M}_{{\cal F}'}$ est encore le commutant de ${\bf A}_{{\cal F}'}^{\nu}$ dans ${\bf G}_{{\cal F  }}$. L'espace ${\bf M}_{{\cal F}'}^{\nu}({\mathbb F}_{q})$ s'identifie \`a ${\bf G}_{{\cal F  }'}^{\nu}({\mathbb F}_{q})$. De nouveau, cette identification provient d'un isomorphisme alg\'ebrique ${\bf M}_{{\cal F}'}^{\nu}\simeq {\bf G}_{{\cal F  }'}^{\nu}$.

\bigskip

 \section{Descente aux sous-groupes de Levi}

  Soit $M$ un sous-groupe de Levi de $G$. Notons $M_{ad}=M/Z(G)$ l'image de $M$ dans $G_{AD}$. L'immeuble $Imm(M_{ad})$ de ce groupe s'identifie \`a un sous-ensemble de $Imm(G_{AD})$, \`a savoir la r\'eunion des appartements $App(A')$ pour les sous-tores d\'eploy\'es maximaux $A'$ de $M$. L'action de $M(F)$ sur $Imm(M_{ad})$ est la restriction de celle sur $Imm(G_{AD})$. Par contre, la d\'ecomposition en facettes de $Imm(M_{ad})$ n'est pas  celle de ce sous-ensemble, elle est moins fine. L'espace ${\cal A}_{M}/{\cal A}_{G}$ agit sur $Imm(M_{ad})$ car ${\cal A}_{M}$ est inclus dans ${\cal A}'=X_{*}(A')\otimes_{{\mathbb Z}}{\mathbb R}$ pour tout tore $A'$ comme ci-dessus. Le quotient de $Imm(M_{ad})$ par l'action de cet espace s'identifie \`a l'immeuble $Imm(M_{AD})$ du groupe adjoint de $M$. On note $p_{M}:Imm(M_{ad})\to Imm(M_{AD})$ cette projection. Elle est \'equivariante pour les actions de $M(F)$. Pour tout tore $A'$ comme ci-dessus, on note $App^M(A')=App(A')/({\cal A}_{M}/{\cal A}_{G})$ l'appartement de $Imm(M_{AD})$ associ\'e \`a $A'$ et on note  encore $p_{M}:App(A')\to App^M(A')$ la projection. Elle est \'equivariante pour les actions de $Norm_{M(F)}(A')$. 
  
  De m\^eme que l'on a introduit le groupe ${\cal N}={\cal N}^G$, cf. paragraphe 3, on introduit le groupe ${\cal N}^M$.  Dire qu'un  \'el\'ement $m\in M(F)$ est compact modulo $Z(G)$, c'est-\`a-dire que  l'image dans de $m$ dans $G_{AD}(F)$  est contenue dans un sous-groupe compact de $G_{AD}(F)$, revient  \`a dire que l'image $m_{ad}$ de  $m$ dans $M_{ad}(F)$ est contenue dans un sous-groupe compact de $M_{ad}(F)$. Pour un tel \'el\'ement, $w_{M_{ad}}(m_{ad})$ est de torsion. On a un homomorphisme naturel ${\cal N}^M\to {\cal N}^{M_{ad}}$. On note ${\cal N}^{M}_{G-comp}$ l'image r\'eciproque par cet homomorphisme du sous-groupe de torsion de ${\cal N}^{M_{ad}}$ (remarquons que, si $M=G$, ${\cal N}^{G}_{G-comp}={\cal N}^G$).  De l'inclusion $Z(\hat{G})\subset Z(\hat{M})$ se d\'eduit un homomorphisme ${\cal N}^M\to {\cal N}^G$. Montrons que

 (1) cet homomorphisme se restreint en un homomorphisme injectif  ${\cal N}^{M}_{G-comp}\to {\cal N}^G$.  
 
On introduit le rev\^etement simplement connexe $\hat{G}_{SC}$ du groupe d\'eriv\'e de $\hat{G}$ et l'image r\'eciproque $\hat{M}_{sc}$ de $\hat{M}$ dans $\hat{G}_{SC}$. Le groupe $\hat{G}_{SC}$, resp. $\hat{M}_{sc}$, est le groupe dual de $G_{AD}$, resp. $M_{ad}$.  Des homomorphismes naturels 
$$\begin{array}{ccc}Z(\hat{G})^{I}&&\\ &\searrow&\\ && Z(\hat{M})^{I}\\ &\nearrow&\\ Z(\hat{M}_{sc})^{I}&&\\ \end{array}$$
se d\'eduisent des homomorphismes
$$\begin{array}{ccc}&&X^*(Z(\hat{G})^{I})\\ &\nearrow&\\ X ^*(Z(\hat{M})^{I})&&\\ &\searrow&\\ &&X^*(Z(\hat{M}_{sc})^{I})\\ \end{array}$$
Les groupes ${\cal N}^M$, ${\cal N}^G$ et ${\cal N}^{M_{ad}}$ sont les sous-groupes de points fixes par $\Gamma/I$ dans les groupes ci-dessus. Notons $Y^*$ le sous-groupe des \'el\'ements de $X ^*(Z(\hat{M})^{I})$ qui s'envoient dans le sous-groupe de torsion de $X^*(Z(\hat{M}_{sc})^{I})$. Alors ${\cal N}^{M}_{G-comp}$ est le sous-groupe des points fixes par $\Gamma/I$ dans $Y^*$. Il suffit de d\'emontrer que $Y^*$ s'envoie injectivement dans $X^*(Z(\hat{G})^{I})$. On introduit le groupe $\hat{M}_{ad}=\hat{M}/Z(\hat{G})$.
  On a une suite exacte
 $$1\to Z(\hat{G})\to Z(\hat{M})\to Z(\hat{M}_{ad})\to 1$$
 Parce que $X_{*}(Z(\hat{M}_{ad}))$ est un module galoisien induit, le groupe $Z(\hat{M}_{ad})^{I}$ est connexe. Donc la suite 
 $$1\to Z(\hat{G})^{I}\to Z(\hat{M})^{I}\to Z(\hat{M}_{ad})^{I}\to 1$$
 est encore exacte. D'o\`u une suite exacte
 $$0\to X^*(Z(\hat{M}_{ad})^{I})\to X^*(Z(\hat{M})^{I})\to X^*(Z(\hat{G})^{I})\to 0$$
 On doit donc montrer que $Y^*\cap X^*(Z(\hat{M}_{ad})^{I}) =\{0\}$. On a l'\'egalit\'e $\hat{M}_{ad}=\hat{M}_{sc}/Z(\hat{G}_{SC})$. On peut donc remplacer ci-dessus $\hat{M}$ et $\hat{G}$ par $\hat{M}_{sc}$ et $\hat{G}_{SC}$ et on  a  une suite exacte
  $$0\to X^*(Z(\hat{M}_{ad})^{I})\to X^*(Z(\hat{M}_{sc})^{I})\to X^*(Z(\hat{G}_{SC})^{I})\to 0$$
  Les premiers homomorphismes des deux suites ci-dessus sont injectifs et le diagramme suivant est commutatif:
  $$\begin{array}{ccc}&&X^*(Z(\hat{M})^{I})\\ &\nearrow&\\ X ^*(Z(\hat{M}_{ad})^{I})&& \downarrow\\ &\searrow&\\ &&X^*(Z(\hat{M}_{sc})^{I})\\ \end{array}$$
  Un \'el\'ement  $y\in Y^*\cap X^*(Z(\hat{M}_{ad})^{I}) $ est un \'el\'ement de $X^*(Z(\hat{M}_{ad})^{I})$ qui s'envoie sur un \'el\'ement de torsion dans $X^*(Z(\hat{M}_{sc})^{I})$. D'apr\`es l'injectivit\'e de la fl\`eche sud-est ci-dessus, $y$ est un \'el\'ement de torsion dans $X^*(Z(\hat{M}_{ad})^{I})$. Or $Z(\hat{M}_{ad})^{I}$ est connexe donc $X^*(Z(\hat{M}_{ad})^{I})$ est sans torsion. Donc $y=0$, ce qui prouve (1). $\square$

 Gr\^ace \`a (1), on identifie   ${\cal N}^{M}_{G-comp}$ \`a un sous-groupe de $ {\cal N}={\cal N}^G$. A toute facette ${\cal F}_{M}$ de $Imm(M_{AD})$ est associ\'e comme au paragraphe 3 un sous-groupe ${\cal N}^M({\cal F}_{M})$ de ${\cal N}^M$. On pose ${\cal N}^{M}_{G-comp}({\cal F}_{M})={\cal N}^M({\cal F}_{M})\cap {\cal N}^{M}_{G-comp}$. L'ensemble $\bigcup_{\nu\in {\cal N}^{M}_{G-comp}({\cal F}_{M})}K_{{\cal F}_{M}}^{\nu}$ est un groupe dont l'image dans $M_{ad}(F)$ est   compacte.
 On note $Fac^*(M)_{G-comp}$ l'ensemble des couples $({\cal F}_{M},\nu_{M})\in Fac^*(M)$ tels que $\nu_{M}\in {\cal N}^M_{G-comp}$. On d\'efinit de m\^eme $Fac^*(M;A)_{G-comp}$. Un \'el\'ement $m\in M(F)$ est compact modulo $Z(G)$ si et seulement s'il existe $({\cal F}_{M},\nu)\in Fac^*(M)_{G-comp}$ tel que $m\in K_{{\cal F}_{M}}^{\nu}$.

 On suppose pour la suite de cette section que $M$ contient $A$. Soit $({\cal F},\nu)\in Fac^*(G;A)$. 
 L'image de ${\cal F}$ par la projection $p_{M}$ est contenue dans une unique facette ${\cal F}^M$ de$App^M(A)$:  ${\cal F}$ \'etant l'ensemble des $x\in App(A)$ v\'erifiant les relations (1) et (2) du paragraphe 4,
  ${\cal F}^M$ est  l'ensemble des $x\in App^M(A)$ v\'erifiant 
  
  (2) $\alpha(x)=c_{\alpha,{\cal F}}$ pour $\alpha\in \Sigma_{{\cal F}}\cap \Sigma^M$;
  
  (3)  $c_{\alpha,{\cal F}}<\alpha(x)<c_{\alpha,{\cal F}}^+$ pour tout $\alpha\in \Sigma^M-\Sigma_{{\cal F}}\cap \Sigma^M$.
  
  On a d\'efini le groupe $M_{{\cal F},\nu}$ au paragraphe 5. Remarquons que, si ce groupe est contenu dans $M$, on a $\Sigma_{{\cal F}}\subset \Sigma^M$, ce qui simplifie les relations ci-dessus.

\ass{Lemme }{Soient $({\cal F},\nu)\in Fac^*(G;A)$ et $M$ un sous-groupe de Levi de $G$ contenant $M_{{\cal F},\nu}$. Alors

(i) L'\'el\'ement $\nu$ appartient \`a ${\cal N}^{M}_{G-comp}({\cal F}^M)$. 

(ii) On a les \'egalit\'es $K_{{\cal F}^M}^+= K_{{\cal F}}^+\cap M(F)$, $K_{{\cal F}^M}^0= K_{{\cal F}}^0\cap M(F)$, $K_{{\cal F}^M}^{\nu}=K_{{\cal F}}^{\nu}\cap M(F)$.

(iii) Soit $P\in {\cal P}(M)$, notons $\bar{P}$ le parabolique oppos\'e. Alors le groupe $K_{{\cal F}}^+$ est produit de ses sous-groupes $K_{{\cal F}}^+\cap U_{P}(F)$,   $K_{{\cal F}^M}^+$,   et $K_{{\cal F}}^+\cap U_{\bar{P}}(F)$. Le groupe $K_{{\cal F}}^0$ est produit de ses sous-groupes $K_{{\cal F}}^+\cap U_{P}(F)$,   $K_{{\cal F}^M}^0$,   et $K_{{\cal F}}^+\cap U_{\bar{P}}(F)$. L'ensemble $K_{{\cal F}}^{\nu}$ est produit du groupe $K_{{\cal F}}^+\cap U_{P}(F)$,    de l'ensemble $K_{{\cal F}^M}^{\nu}$ et du groupe $K_{{\cal F}}^+\cap U_{\bar{P}}(F)$. Ces produits poss\`edent les propri\'et\'es usuelles des d\'ecompositions "de type Iwahori".

 (iv) ${\cal A}_{{\cal F}^M,\nu}={\cal A}_{{\cal F},\nu}$ et $M_{{\cal F}^M,\nu}=M_{{\cal F},\nu}$;
 
 (v) ${\cal F}$ est un sous-ensemble ouvert de $p_{M}^{-1}({\cal F}^M)$;
 
 (vi) ${\cal F}^{\nu}={\cal F}\cap p_{M}^{-1}({\cal F}^{M,\nu})$ et cet ensemble est ouvert dans $p_{M}^{-1}({\cal F}^{M,\nu})$.}

 Preuve.  Bruhat et Tits on d\'efini un   sous-groupe distingu\'e $M_{min}(F)^+$ de $M_{min}(F)$ qui est pro-$p$-unipotent et tel que le groupe $K_{{\cal F}}^+$ soit le produit  de $ M_{min}(F)^+$    et des groupes radiciels $K_{{\cal F}}^+\cap U_{\alpha}(F)$. 
  L'\'egalit\'e $K_{{\cal F}^M}^+= K_{{\cal F}}^+\cap M(F)$ se d\'eduit alors de la description des facettes ${\cal F}$ et ${\cal F}^M$ et de 4(7). On obtient de m\^eme 
  la d\'ecomposition en produit de (iii) du groupe $K_{{\cal F}}^+$. 
 
Notons $M_{min}(F)^0$ le sous-groupe des $m\in M_{min}(F)$ tels que $w_{M_{min}}(m)=0$.  Le groupe $K_{{\cal F}}^0$ est engendr\'e par $M_{min}(F)^0$   et par les sous-groupes radiciels $K_{{\cal F}}^0\cap U_{\alpha}(F)$  (il ne s'agit plus d'une d\'ecomposition en produit).    Il s'ensuit encore l'\'egalit\'e $K_{{\cal F}^M}^0= K_{{\cal F}}^0\cap M(F)$ et la d\'ecomposition en produit de (iii) du groupe $K_{{\cal F}}^0$. 
 
  On a vu en 5(1)  que $K_{{\cal F}}^{\nu}$ etait engendr\'e par $K_{{\cal F}}^0$ et  par un \'el\'ement quelconque de $K_{{\cal F}}^{\nu}\cap  Norm_{G(F)}(A)$. Fixons un tel \'el\'ement $n$. D'apr\`es 5(2), $n$ appartient \`a $M_{{\cal F},\nu}(F)$, a fortiori \`a $M(F)$, donc \`a $Norm_{M(F)}(A)$. L'action de $n$ sur $Imm(G_{AD})$  conserve ${\cal F}$. D'apr\`es l'\'equivariance de $p_{M}$ pour les actions de $Norm_{M(F)}(A)$, l'action de $n$  sur $Imm(M_{AD})$ conserve ${\cal F}^M$. Donc $n\in K_{{\cal F}^M}^{\dag}$.   Puisque $n$ appartient \`a $K_{{\cal F}}^{\dag}$, son image dans $G_{AD}(F)$ est contenue dans un sous-groupe compact, donc son image dans $M_{ad}(F)$ v\'erifie la m\^eme propri\'et\'e. Donc $w_{M}(n)$ appartient \`a ${\cal N}^{M}_{G-comp}$. Puisque l'image de $w_{M}(n)$ dans ${\cal N}$ est $\nu$, cela entra\^{\i}ne que $\nu$ appartient \`a ${\cal N}^{M}_{G-comp}$ puis que $n\in K_{{\cal F}^M}^{\nu}$. Ces deux derni\`eres propri\'et\'es entra\^{\i}nent     le (i) de l'\'enonc\'e. Puisque $n$ \'etait un \'el\'ement quelconque de $K_{{\cal F}}^{\nu}\cap  Norm_{G(F)}(A)$, on vient de prouver que
   $K_{{\cal F}}^{\nu}\cap  Norm_{G(F)}(A) \subset K_{{\cal F}^M}^{\nu}\cap Norm_{M(F)}(A)$. Puisque $K_{{\cal F}}^{\nu}$ est engendr\'e par $K_{{\cal F}}^0$ et  par un \'el\'ement quelconque du premier ensemble et que, de m\^eme,  $K_{{\cal F}^M}^{\nu}$ est engendr\'e par $K_{{\cal F}^M}^0$ et  par un \'el\'ement quelconque du second, les assertions de (ii) et (iii) relatives \`a $K_{{\cal F}}^{\nu}$ se d\'eduisent de celles d\'ej\`a vues relatives \`a $K_{{\cal F}}^{0}$.
  
   La facette ${\cal F}$ est l'ensemble des $x\in App(A)$ v\'erifiant   les relations (1) et (2) du paragraphe 4 et $ {\cal F}^M$ est l'ensemble des $x\in App^M(A)$ v\'erifiant les relations (2) et (3) ci-dessus (et on a $\Sigma_{{\cal F}}\subset \Sigma^M$). Par d\'efinition, $p_{M}^{-1}({\cal F}^M)$ est 
   l'ensemble des $x\in App(A)$ v\'erifiant 
   ces relations (2) et (3). 
   Donc ${\cal F}$ est le sous-ensemble de $p_{M}^{-1}({\cal F}^M)$ d\'efini par les relations 
   
       $c_{\alpha,{\cal F}}<\alpha(x)<c_{\alpha,{\cal F}}^+$ pour tout $\alpha\in \Sigma-\Sigma^M$.
 
 Ces relations sont ouvertes. Donc ${\cal F}$ est un sous-ensemble ouvert de $p_{M}^{-1}({\cal F}^M)$. Une racine $\alpha\in \Sigma^M$ qui est constante sur ${\cal F}^M$ est aussi constante sur $p_{M}^{-1}({\cal F}^M)$ donc aussi sur ${\cal F}$. Inversement, soit $\alpha\in \Sigma$ qui est constante sur ${\cal F}$. Puisque $M$ contient $M_{{\cal F},\nu}$, a fortiori contient $M_{{\cal F}}$, on a $\alpha\in \Sigma^M$. La valeur de $\alpha$ en un point de $ App(A)$ est aussi sa valeur sur la projection de ce point dans $App^M(A)$. Donc $\alpha$ est constante sur la projection de ${\cal F}$. Or celle-ci est un ouvert non vide de ${\cal F}^M$. Donc $\alpha$ est constante sur ${\cal F}^M$.  Cela d\'emontre que l'ensemble des $\alpha\in \Sigma^M$ qui sont constantes sur ${\cal F}^M$ est \'egal \`a celui des $\alpha\in \Sigma$ qui sont constantes sur ${\cal F}$. Cela entra\^{\i}ne ${\cal A}_{{\cal F}}={\cal A}_{{\cal F}^M}$ et $M_{{\cal F}}=M_{{\cal F}^M}$.  Fixons $n\in Norm_{M(F)}(A)\cap K_{{\cal F}^M}^{\nu}$. D'apr\`es le  (iii) d\'ej\`a d\'emontr\'e, ${\cal F}^{\nu}$ est l'ensemble des points fixes dans ${\cal F}$ de l'action de $n$ tandis que ${\cal F}^{M,\nu}$ est celui des points fixes dans ${\cal F}^{M}$ de l'action du m\^eme \'el\'ement. Puisque la projection est \'equivariante pour l'action de $n$, la projection de ${\cal F}^{\nu}$
est incluse dans ${\cal F}^{M,\nu}$, c'est-\`a-dire ${\cal F}^{\nu}\subset p_{M}^{-1}({\cal F}^{M,\nu})$.  
 Fixons $x\in {\cal F}^{\nu}$ et soit $y\in p_{M}^{-1}({\cal F}^{M,\nu})\cap {\cal F}$. Parce que $y$ se projette en un point de ${\cal F}^{M,\nu}$, on a $n(y)-y\in {\cal A}_{M}/{\cal A}_{G}$. Soit $e\in {\cal A}_{M}/{\cal A}_{G}$ tel que $n(y)-y=e$. Posons $f=y-x$ et notons $w$ l'image de $n$ dans $W^M$. Puisque $n(x)=x$, on a $w(f)=f+e$. Mais $e$ appartient \`a l'espace des points fixes de l'action de $w$ dans ${\cal A}/{\cal A}_{G}$ donc $w^k(f)=f+ke$ pour tout $k\in {\mathbb N}$.  Puisque $w$ est d'ordre fini, cela entra\^{\i}ne $e=0$. Donc $n(y)=y$ et $y$ appartient \`a ${\cal F}^{\nu}$. Cela prouve l'\'egalit\'e du (vi) de l'\'enonc\'e. Puisque ${\cal F}$ est ouvert dans $p_{M}^{-1}({\cal F}^M)$, son intersection avec le sous-ensemble $p_{M}^{-1}({\cal F}^{M,\nu})$ est aussi ouverte dans ce sous-ensemble. Mais alors, les m\^emes arguments que ci-dessus montrent que l'ensemble des $\alpha\in \Sigma^M$ qui sont constantes sur ${\cal F}^{M,\nu}$ est \'egal \`a celui des $\alpha\in \Sigma$ qui sont constantes sur ${\cal F}^{\nu}$. D'o\`u le (iv) de l'\'enonc\'e. $\square$

 Une cons\'equence du lemme est que les deux espaces $K_{{\cal F}}^{\nu}/K_{{\cal F}}^{+}={\bf G}_{{\cal F  }}^{\nu}({\mathbb F}_{q})$ et $K_{{\cal F}^M}^{\nu}/K_{{\cal F}^M}^{+}={\bf M}_{{\cal F}^{M}}^{\nu}({\mathbb F}_{q})$ sont naturellement isomorphes. De nouveau, on montre que cet isomorphisme ainsi que son analogue pour $\nu=0$ proviennent d'isomorphismes alg\'ebriques compatibles ${\bf G}_{{\cal F  }}\simeq {\bf M}_{{\cal F}^{M}}$, ${\bf G}_{{\cal F  }}^{\nu}\simeq {\bf M}_{{\cal F}^M}^{\nu}$.

 \bigskip
 
 \section{Rel\`evement de facettes de sous-groupes de Levi}
 
   Soit $M$ un sous-groupe de Levi de $G$ contenant $A$. Soit $({\cal F}_{M},\nu)\in Fac^*(M;A)_{G-comp}$.   L'ensemble $p_{M}^{-1}({\cal F}_{M})$ est l'ensemble des $x\in App(A)$ qui v\'erifient les relations
  
      $\alpha(x)=c_{\alpha,{\cal F}_{M}}$ pour tout $\alpha\in \Sigma^M_{{\cal F}_{M}}$;
 
  $c_{\alpha,{\cal F}_{M}}<\alpha(x)<c_{\alpha,{\cal F}_{M}}^+$ pour tout $\alpha\in \Sigma^M-\Sigma^M_{{\cal F}_{M}}$.

Il est clair que c'est une r\'eunion de facettes de $App(A)$. 

\ass{Lemme }{(i) Pour toute facette ${\cal F}\in Fac(G;A)$ qui coupe $p_{M}^{-1}({\cal F}_{M}^{\nu})$, resp. $p_{M}^{-1}(\bar{{\cal F}}_{M}^{\nu})$, on a $\nu\in {\cal N}({\cal F})$ et ${\cal F}^{\nu}={\cal F}\cap p_{M}^{-1}({\cal F}_{M}^{\nu})$, resp. ${\cal F}^{\nu}={\cal F}\cap p_{M}^{-1}(\bar{{\cal F}}_{M}^{\nu})$.

(ii) Pour toute facette ${\cal F}\in Fac(G;A)$, les conditions suivantes sont \'equivalentes:

(a) ${\cal F}\cap p_{M}^{-1}({\cal F}_{M}^{\nu})$ est un ouvert non vide de $p_{M}^{-1}({\cal F}^{M,\nu})$;

(b) ${\cal F}\cap p_{M}^{-1}(\bar{{\cal F}}_{M}^{\nu})$ est un ouvert non vide de $p_{M}^{-1}(\bar{{\cal F}}_{M}^{\nu})$;

(c) $\nu\in {\cal N}({\cal F})$, $M$ contient $M_{{\cal F},\nu}$ et ${\cal F}^M={\cal F}_{M}$. 
 }
 
 Preuve.  Fixons $n\in Norm_{M(F)}(A)\cap K_{{\cal F}_{M}}^{\nu}$. Les actions de $n$ sur $App(A)$ et sur $App^M(A)$ sont compatibles. Donc l'action de $n$ sur $App(A)$ conserve $p_{M}^{-1}({\cal F}_{M}^{\nu})$. Soit $x$ dans cet ensemble. On a $n(x)=x+e$ avec $e\in {\cal A}_{M}/{\cal A}_{G}$. Puisque $n\in Norm_{M(F)}(A)$, son image dans $W$ appartient \`a $W^M$ et fixe tout point de ${\cal A}_{M}/{\cal A}_{G}$ Donc $n^k(x)=x+ke$ pour tout $k\in {\mathbb N}$. Parce que $\nu$ appartient \`a ${\cal N}({\cal F}_{M})_{G-comp}$, les images des $n^k$ dans $M_{ad}(F)$, donc aussi leurs images dans $G_{AD}(F)$, restent dans un groupe compact. Il en r\'esulte que les $n^k(x)$ restent dans un sous-ensemble compact de $App(A)$. Cela implique $e=0$. 
  Donc l'action de $n$ fixe tout point de $p_{M}^{-1}({\cal F}_{M}^{\nu})$. Un point de $p_{M}^{-1}({\cal F}_{M})$ qui est fixe par $n$ se projette forc\'ement  dans ${\cal F}_{M}^{\nu}$. Donc  $p_{M}^{-1}({\cal F}_{M}^{\nu})$ est plus pr\'ecis\'ement le sous-ensemble des \'el\'ements de $p_{M}^{-1}({\cal F}_{M})$ qui sont fixes par $n$. Pour une facette ${\cal F}\in Fac(G;A)$ coupant $p_{M}^{-1}({\cal F}_{M}^{\nu})$, $n$ fixe un point de ${\cal F}$ donc conserve cette facette, c'est-\`a-dire que $n\in K_{{\cal F}}^{\dag}$. Puisque $w_{G}(n)=\nu$, on a $n\in K_{{\cal F}}^{\nu}$, a fortiori cet ensemble est non vide. D'autre part, puisque $p_{M}^{-1}({\cal F}_{M})$ est r\'eunion de facettes, on a ${\cal F}\subset p_{M}^{-1}({\cal F}_{M})$. Donc ${\cal F}^{\nu}$, qui est l'ensemble des points fixes par l'action de $n$ dans ${\cal F}$, est \'egal \`a ${\cal F}\cap p_{M}^{-1}({\cal F}_{M}^{\nu})$. Cela prouve les assertions du (i) de l'\'enonc\'e concernant $p_{M}^{-1}({\cal F}_{M}^{\nu})$. D'apr\`es 5(4), l'ensemble $\bar{{\cal F}}_{M}^{\nu}$ est r\'eunion des ${\cal F}_{M,1}^{\nu}$  pour les facettes ${\cal F}_{M,1}$ de $Fac(M;A)$ telles que   ${\cal F}_{M,1}\subset \bar{{\cal F}}_{M}$ et $\nu\in {\cal N}^M({\cal F}_{M,1})$. Cette derni\`ere condition \'equivaut \`a $\nu\in {\cal N}^M({\cal F}_{M,1})_{G-comp}$ puisque $\nu\in {\cal N}^M_{G-comp}$.  En appliquant ce que l'on vient de prouver \`a chacune de ces facettes, on obtient les assertions restantes  de (i).
 
 L'ensemble $ p_{M}^{-1}({\cal F}_{M}^{\nu})$ est un ouvert dense de $p_{M}^{-1}(\bar{{\cal F}}_{M}^{\nu})$. Cela entra\^{\i}ne que les conditions (a) et (b) du (ii) sont \'equivalentes. 
  Soit ${\cal F}\in Fac(G;A)$ v\'erifiant le (a) du (ii). D'apr\`es le (i) d\'ej\`a prouv\'e, on a  $\nu\in {\cal N}({\cal F})$ et ${\cal A}_{{\cal F}}^{\nu}/{\cal A}_{G}$ est engendr\'e par les \'el\'ements $x-y$ pour $x,y\in {\cal F}\cap p_{M}^{-1}({\cal F}_{M}^{\nu})$. Puisque cet ensemble est ouvert dans $p_{M}^{-1}({\cal F}_{M}^{\nu})$ par hypoth\`ese et puisque $p_{M}^{-1}({\cal F}_{M}^{\nu})$ est invariant par translation par ${\cal A}_{M}/{\cal A}_{G}$ par construction, l'ensemble de ces \'el\'ements $x-y$ contient un voisinage de $0$ dans ${\cal A}_{M}/{\cal A}_{G}$. A fortiori le sous-espace engendr\'e par ces \'el\'ements contient ${\cal A}_{M}/{\cal A}_{G}$, donc ${\cal A}_{{\cal F}}^{\nu}$ contient ${\cal A}_{M}$. Cela entra\^{\i}ne que $M$ contient $M_{{\cal F},\nu}$.  Par d\'efinition de la facette ${\cal F}^M$, la projection de ${\cal F}$ est contenue dans ${\cal F}^M$. Cette projection contient celle de ${\cal F}\cap p_{M}^{-1}({\cal F}_{M}^{\nu})$, laquelle est contenue dans ${\cal F}_{M}$. D'o\`u l'\'egalit\'e ${\cal F}^M={\cal F}_{M}$, ce qui d\'emontre que le (a) du (ii) implique le (c). La r\'eciproque r\'esulte du (vi) du lemme 6. $\square$
  
  \bigskip
  
  \section{Distributions invariantes}
  
  Pour tout sous-groupe de Levi $M$ de $G$, on fixe une mesure de Haar sur $M(F)$. 
  
  Notons $C_{c}^{\infty}(G(F))$ l'espace des fonctions sur $G(F)$, \`a  valeurs complexes, localement constantes et \`a support compact. Pour $f\in C_{c}^{\infty}(G(F))$ et $g\in G(F)$, on note $^gf$ la fonction $x\mapsto f(g^{-1}xg)$ sur $G(F)$. On appelle distribution sur $G(F)$  une forme lin\'eaire sur $C_{c}^{\infty}(G(F))$ et distribution invariante une telle forme lin\'eaire $D$ telle que $D(^gf)=D(f)$ pour toute $f\in C_{c}^{\infty}(G(F))$ et tout $g\in G(F)$.

  Soit $M$ un sous-groupe  de Levi de $G$. Fixons un sous-groupe parabolique $P\in {\cal P}(M)$  et  munissons $U_{P}(F)$ d'une mesure de Haar.  Pour $f\in C_{c}^{\infty}(G(F))$, on d\'efinit une fonction $f_{[P]}\in C_{c}^{\infty}(M(F))$  par $f_{[P]}(m)=\delta_{P}(m)^{1/2}\int_{U_{P}(F)} f(mu)\,du$ pour tout $m\in M(F)$,  o\`u $\delta_{P}$ est le module usuel.  Fixons une facette sp\'eciale ${\cal F}_{sp}\in Fac(G)$, posons $K_{sp}^0=K_{{\cal F}_{sp}}^0$.   On a la d\'ecomposition $G(F)=M(F)U_{P}(F)K_{sp}^0$.  Il existe une constante $c>0$ telle  que l'on ait l'\'egalit\'e
 $$\int_{G(F)}f(g)\,dg=c\int_{M(F)\times U_{P}(F)\times K_{sp}^0}f(muk)\,dk\,du\,dm$$
 pour tout $f\in C_{c}^{\infty}(G(F))$. Pour toute telle fonction, on d\'efinit $f_{P}\in C_{c}^{\infty}(M(F))$ par
 $$f_{P}(m)=c\delta_{P}(m)^{1/2}\int_{U_{P}(F)\times K_{sp}^0}f(k^{-1}muk)\,dk\,du$$
 $$=c\int_{K_{sp}^0}(^kf)_{[P]}(m).$$
 Cette fonction ne d\'epend pas de la mesure sur $U_{P}(F)$. Elle d\'epend de $P$ et ${\cal F}_{sp}$ mais on sait que, pour toute distribution invariante $D^M$ sur $M(F)$, $D^M(f_{P})$ ne d\'epend pas de ces choix (par contre, ce terme d\'epend des mesures fix\'ees sur $M(F)$ et $G(F)$). Pour toute telle distribution $D^M$, on d\'efinit la distribution induite $ind_{M}^G(D^M)$ \`a $G(F)$ par l'\'egalit\'e
 $$ind_{M}^G(D^M)(f)=D^M(f_{P})$$
 pour toute $f\in C_{c}^{\infty}(G(F))$.  On peut formuler autrement cette d\'efinition. On oublie le groupe $K_{sp}^0$, on munit   $P(F)=M(F)U_{P}(F)$ de la mesure produit et $P(F)\backslash G(F)$ de la mesure invariante \`a droite compatible avec celles fix\'ees sur $P(F)$ et $G(F)$. Rappelons que cette mesure s'applique non pas \`a des fonctions invariantes \`a gauche  par $P(F)$ mais \`a des fonctions $\varphi$ sur $G(F)$ v\'erifiant $\varphi(mug)=\delta_{P}(m)\varphi(g)$ pour tous $m\in M(F)$, $u\in U_{P}(F)$ et $g\in G(F)$. Pour $f\in C_{c}^{\infty}(G(F))$, on a
    $$ind_{M}^G(D^M)(f)=\int_{P(F)\backslash G(F)}D^M((^gf)_{[P]})\, dg.$$

 \bigskip
 
  \section{Fonctions tr\`es cuspidales}

 Soit $P$ un sous-groupe parabolique de $G$. Fixons une composante de Levi $M$ de $P$ et une mesure de Haar sur $U_{P}(F)$.  Pour $f\in C_{c}^{\infty}(G(F))$, on a d\'efini une fonction $f_{[P]}$ sur $M(F)$.  Remarquons que la condition $f_{[P]}=0$ ne d\'epend ni du choix de $M$, ni de celui  de la mesure sur $U_{P}(F)$. On dit qu'une fonction $f\in C_{c}^{\infty}(G(F))$ est tr\`es cuspidale si $f_{[P]}=0$ pour tout sous-groupe parabolique propre $P\subsetneq G$. 
 
 Soient $f,\varphi\in C_{c}^{\infty}(G(F))$. Pour $g\in G(F)$, posons
 $$I(f,\varphi,g)=\int_{G(F)}f(g^{-1}xg)\varphi(x)\,dx.$$
 
 \ass{Lemme}{Supposons $\varphi$ tr\`es cuspidale. Alors l'image dans $A_{G}(F)\backslash G(F)$ du support de la fonction $g\mapsto I(f,\varphi,g)$ est compacte.}
 
 Ce lemme est bien connu mais il est plus rapide d'en donner une d\'emonstration que de trouver une r\'ef\'erence. Auparavant, introduisons une notation qui nous sera utile dans la suite. Pour $P\in {\cal F}(M_{min})$, on note $\Sigma(U_{P})$ le sous-ensemble des racines $\alpha\in \Sigma$ telles que $U_{\alpha}\subset U_{P}$. Si $M$ est la composante de Levi de $P$ contenant $M_{min}$, on note aussi $\Sigma^M$ l'analogue de $\Sigma$ quand on remplace $G$ par $M$, autrement dit l'ensemble des $\alpha\in \Sigma$ tels que $U_{\alpha}\subset M$.

 Preuve du lemme. Fixons un sous-groupe parabolique minimal $P_{min}\in {\cal P}(M_{min})$. On note $\Sigma^+=\Sigma(U_{P_{min}})$ et $A(F)^+$ l'ensemble des $a\in A(F)$ tels que $\vert \alpha(a)\vert _{F}\geq1$ pour tout $\alpha\in \Sigma^+$, o\`u $\vert .\vert _{F}$ est la valeur absolue usuelle de $F$. On sait qu'il existe un sous-ensemble compact $C$ de $G(F)$ tel que $G(F)=CA(F)^+C$.  Puisque les ensembles de fonctions $\{^cf; c\in C\}$ et $\{^{c^{-1}}\varphi; c\in C\}$ sont finis et que $^{c^{-1}}\varphi$ est tout aussi cuspidale que $\varphi$, il suffit de prouver que  le support de la fonction $a\mapsto I(f,\varphi,a)$ sur $A(F)^+$ est d'image compacte dans $A_{G}(F)\backslash A(F)^+$. Pour tout $P\in {\cal P}(M_{min})$ et tout entier $N>0$, notons $A(F)^+(P,N)$ l'ensemble des $a\in A(F)^+$  tels que $\vert \alpha(a)\vert _{F}>N$ pour tout $\alpha\in \Sigma(U_{P})$. Pour $N$ fix\'e, l'ensemble $A_{G}(F)\backslash A(F)^+$ est r\'eunion des $A_{G}(F)\backslash A(F)^+(P,N)$ quand $P$ parcourt les sous-groupes paraboliques  $P\in {\cal P}(M_{min})$ qui sont propres et maximaux, et d'un sous-ensemble compact. Il nous suffit donc de fixer $P\in {\cal P}(M_{min})$, $P\not=G$ et de prouver qu'il existe $N$ tel que $I(f,\varphi,a)=0$ pour tout $a\in A(F)^+(P,N)$.  Fixons une facette sp\'eciale ${\cal F}_{sp}\in Fac(G;A)$ et  posons $K=K_{{\cal F}_{sp}}^0$. Fixons un sous-groupe ouvert  et distingu\'e $K'$ de $K$ tel que $\varphi$ soit biinvariante par $K'$ et que $K'=(K'\cap U_{\bar{P}}(F))(K'\cap P(F))$. Montrons que
 
 (1) si $N$ est assez grand, pour tout $a\in A(F)^+(P,N)$, le support de la fonction $x\mapsto f(a^{-1}xa)\varphi(x)$ est contenu dans $(K'\cap U_{\bar{P}}(F))P(F)$.

 Posons  $K^+=K_{{\cal F}_{sp}}^+$, $\bar{U}_{K}^+=K^+\cap U_{\bar{P}_{min}}(F)$, $M_{min,K}=K\cap M_{min}(F)$, $U_{K}=K\cap U_{P_{min}}(F)$.  Pour tout $w\in W$, on peut relever $w$ en un \'el\'ement $n_{w}\in Norm_{G(F)}(A)\cap K$. L'image dans ${\bf G}_{{\cal F  }_{sp}}^0({\mathbb F}_{q})$ de $M_{min,K}U_{K}$ est l'ensemble des points sur ${\mathbb F}_{q}$ d'un sous-groupe de Borel de ${\bf G}_{{\cal F  }_{sp}}^0$. La d\'ecomposition de Bruhat dans ${\bf G}_{{\cal F}_{sp}}^0({\mathbb F}_{q})$ entra\^{\i}ne l'\'egalit\'e   $K=\cup_{w\in W}U_{K} \bar{U}_{K}^+n_{w}M_{min,K}U_{K}$.    Il est  \'egalement connu que $G(F)=KP_{min}(F)$. D'o\`u $G(F)=\cup_{w\in W}U_{K}\bar{U}_{K}^+n_{w} M_{min}(F)U_{P_{min}}(F)$. Soit $x\in G(F)$ que l'on \'ecrit $x=u\bar{u}n_{w}m_{min}u'$ conform\'ement \`a cette d\'ecomposition. Supposons $\varphi(x)\not=0$. Alors $m_{min}u'$ appartient au produit de $K$ et du support de $\varphi$, qui est compact. Cela entra\^{\i}ne que $m_{min}$, resp. $u'$, appartient \`a un certain sous-ensemble compact $C_{M_{min}}\subset M_{min}(F)$, resp. $C_{U_{min}}\subset U_{P_{min}}(F)$. Soit $a\in A(F)^+(P,N) $ et, pour tout $g\in G(F)$, posons $g^{a}=a^{-1}ga$. Alors $x^{a}=u^{a}\bar{u}^{a}a^{-1}w(a)n_{w}m_{min}^{a}{u'}^{a}$. Supposons que $f(x^{a})\not=0$. Alors $x^{a}$ appartient au support de $f$. 
 Parce que $a\in A(F)^+$, la conjugaison par $a^{-1}$ contracte $U_{P_{min}}(F)$, fixe $M_{min}(F)$ et dilate $U_{\bar{P}_{min}}(F)$. Puisque $u$, $m_{min}$ et $u'$ sont respectivement dans les compacts $U_{K}$, $C_{M_{min}}$ et $C_{U_{min}}$, les \'el\'ements $u^{a}$, $m_{min}^{a}$ et ${u'}^{a}$ restent dans des compacts (on entend par l\`a des compacts ind\'ependants de $x$ et $a$). Il en r\'esulte que $\bar{u}^{a}a^{-1}w(a)$ reste dans un compact. Mais c'est un \'el\'ement de $\bar{P}_{min}(F)$, donc les deux composantes  $\bar{u}^{a}$ et $a^{-1}w(a)$ restent dans des compacts. On peut \'ecrire $\bar{u}=\bar{u}_{\bar{P}}\bar{u}^{\bar{P}}$, o\`u $\bar{u}_{\bar{P}}\in  K^+\cap U_{\bar{P}}(F)$ et $\bar{u}^{\bar{P}}\in K^+\cap  M(F)\cap U_{\bar{P}_{min}}(F)$. De nouveau, les deux composantes $\bar{u}_{\bar{P}}^{a}$ et $\bar{u}^{\bar{P},a}$ restent dans des compacts. Notons $C_{U_{\bar{P}}}$ le premier. Alors $\bar{u}\in aC_{U_{\bar{P}}}a^{-1}$.   La conjugaison par $a$ agit dans $U_{\bar{P}}(F)$ par des racines de  valeurs absolues $<N^{-1}$ puisque $a\in A(F)^+(P,N)$. Si $N$ est assez grand, on a donc $aC_{U_{\bar{P}}}a^{-1}\subset K'\cap U_{\bar{P}}(F)$. Il en r\'esulte que $\bar{u}_{\bar{P}}\in K'\cap U_{\bar{P}}(F)$. Consid\'erons maintenant la deuxi\`eme condition: $a^{-1}w(a)$ reste dans un compact. Il est bien connu qu'il existe un sous-groupe de Levi $L\in {\cal F}(M_{min})$ de sorte que l'ensemble des points fixes de $w$ agissant sur ${\cal A}$ est ${\cal A}_{L}$ et que $w$ appartienne \`a $W^L$. La condition "$a^{-1}w(a)$ reste dans un compact" entra\^{\i}ne qu'il existe un sous-ensemble compact $C_{A}$ de $A(F)$ tel que $a$ reste dans $A_{L}(F)C_{A}$. Pour toute racine $\alpha\in \Sigma^L$, $\alpha(a)$ reste donc dans un compact de ${\mathbb R}_{>0}$. Si $N$ est assez grand, une telle racine ne peut intervenir ni dans $U_{P}$, ni dans $U_{\bar{P}}$ (car alors, on aurait $\vert \alpha(a)\vert _{F}>N$, resp. $\vert \alpha(a)\vert _{F}<N^{-1}$). Elle appartient donc \`a $\Sigma^M$, d'o\`u l'inclusion $\Sigma^L\subset \Sigma^M$, puis $L\subset M$. Donc $w\in W^M$. On suppose $N$ assez grand pour que les deux propri\'et\'es ci-dessus soient v\'erifi\'ees. On a \'ecrit $x=u\bar{u}n_{w}m_{min}u'=u\bar{u}_{\bar{P}}\bar{u}^{\bar{P}}n_{w}m_{min}u'$. On a vu que $\bar{u}_{\bar{P}}$ appartenait \`a $K'\cap U_{\bar{P}}(F)$. Les propri\'et\'es de $K'$ entra\^{\i}nent que $u\bar{u}u^{-1}$ appartient \`a $K'$ et peut s'\'ecrire $u\bar{u}u^{-1}=\bar{u}'p$, avec $\bar{u}'\in K'\cap U_{\bar{P}}(F)$ et $p\in K'\cap P(F)$. Alors $x=\bar{u}'pun_{w}m_{min}u'$. Le premier terme est dans $K'\cap U_{\bar{P}}(F)$ et tous les autres appartiennent \`a $P(F)$.  Cela d\'emontre (1).

 On peut fixer des mesures sur $U_{\bar{P}}(F)$ et $U_{P}(F)$ de sorte que l'on ait l'\'egalit\'e
 $$\int_{G(F)}f(a^{-1}xa)\varphi(x)\,dx=\int_{U_{\bar{P}}(F)\times M(F)\times U_{P}(F)}f(a^{-1}\bar{u}mua)\varphi(\bar{u}mu)\delta_{P}(m)\,d\bar{u}\,dm\,du.$$
 Supposons $N$ assez grand pour que (1) soit v\'erifi\'ee. On peut alors   restreindre l'int\'egration sur $U_{\bar{P}}(F)$ en une int\'egration sur $K'\cap U_{\bar{P}}(F)$. Puisque $\varphi$ est invariante par ce groupe, on a alors $\varphi(\bar{u}mu)=\varphi(mu)$. Il existe des compacts $C_{M}\subset M(F)$ et $C_{U_{P}}\subset U_{P}(F)$ tels que $\varphi(mu)\not=0$ entra\^{\i}ne $m\in C_{M}$ et $u\in C_{U_{P}}$. On a alors
  $$\int_{G(F)}f(a^{-1}xa)\varphi(x)\,dx=\int_{K'\cap U_{\bar{P}}(F)\times C_{M}\times C_{U_{P}}}f(a^{-1}\bar{u}mua)\varphi(mu)\,d\bar{u}\,dm\,du.$$
Quand $N$ est assez grand, $a^{-1}C_{U_{P}}a$  est aussi petit que l'on veut. On peut donc supposer que $f$ est invariante \`a droite par $a^{-1}C_{U_{P}}a$. Alors $f(a^{-1}\bar{u}mua)=f(a^{-1}\bar{u}ma)$. Mais alors, dans la formule ci-dessus, on voit appara\^{\i}tre l'int\'egrale
$$\int_{C_{U_{P}}}\varphi(mu)\,du,$$
qui n'est autre que
$$\int_{U_{P}(F)}\varphi(mu)\,du$$
par d\'efinition de $C_{U_{P}}$. Or cette derni\`ere int\'egrale est nulle puisque $\varphi$ est tr\`es cuspidale. Cela ach\`eve la d\'emonstration. $\square$

Fixons une mesure de Haar sur $A_{G}(F)$. Soit $\varphi\in C_{c}^{\infty}(G(F))$, supposons $\varphi$ tr\`es cuspidale. Le lemme permet d'associer \`a $\varphi$ la distribution $D_{\varphi}$ d\'efinie par
$$D_{\varphi}(f)=\int_{A_{G}(F)\backslash G(F)}\int_{G(F)}f(g^{-1}xg)\varphi(x)\,dx\,dg$$
pour toute $f\in C_{c}^{\infty}(G(F))$ (cette d\'efinition est similaire \`a celle de \cite{Co} paragraphe 3.2). C'est une distribution invariante. 

 Soit $M\in {\cal L}(M_{min})$. On fixe sur $A_{M}(F)$ la mesure de Haar telle que la mesure du plus grand sous-groupe compact $A_{M}(F)^c$ de ce groupe vaille $1$.
  Notons ${\cal A}_{M}^*$ l'espace vectoriel r\'eel dual de ${\cal A}_{M}$.  Il existe un r\'eseau $R_{M}\subset {\cal A}_{M}^*$ tel que $i{\cal A}_{M}^*/iR_{M}$ s'identifie au groupe des caract\`eres unitaires de $A_{M}(F)/A_{M}(F)^c$. On munit ${\cal A}_{M}^*$ de la mesure de Haar telle que la mesure de ${\cal A}_{M}^*/R_{M}$ vaille $1$. Notons $M(F)_{ell}$ l'ensemble des \'el\'ements semi-simples fortement r\'eguliers de $M(F)$ qui sont elliptiques, c'est-\`a-dire contenus dans un sous-tore maximal $T$ de $M$ tel que $A_{T}=A_{M}$. Fixons une facette sp\'eciale ${\cal F}_{sp}\subset App(A)$ et posons $K_{sp}^0=K_{{\cal F}_{sp}}^0$.  Pour une fonction $f\in C_{c}^{\infty}(G(F))$ et pour $x\in M(F)_{ell}$, on  d\'efinit l'int\'egrale orbitale pond\'er\'ee
  $$J_{M}^G(x,f)=D^G(x)^{1/2}\int_{A_{M}(F)\backslash G(F)}f(g^{-1}xg)v_{M}^G(g)\,dg,$$
  o\`u $D^G$ est un d\'eterminant usuel et $v_{M}^G$ est le poids  d\'efini par Arthur  \`a l'aide du groupe $K_{sp}^0$ et de la mesure fix\'ee sur ${\cal A}_{M}^*/{\cal A}_{G}^*$. La fonction $x\mapsto J_{M}^G(x,f)$ est invariante par conjugaison par $M(F)$.  Dans le cas o\`u $f$ est tr\`es cuspidale, on sait que cette fonction est ind\'ependante du choix de ${\cal F}_{sp}$, cf. \cite{BP} lemme 5.2.1.
  
 Soit $\varphi\in C_{c}^{\infty}(G(F))$, supposons $\varphi$ tr\`es cuspidale. La formule des traces locale (\cite{A3}) permet de r\'ecrire la d\'efinition de la distribution $D_{\varphi}$ sous la forme suivante:
$$(1) \qquad D_{\varphi}(f)=\sum_{M\in {\cal L}(M_{min})}\vert W^M\vert \vert W\vert ^{-1}(-1)^{dim(A_{M})-dim(A_{G})}\int_{M(F)_{ell}}f_{P_{M}}(m)J_{M}^G(m, \varphi)\,dm,$$
  o\`u, pour tout $M$, on a fix\'e   un \'el\'ement  quelconque $P_{M}\in {\cal P}(M)$. 

 \bigskip
 
 \section{Repr\'esentations de niveau $0$ et fonctions cuspidales associ\'ees}
 
  Soit $\pi$ une repr\'esentation admissible  de $G(F)$ dans un espace complexe $V$. On  suppose que $\pi$ est de longueur finie et de niveau $0$, ce qui signifie que $V$ est engendr\'e par les sous-espaces d'invariants $V^{K_{{\cal F}}^+}$ pour ${\cal F}\in Fac(G)$. On conservera cette repr\'esentation jusqu'\`a la fin de l'article.
 
      Pour toute facette ${\cal F}\in Fac(G)$, notons $V^{K_{{\cal F}}^+}$ le sous-espace des \'el\'ements de $V$ qui sont invariants par $K_{{\cal F}}^+$.  Il est de dimension finie. La restriction de $\pi$ au groupe $K_{{\cal F}}^{\dag}$ conserve cet espace. Il s'en d\'eduit une repr\'esentation de dimension finie du groupe $K_{{\cal F}}^{\dag}/K_{{\cal F}}^+=
       \sqcup_{\nu\in {\cal N}({\cal F})}{\bf G}_{{\cal F}}^{\nu}({\mathbb F}_{q})$.  On note $\pi_{{\cal F}}$ cette repr\'esentation et $trace\,\pi_{{\cal F}}$ son caract\`ere habituel. Soit $\nu\in {\cal N}({\cal F})$. La fonction $trace\,\pi_{{\cal F}}$ se restreint \`a ${\bf G}_{{\cal F}}^{\nu}({\mathbb F}_{q})$ en un \'el\'ement de $C^{inv}({\bf G}_{{\cal F}}^{\nu})$ que l'on note $\phi_{{\cal F},\nu}$ (ou plus pr\'ecis\'ement $\phi_{\pi,{\cal F},\nu}$). On note $\phi_{{\cal F},\nu,cusp}$  sa projection cuspidale dans $C^{inv}_{cusp}({\bf G}_{{\cal F}}^{\nu})$. Soit ${\cal F}'$ une facette dont l'adh\'erence contient ${\cal F}$ et telle que $\nu\in {\cal N}({\cal F}')$. Il correspond \`a ${\cal F}'$ un sous-espace parabolique ${\bf P}_{{\cal F}'}^{\nu}$ de ${\bf G}_{{\cal F}}^{\nu}$. Une composante de Levi ${\bf M}_{{\cal F}'}^{\nu}$ de cet espace s'identifie \`a ${\bf G}_{{\cal F}'}^{\nu}$. Modulo cette identification, on v\'erifie que
      
      (1) $\phi_{{\cal F}',\nu}=res_{{\bf M}_{{\cal F}'}^{\nu}}^{{\bf G}_{{\cal F}}^{\nu}}(\phi_{{\cal F},\nu})$.

Soit $M$ un sous-groupe de Levi de $G$. Soit $P\in {\cal P}(M)$. On d\'efinit le module de Jacquet $V_{P}$, quotient de $V$ par le sous-espace engendr\'e par les $\pi(u)(v)-v$ pour $v\in V$ et $u\in U_{P}(F)$. On d\'efinit la repr\'esentation $\pi_{P}$ de $M(F)$ dans $V_{P}$ de sorte que, pour $m\in M(F)$ et $v\in V$, la projection de $\pi(m)(v)$ dans $V_{P}$ soit $\delta_{P}(m)^{1/2}\pi_{P}(m)(v_{P})$, o\`u $v_{P}$ est la projection de $v$. Ainsi, pour tout $({\cal F}_{M},\nu_{M})\in Fac^*(M)$, on d\'efinit des fonctions $\phi_{\pi_{P},{\cal F}_{M},\nu_{M}}$ et $\phi_{\pi_{P},{\cal F}_{M},\nu_{M},cusp}$ sur ${\bf M}_{{\cal F}_{M}}^{\nu_{M}}({\mathbb F}_{q})$. Supposons $\nu_{M}\in {\cal N}^M_{G-comp}$, notons simplement $\nu$ cet \'el\'ement. 
  Choisissons une facette ${\cal F}\in Fac(G)$ telle que ${\cal F}$ coupe $p_{M}^{-1}({\cal F}_{M}^{\nu})$   selon un ensemble ouvert. D'apr\`es les lemmes 7 et 6 (que l'on peut appliquer car on ne perd rien \`a supposer que $A\subset M$ et ${\cal F}\subset App(A)$), $\nu$ appartient \`a ${\cal N}({\cal F})$ et on a un isomorphisme canonique ${\bf G}_{{\cal F}}^{\nu}({\mathbb F}_{q})\simeq  {\bf M}_{{\cal F}_{M}}^{\nu}({\mathbb F}_{q})$. Les fonctions $\phi_{\pi,{\cal F},\nu}$ et $\phi_{\pi,{\cal F},\nu,cusp}$ peuvent \^etre consid\'er\'ees comme des fonctions sur ce dernier groupe ${\bf M}_{{\cal F}_{M}}^{\nu}({\mathbb F}_{q})$. Montrons que
  
  (2) on a les \'egalit\'es $\phi_{\pi,{\cal F},\nu}=\phi_{\pi_{P},{\cal F}_{M},\nu}$ et  $\phi_{\pi,{\cal F},\nu,cusp}=\phi_{\pi_{P},{\cal F}_{M},\nu,cusp}$. 
  
  Preuve. D'apr\`es Moy et Prasad, la projection de $V$ dans $V_{P}$ se restreint en un isomorphisme de $V^{K_{{\cal F}}^+}$ sur $V_{P}^{K_{{\cal F}_{M}}^+}$, cf. \cite{MP} proposition 6.7. Si on oublie le $\delta_{P}$ de la d\'efinition du module de Jacquet, cet isomorphisme entrelace l'action  de ${\bf G}_{{\cal F}}^{\nu}({\mathbb F}_{q})$ sur le premier espace et celle de ${\bf M}_{{\cal F}_{M}}^{\nu}({\mathbb F}_{q})$ sur le second. Mais le $\delta_{P}$ vaut $1$ sur $K_{{\cal F}}^{\nu}\cap M(F)$ car l'image dans $M_{ad}(F)$ d'un \'el\'ement de cet ensemble est contenu dans un sous-groupe compact. La premi\`ere \'egalit\'e de (2) s'en d\'eduit. La seconde se d\'eduit de la premi\`ere par projection cuspidale. $\square$
  
  En cons\'equence,
  
  (3) les fonctions $\phi_{\pi,{\cal F},\nu}$ et $\phi_{\pi,{\cal F},\nu,cusp}$ ci-dessus ne d\'ependent pas du choix de ${\cal F}$ et les fonctions $\phi_{\pi_{P},{\cal F}_{M},\nu}$ et $\phi_{\pi_{P},{\cal F}_{M},\nu,cusp}$ ne d\'ependent pas du choix de $P$. 
  
  On notera simplement $\phi_{\pi,{\cal F}_{M},\nu}$ et $\phi_{\pi,{\cal F}_{M},\nu}$, ou encore plus simplement $\phi_{{\cal F}_{M},\nu}$ et $\phi_{{\cal F}_{M},\nu,cusp}$ ces fonctions. 
  
  Soient $M$ un sous-groupe de Levi de $G$ et $g\in G(F)$. Posons $M'=gMg^{-1}$. L'action de $g$ sur l'immeuble $Imm(G_{AD})$ envoie $Imm(M_{ad})$ sur $Imm(M'_{ad})$ et elle se descend en une bijection de $Imm(M_{AD})$ sur $Imm(M'_{AD})$ qui est compatible aux d\'ecompositions en facettes de ces immeubles. Soit $({\cal F}_{M},\nu)\in Fac^*(M)_{G-comp}$ et notons ${\cal F}_{M'}=g({\cal F}_{M})$ l'image de ${\cal F}_{M}$ par cette bijection. Il est clair que $({\cal F}_{M'},\nu)\in Fac^*(M')_{G-comp}$ et que, de la conjugaison par $g$ se d\'eduit un isomorphisme ${\bf M}_{{\cal F}_{M}}^{\nu}({\mathbb F}_{q})\simeq {\bf M}_{{\cal F}_{M'}}^{_{'}\nu}({\mathbb F}_{q})$. On a

 (4) cet isomorphisme transporte les fonctions $\phi_{{\cal F}_{M},\nu}$ et $\phi_{{\cal F}_{M},\nu,cusp}$ en les fonctions  $\phi_{{\cal F}_{M'},\nu}$ et $\phi_{{\cal F}_{M'},\nu,cusp}$.
 
 En effet, si on "rel\`eve" ${\cal F}_{M}$ en ${\cal F}$ comme ci-dessus, on peut relever ${\cal F}_{M'}$ en ${\cal F}'=g({\cal F})$. De la conjugaison par $g$ se d\'eduit un isomorphisme ${\bf G}_{{\cal F}}^{\nu}({\mathbb F}_{q})\simeq {\bf G}_{{\cal F}'}^{\nu}({\mathbb F}_{q})$. Celui-ci transporte la fonction $\phi_{{\cal F},\nu}$ en la fonction $\phi_{{\cal F}',\nu}$ parce que l'op\'erateur $\pi(g)$ entrelace les repr\'esentations $\pi_{{\cal F}}$ et $\pi_{{\cal F}'}$. L'assertion (3) r\'esulte alors de (2). $\square$

Pour $({\cal F},\nu)\in Fac^*(G)$, on a d\'efini ci-dessus des fonctions $\phi_{{\cal F},\nu}$ et $\phi_{{\cal F},\nu,cusp}$ sur l'espace ${\bf G}_{{\cal F}}^{\nu}({\mathbb F}_{q})$. Dans la suite, on les consid\'erera la plupart du temps comme des fonctions sur $G(F)$, \`a support dans $K_{{\cal F}}^{\nu}$ et biinvariantes par $K_{{\cal F}}^+$.  Notons  $Fac^*_{max}(G)$ l'ensemble des $({\cal F},\nu)\in Fac^*(G)$ tels que ${\cal F}^{\nu}$ est r\'eduit \`a un point.

\ass{Lemme}{Soit $({\cal F},\nu)\in Fac^*_{max}(G)$. Alors la fonction $\phi_{{\cal F},\nu,cusp}$ est tr\`es cuspidale.}

Preuve (cf. \cite{Co} paragraphe 3.2).  On peut supposer ${\cal F}\subset App(A)$. Fixons une facette ouverte ${\cal F}_{min}\subset App(A)$ dont l'adh\'erence contient ${\cal F}$. Comme on le sait, on peut fixer une facette  sp\'eciale ${\cal F}_{sp}$ contenue dans l'adh\'erence de ${\cal F}_{min}$. On note $K_{sp}^0=K_{{\cal F}_{sp}}^0$, $K_{sp}^+=K_{{\cal F}_{sp}}^+$  et $K_{min}^0=K_{{\cal F}_{min}}^0$. Le groupe $K_{min}^0$ est contenu dans $K_{sp}^0$ et il existe une unique sous-groupe parabolique $P_{min}\in {\cal P}(M_{min})$ de sorte que $K_{min}^0=(K_{sp}^0\cap P_{min}(F))K_{sp}^+$.  Soit $P$ un sous-groupe parabolique propre de $G$. On peut fixer $g\in G(F)$ tel que $g^{-1}(P)$ contient $P_{min}$. On sait que $G(F)=K_{sp}^0P_{min}(F)$. On peut donc supposer $g\in K_{sp}^0$. On a la d\'ecomposition $K_{sp}^0=\cup_{w\in W}K_{min}^0n_{w}K_{min}^0$ o\`u, pour tout $w\in W$, $n_{w}$ est un rel\`evement de $w$ dans $K_{sp}^0\cap Norm_{G(F)}(A)$. En vertu de l'\'egalit\'e $K_{min}^0=(K_{sp}^0\cap P_{min}(F))K_{sp}^+$, on a aussi bien $K_{sp}^0=\cup_{w\in W}K_{min}^0n_{w}(K_{sp}^0\cap P_{min}(F))$. Puisque l'on peut multiplier $g$ \`a droite par un \'el\'ement de $P_{min}(F)$, on peut donc supposer $g=kn_{w}$ pour un $w\in W$ et un $k\in K_{min}^0$. En posant $Q=k^{-1}(P)$, on a alors $Q\in {\cal F}(M_{min})$. Notons $L$ la composante de Levi de $Q$ contenant $M_{min}$ et $M=k(L)$. Pour $m\in M(F)$, on a $\phi_{{\cal F},\nu,cusp,[P]}(m)=\int_{U_{Q}(F)}\phi_{{\cal F},\nu,cusp}(kluk^{-1})\,du$, o\`u $l=k^{-1}mk$. Mais $K_{min}^0\subset K_{{\cal F}}^0$ et la fonction $\phi_{{\cal F},\nu,cusp}$ est invariante par conjugaison par ce groupe. On obtient que la relation $\phi_{{\cal F},\nu,cusp,[P]}=0$ \'equivaut \`a $\phi_{{\cal F},\nu,cusp,[Q]}=0$. Autrement dit, en oubliant cette construction, on peut supposer que $P$ et $M$ contiennent $M_{min}$. Dans ce cas, $P$ d\'etermine un sous-groupe parabolique ${\bf P}$ de ${\bf G}_{{\cal F}}$ de sorte que l'image de $P(F)\cap K_{{\cal F}}^0$ dans ${\bf G}_{{\cal F}}({\mathbb F}_{q})$ soit ${\bf P}({\mathbb F}_{q})$. On veut prouver que $\phi_{{\cal F},\nu,cusp,[Q]}=0$. Puisque $\phi_{{\cal F},\nu,cusp}$ est \`a support dans $K_{{\cal F}}^{\nu}$, c'est clair si cet ensemble ne coupe pas $P(F)$. Supposons donc que $P(F)\cap K_{{\cal F}}^{\nu}\not=\emptyset$. Dans ce cas, il existe un sous-espace parabolique ${\bf P}^{\nu}$ de ${\bf G}_{{\cal F}}^{\nu}$ tel que l'image de $P(F)\cap K_{{\cal F}}^{\nu}$ dans ${\bf G}_{{\cal F}}^{\nu}({\mathbb F}_{q})$ soit ${\bf P}^{\nu}({\mathbb F}_{q})$. On voit qu'il existe un espace de Levi ${\bf M}^{\nu}$ de ${\bf P}^{\nu}$ de sorte que la fonction $ \phi_{{\cal F},\nu,cusp,[Q]}$ s'identifie \`a $res_{{\bf M}^{\nu}}^{{\bf G}_{{\cal F}}^{\nu}}(\phi_{{\cal F},\nu,cusp})$, \`a une constante pr\`es provenant des mesures, o\`u ici, $\phi_{{\cal F},\nu,cusp}$ est vue comme une fonction sur ${\bf G}_{{\cal F}}^{\nu}({\mathbb F}_{q})$. La nullit\'e cherch\'ee r\'esulte de la cuspidalit\'e de  cette fonction pourvu que ${\bf P}^{\nu}$ soit un sous-espace parabolique propre. Mais supposons que ${\bf P}^{\nu}={\bf G}_{{\cal F}}^{\nu}$. Alors $K_{{\cal F}}^{\nu}=(P(F)\cap K_{{\cal F}}^{\nu})K_{{\cal F}}^+$. On sait que $K_{{\cal F}}^+=(P(F)\cap K_{{\cal F}}^+)(U_{\bar{P}}(F)\cap K_{{\cal F}}^+)$ (o\`u $\bar{P}$ est le groupe contenant $M$ et oppos\'e \`a $P$). On a donc $K_{{\cal F}}^{\nu}=(P(F)\cap K_{{\cal F}}^{\nu})(U_{\bar{P}}(F)\cap K_{{\cal F}}^+)$. Consid\'erons $n\in K_{{\cal F}}^{\nu}\cap Norm_{G(F)}(A)$. L'\'egalit\'e pr\'ec\'edente entra\^{\i}ne $n\in P(F)U_{\bar{P}}(F)$. Mais l'intersection de $P(F)U_{\bar{P}}(F)$ et de $Norm_{G(F)}(A)$ est $Norm_{M(F)}(A)$. Donc $n\in M(F)$. 
  Introduisons le sous-tore ${\bf A}$ de ${\bf G}_{{\cal F}}^0$ comme en 4(8) et (9) et identifions ${\cal A}$ \`a $X_{*}({\bf A})\otimes_{{\mathbb Z}}{\mathbb R}$. L'\'egalit\'e ${\bf P}^{\nu}={\bf G}_{{\cal F}}^{\nu}$ entra\^{\i}ne que l'image de $M(F)\cap K_{{\cal F}}^0$ dans ${\bf G}_{{\cal F}}({\mathbb F}_{q})$ est ce groupe tout entier. Donc ${\cal A}_{{\cal F}}$ contient ${\cal A}_{M}$. Rappelons que ${\cal A}_{{\cal F}}^{\nu}$ est l'ensemble des points fixes de l'action de ${\bf n}$ dans ${\cal A}_{{\cal F}}$. Puisque $n\in M(F)$, cette action est triviale sur ${\cal A}_{M}$ donc ${\cal A}_{{\cal F}}^{\nu}$ contient ${\cal A}_{M}$. Or ${\cal A}_{M}\not={\cal A}_{G}$ puisque $P$ est propre alors que ${\cal A}_{{\cal F}}^{\nu}={\cal A}_{G}$ d'apr\`es l'hypoth\`ese que ${\cal F}^{\nu}$ est r\'eduit \`a un point. Cette contradiction ach\`eve la preuve. $\square$

   \bigskip
   
      \section{D\'efinition de distributions}
   
    Fixons un sous-ensemble $\underline{Fac}^*_{max}(G)$ des classes de conjugaison par $G(F)$ dans $Fac^*_{max}(G)$. 
   
   \ass{Lemme}{Soit $f\in C_{c}^{\infty}(G(F))$. Alors l'ensemble des $({\cal F},\nu)\in Fac^*_{max}(G)$ telles que l'int\'egrale 
   $$\int_{G(F)}f(x)\phi_{{\cal F},\nu,cusp}(x)\,dx$$
   soit non nulle est fini.}
   
   Preuve (cf. \cite{Co}). On peut \'evidemment se limiter aux $({\cal F},\nu)$ tels que $\nu$ appartienne \`a l'image finie du support de $f$ par l'application $w_{G}$. Ces \'el\'ements sont contenus dans  un ensemble fini de classes de conjugaison par $G(F)$. On peut donc fixer $(\underline{{\cal F}},\nu)\in \underline{Fac}^*_{max}(G)$ et se limiter \`a l'ensemble des $(g^{-1}(\underline{{\cal F}}),\nu)$ pour $g\in G(F)$, ou encore $g\in K_{\underline{{\cal F}}}^{\dag}\backslash G(F)$. Mais $\phi_{g^{-1}(\underline{{\cal F}}),\nu,cusp}(x)=\phi_{\underline{{\cal F}},\nu,cusp}(gxg^{-1})$, donc
   $$\int_{G(F)}f(x)\phi_{g^{-1}(\underline{{\cal F}}),\nu,cusp}(x)\,dx=I(f,\phi_{\underline{{\cal F}},\nu,cusp},g)$$
   avec la notation du paragraphe 9. D'apr\`es les lemmes 9 et 10, L'ensemble des $g$ tels que ce terme soit non nul est compact modulo $A_{G}(F)$, donc fini modulo $K_{\underline{{\cal F}}}^{\dag}$. Cela prouve le lemme. $\square$
   
   On d\'efinit une distribution invariante $\Theta_{\pi,cusp}$ par
   $$\Theta_{\pi,cusp}(f)=\sum_{({\cal F},\nu)\in Fac^*_{max}(G)}\int_{G(F)}f(x)\phi_{{\cal F},\nu,cusp}(x)\,dx.$$
  Gr\^ace au lemme 10, on peut d\'efinir comme au paragraphe 9 une distribution $D_{\phi_{{\cal F},\nu,cusp}}$ pour tout $({\cal F},\nu)\in Fac^*_{max}(G)$.  Il r\'esulte de la preuve ci-dessus que, pour tout $f\in C_{c}^{\infty}(G(F))$, l'ensemble des $({\cal F},\nu)\in \underline{Fac}^*_{max}(G)$ tels que $D_{\phi_{{\cal F},\nu,cusp}}(f)\not=0$ est fini qu'on a l'\'egalit\'e
  $$(1) \qquad \Theta_{\pi,cusp}(f)=\sum_{({\cal F},\nu)\in \underline{Fac}^*_{max}(G)} mes(A_{G}(F)\backslash K_{{\cal F}}^{\dag})^{-1}D_{\phi_{{\cal F},\nu,cusp}}(f).$$

 Soit $M$ un sous-groupe de Levi de $G$. Fixons un sous-groupe parabolique $P\in {\cal P}(M)$. On d\'efinit comme ci-dessus la distribution $\Theta_{\pi_{P},cusp}$ sur $M(F)$. On d\'efinit de fa\c{c}on \'evidente les ensembles $Fac^*_{max}(M)_{G-comp}$ et $\underline{Fac}^*_{max}(M)_{G-comp}$. On d\'efinit une autre distribution $\Theta^M_{\pi,cusp}$ sur $M(F)$ par
 $$\Theta^M_{\pi,cusp}(f)=\sum_{({\cal F}_{M},\nu)\in Fac^*_{max}(M)_{G-comp}}\int_{M(F)}f(x)\phi_{{\cal F}_{M},\nu,cusp}(x)\,dx$$
 pour $f\in C_{c}^{\infty}(M(F))$. Le groupe $P$ n'intervient pas dans cette d\'efinition. 
 On a
 
 (2) soit $f\in C_{c}^{\infty}(M(F))$; si le support de $f$ est form\'e d'\'el\'ements  compacts modulo $Z(G)$, on a l'\'egalit\'e  $\Theta^M_{\pi_{P},cusp}(f)=\Theta^M_{\pi,cusp}(f)$. 
 
 En effet, pour $({\cal F}_{M},\nu_{M})\in Fac^*_{max}(M)$ tel que $\nu_{M}\not\in {\cal N}^M_{G-comp}$, le support de $\phi_{\pi_{P},{\cal F}_{M},\nu_{M},cusp}$ est form\'e d'\'el\'ements qui ne sont pas compacts modulo $Z(G)$. Pour une fonction $f$ comme ci-dessus, l'int\'egrale 
 $$\int_{M(F)}f(x)\phi_{\pi_{P},{\cal F}_{M},\nu,cusp}(x)\,dx$$
 est donc nulle. Et si $\nu_{M}\in {\cal N}^M_{G-comp}$, les fonctions $\phi_{\pi_{P},{\cal F}_{M},\nu_{M},cusp}$ et $\phi_{{\cal F}_{M},\nu_{M},cusp}$ sont \'egales d'apr\`es 10(2). 
 
 \bigskip
 
 \section{Le th\'eor\`eme principal}
 On d\'efinit le caract\`ere-distribution de $\pi$ par $\Theta_{\pi}(f)=trace\,\pi(f)$ pour tout $f\in C_{c}^{\infty}(G(F))$. 
 
 \ass{Th\'eor\`eme}{Soit $f\in C_{c}^{\infty}(G(F))$. Supposons que le support de $f$ soit form\'e d'\'el\'ements compacts modulo $Z(G)$. Alors on a l'\'egalit\'e
 $$\Theta_{\pi}(f)=	\sum_{M\in {\cal L}(M_{min})}\vert W^M\vert \vert W^G\vert ^{-1} ind_{M}^G(\Theta_{\pi,cusp}^M)(f).$$}
 
 La d\'emonstration de ce th\'eor\`eme occupe les paragraphes 13 \`a 17.

  \bigskip
 
 \section{R\'esolutions de Schneider-Stuhler et Mayer-Solleveld}

  D'apr\`es Schneider et Stuhler, on dispose d'une r\'esolution de $V$:
 
$$(1) \qquad  0\to C_{n}\to C_{n-1}\to ... \to C_{0}\stackrel{\delta_{0}}{\to} V$$

\noindent dont nous allons rappeler la d\'efinition, cf. \cite{SS} th\'eor\`eme II.3.1. Pour $i=0,...,n$, $C_{i}$ est la somme directe sur les facettes ${\cal F}\in Fac(G)$ telles que $dim({\cal F})=i$ des espaces d'invariants $V^{K_{{\cal F}}^+}$. L'homomorphisme $C_{0}\to V$ est la somme des inclusions $V^{K_{{\cal F}}^+}\subset V$. Pour d\'efinir l'homomorphisme $C_{i}\to C_{i-1}$ pour $i>0$, on munit arbitrairement chaque facette de dimension positive d'une orientation. Autrement dit,  $E_{{\cal F}} ={\cal A}_{{\cal F}}/{\cal A}_{G}$ \'etant  un espace vectoriel r\'eel supportant ${\cal F}$, on choisit une demi-droite dans la droite r\'eelle $\bigwedge^{dim({\cal F})}(E_{{\cal F}})$. Si $dim({\cal F})>1$, l'orientation de ${\cal F}$ d\'etermine une orientation de toute facette  ${\cal F}'$ de dimension $dim({\cal F})-1$ dans l'adh\'erence de ${\cal F}$, d'o\`u un signe $\epsilon_{{\cal F},{\cal F}'}$ qui vaut $1$ si cette orientation est celle fix\'ee pour ${\cal F}'$ et $-1$ dans le cas contraire. Si $dim({\cal F})=1$ et ${\cal F}'$ est l'un des deux points limites de ${\cal F}$, $\epsilon_{{\cal F},{\cal F}'}$ vaut $1$ si l'orientation de ${\cal F}$ pointe dans la direction de ${\cal F}'$ et $-1$ dans le cas contraire. D'autre part, pour ${\cal F}$ et ${\cal F}'$ comme ci-dessus, on a $K_{{\cal F}'}^+\subset K_{{\cal F}}^+$, d'o\`u $V^{K_{{\cal F}}^+}\subset V^{K_{{\cal F}'}^+}$. Alors l'homomorphisme $C_{i}\to C_{i-1}$ est la somme sur ${\cal F}$ de dimension $i$  et ${\cal F}'$ dans l'adh\'erence de ${\cal F}$ des inclusions $V^{K_{{\cal F}}^+}\subset V^{K_{{\cal F}'}^+}$ multipli\'ees par le signe $\epsilon_{{\cal F},{\cal F}'}$. Chaque $C_{i}$ est muni d'une action de $G(F)$. Un \'el\'ement $g\in G(F)$ envoie un \'el\'ement $v\in V^{K_{{\cal F}}^+}$ sur $\pi(g)v\in V^{K_{g{\cal F}}^+}$ multipli\'e, si $i>0$, par $+1$ si l'action de $g$ envoyant ${\cal F}$ sur $g{\cal F}$ conserve l'orientation, par $-1$ sinon. La suite (1) est une suite exacte de repr\'esentations de $G(F)$.

  Fixons  ${\cal F}_{\star}\in Fac(G;A)$. On pose simplement $K_{\star}^{\dag}=K_{{\cal F}_{\star}}^{\dag}$, $K_{\star}^0=K_{{\cal F}_{\star}}^0$, $K_{\star}^+=K_{{\cal F}_{\star}}^+$. Pour tout entier $R>0$, on fixe un sous-ensemble $B_{R}\subset Imm(G_{AD})$ de sorte que les conditions suivantes soient v\'erifi\'ees:
 
 (2) $B_{R}$ est compact et est r\'eunion (finie) de facettes;
 
 (3) $B_{R}$ est invariant par $K_{\star}^{\dag}$;
 
 (4) l'intersection $B_{R}(A)= B_{R}\cap App(A)$ est convexe;
 
 (5) pour $R<R'$, on a $B_{R}\subset B_{R'}$ et $Imm(G_{AD})$ est r\'eunion des $B_{R}$ pour $R\in {\mathbb N}_{>0}$.
 
 On peut construire de tels ensembles de la fa\c{c}on suivante. D\'ecrivons ${\cal F}_{\star}$ par les relations (1) et (2) du paragraphe 4.   Rappelons que tous les ensembles $\Gamma_{\alpha}$  sont invariants par translations par  ${\mathbb Z}$. On d\'efinit $B_{R}(A)$ comme l'ensemble des $x\in App(A)$ v\'erifiant les relations  
 
   $c_{\alpha,{\cal F}_{\star}}-R\leq \alpha(x)\leq c_{\alpha,{\cal F}_{\star}}+R$ pour tout $\alpha\in \Sigma_{{\cal F}_{\star}}$;
 
  $c_{\alpha,{\cal F}_{\star}}-R\leq \alpha(x)\leq c_{\alpha,{\cal F}_{\star}}^++R$ pour tout $\alpha\in \Sigma-\Sigma_{{\cal F}_{\star}}$.
  
  Ces ensembles sont invariants par $K_{\star}^{\dag}\cap Norm_{G(F)}(A)$. On d\'efinit $B_{R}$ comme l'ensemble des $x\in Imm(G_{AD})$ tels qu'il existe $k\in K_{\star}^{\dag}$ de sorte que $kx\in B_{R}(A)$. On a bien $B_{R}(A)=B_{R}\cap App(A)$ en vertu de la propri\'et\'e suivante:
  
  (6) soient deux facettes ${\cal F},{\cal F}'\subset App(A)$; supposons qu'il existe $x\in K_{\star}^{\dag}$ telles que $x({\cal F})={\cal F}'$; alors il existe $n\in K_{\star}^{\dag}\cap Norm_{G(F)}(A)$ tel que $n({\cal F})={\cal F}'$.
  
  Preuve. On \'ecrit $x=n_{x}x'$ avec $n_{x}\in  K_{\star}^{\dag}\cap Norm_{G(F)}(A)$  et $x'\in K_{\star}^0$. En posant ${\cal F}''=n_{x}^{-1}({\cal F}')$, on a $x'({\cal F})={\cal F}''$ et il suffit de montrer qu'il existe $n''\in K_{\star}^{0}\cap Norm_{G(F)}(A)$ tel que $n''({\cal F})={\cal F}''$. En oubliant cette construction, on suppose $x\in K_{\star}^0$. On sait que s'il existe $g\in G(F)$ tel que $g({\cal F})={\cal F}'$, alors il existe $n_{1}\in Norm_{G(F)}(A)$ tel que $n_{1}({\cal F})={\cal F}'$ (cf. \cite{BT1} 7.4.1).  Fixons donc un tel $n_{1}$.  On a alors $n_{1}^{-1}x\in K_{{\cal F}}^{\dag}$ et on peut \'ecrire $n_{1}^{-1}x=n_{2}y$, avec $n_{2}\in Norm_{G(F)}(A)\cap K_{{\cal F}}^{\dag}$ et $y\in K_{{\cal F}}^0$. On sait que, pour $n',n''\in Norm_{G(F)}(A)$, l'\'egalit\'e $K_{\star}^0n'K_{{\cal F}}^0=K_{\star}^0n''K_{{\cal F}}^0$ entra\^{\i}ne que $n'$ appartient \`a $(K_{\star}^{0}\cap Norm_{G(F)}(A))n''(K_{{\cal F}}^0\cap Norm_{G(F)}(A))$, cf. \cite{HR} remarque 9. L'\'egalit\'e $x=n_{1}n_{2}y$ entra\^{\i}ne donc qu'il existe $n\in  K_{\star}^{0}\cap Norm_{G(F)}(A)$ et $n_{3}\in K_{{\cal F}}^0\cap Norm_{G(F)}(A)$ tels que $n_{1}n_{2}=nn_{3}$. Alors $x=nn_{3}y\in nK_{{\cal F}}^{0}$ et ${\cal F}'=x({\cal F})=n({\cal F})$. Cela d\'emontre (6). $\square$

 On voit alors que les ensembles $B_{R}$ v\'erifient les conditions requises.

 Pour $i=0,...,n$, on note $C_{i,R}$ le sous-espace de $C_{i}$ qui est la somme directe des $V^{K_{{\cal F}}^+}$ pour ${\cal F}$ de dimension $i$ et ${\cal F}\subset B_{R}$. On note $V_{R}=\delta_{0}(C_{0,R})$. D'apr\`es Mayer et Solleveld,  la suite
 
 $$0\to C_{n,R}\to...\to C_{0,R}\to V_{R}$$
 est une suite exacte de repr\'esentations de $K_{\star}^{\dag}$, cf. \cite{MS} th\'eor\`eme 2.4. La condition (5) entra\^{\i}ne que $V$ est r\'eunion des $V_{R}$ pour $R\in {\mathbb N}_{>0}$. 
 
Puisqu'on a fix\'e une mesure de Haar sur $G(F)$, l'espace $C_{c}^{\infty}(G(F))$ est une alg\`ebre de convolution. Cette alg\`ebre   agit dans $V$ par $\pi(f)v=\int_{G(F)}\pi(g)v\,f(g)\,dg$ pour $f\in C_{c}^{\infty}(G(F))$ et $v\in V$. Fixons un sous-groupe ouvert $H$ de $K_{\star}^0$ qui est distingu\'e dans $K_{\star}^{\dag}$. On note $C_{c}(K_{\star}^{\dag}/H)$ la sous-alg\`ebre des \'el\'ements de $C_{c}^{\infty}(G(F))$ qui sont biinvariantes par $H$ et \`a support dans $K_{\star}^{\dag}$. L'action de cette alg\`ebre conserve chaque sous-espace $V_{R}$. L'espace d'invariants $V^H$ \'etant de dimension finie, il est inclus dans $V_{R}$ si $R$ est assez grand. Dans ce cas, pour $f\in C_{c}(K_{\star}^{\dag}/H)$, la trace de l'op\'erateur $\pi(f)$ dans $V$ est \'egale \`a celle de sa restriction \`a $V_{R}$. On obtient que, pour $R$ assez grand,
 $$(7) \qquad trace\,\pi(f)=\sum_{i=0,...,n}(-1)^{i}trace\,\pi_{i,R}(f),$$
 o\`u $\pi_{i,R}$ est la repr\'esentation de $K_{\star}^{\dag}$ (ou de $C_{c}(K_{\star}^{\dag}/H)$) dans $C_{i,R}$. Il est facile de calculer chacune de ces traces. Soit ${\cal F}\in Fac(G)$. La restriction de $\pi$ au groupe $K_{{\cal F}}^{\dag}$  conserve l'espace $V^{K_{{\cal F}}^+}$. Comme on l'a dit au paragraphe 10, il s'en d\'eduit une repr\'esentation de dimension finie $\pi_{{\cal F}}$ du groupe  $K_{{\cal F}}^{\dag}/K_{{\cal F}}^+=\sqcup_{\nu\in {\cal N}({\cal F})}{\bf G}_{{\cal F  }}^{\nu}({\mathbb F}_{q})$, dont on note $trace\,\pi_{{\cal F}}$ la trace habituelle. On la  consid\`ere comme une fonction sur $K_{{\cal F}}^{\dag}$ invariante par $K_{{\cal F}}^+$. A cause des orientations intervenant dans la d\'efinition des actions de $G(F)$ sur les $C_{i}$, on doit introduire le caract\`ere $\epsilon_{{\cal F}}$ de $K_{{\cal F}}^{\dag}$ d\'efini ainsi: si $dim({\cal F})=0$, ce caract\`ere est trivial; si $dim({\cal F})>0$, pour $x\in K_{{\cal F}}^{\dag}$, $\epsilon_{{\cal F}}(x)$ vaut $1$ si l'action de $x$ sur ${\cal F}$ conserve l'orientation, $-1$ sinon.  Alors
 $$(8) \qquad trace\,\pi_{i,R}(f)=\sum_{{\cal F}\subset B_{R}, dim({\cal F})=i}\int_{K_{{\cal F}}^{\dag}}\epsilon_{{\cal F}}(x)f(x)trace\,\pi_{{\cal F}}(x)\,dx.$$
 Pour simplifier cette formule, montrons que
 
 (9) soient ${\cal F}\in Fac(G)$ et $\nu\in {\cal N}({\cal F})$; alors $(-1)^{dim({\cal F})}\epsilon_{{\cal F}}(k)=(-1)^{dim({\cal F}^{\nu})}$ pour tout $k\in K_{{\cal F}}^{\nu}$.
 
 Preuve. Si $dim({\cal F})=0$, tous ces signes valent $1$. Supposons $dim({\cal F})>0$. On ne perd rien \`a supposer que ${\cal F}\subset App(A)$. Comme $K_{{\cal F}}^0$ fixe tout point de ${\cal F}$, le caract\`ere $\epsilon_{{\cal F}}$ est trivial sur ce sous-groupe. On peut donc se limiter \`a prouver la formule de l'assertion pour $k\in K_{{\cal F}}^{\nu}\cap Norm_{G(F)}(A)$. Notons $x_{{\cal F}}$ le barycentre de ${\cal F}$. L'action de $k$ conserve ce point et il existe un \'el\'ement $w\in W$ tel que $ky-x_{{\cal F}}=w(y-x_{{\cal F}})$ pour tout $y\in App(A)$.  Le sous-espace 
 ${\cal A}_{{\cal F}}/{\cal A}_{G}$ de $ {\cal A}/{\cal A}_{G}$ est  engendr\'e par les $y-x_{{\cal F}}$ pour $y\in {\cal F}$. Il est conserv\'e par $w$. Par d\'efinition, $\epsilon_{{\cal F}}(k)$ est le d\'eterminant de l'action de $w$ dans ${\cal A}_{{\cal F}}/{\cal A}_{G}$ tandis que $dim({\cal F}^{\nu})$ est la dimension de l'espace des points fixes de cette action. Or $w$ est une isom\'etrie. Il est \'el\'ementaire de prouver que, pour une isom\'etrie $w$ d'un espace euclidien $E$, son d\'eterminant est \'egal \`a $(-1)^{dim(E)-dim(E^w)}$, $E^w$ \'etant l'espace des points fixes de $w$. Cela prouve (8). $\square$.

 En d\'ecomposant chaque $K_{{\cal F}}^{\dag}$ en r\'eunion des $K_{{\cal F}}^{\nu}$ pour $\nu\in {\cal N}({\cal F})$, on d\'eduit de (7), (8) et (9) l'\'egalit\'e
 
 $$(10) \qquad trace\, \pi(f)=\sum_{({\cal F},\nu)\in Fac^*(G); {\cal F}\subset B_{R}} (-1)^{dim({\cal F}^{\nu})} \int_{K_{{\cal F}}^{\nu}} f(k)trace\,\pi_{{\cal F}}(k)\,dk.$$
 
 La somme en ${\cal F}$ est finie puisque ${\cal F}\subset B_{R}$. La somme en $\nu$ ne l'est pas forc\'ement, mais on peut imposer de plus \`a $\nu$ d'appartenir \`a l'image par $w_{G}$ du support de $f$. Alors la somme devient finie. Rappelons que cela vaut pour toute $f\in C_{c}(K_{\star}^{\dag}/H)$ et tout $R$ assez grand.

 \bigskip
 
 \section{Une variante de la formule de Meyer-Solleveld}
  Consid\'erons l'ensemble $Trip(G)$ des triplets $({\cal F},{\cal F}',\nu)$ tels que:

- ${\cal F}$ et ${\cal F}'$ sont deux facettes de l'immeuble et ${\cal F}$ appartient \`a l'adh\'erence $\overline{{\cal F}'}$ de ${\cal F}'$;

- $\nu$ appartient \`a ${\cal N}({\cal F})\cap {\cal N}({\cal F}')$.

Deux tels triplets $({\cal F},{\cal F}',\nu)$ et $({\cal F}_{1},{\cal F}_{1}',\nu_{1})$ sont dits conjugu\'es si $\nu_{1}=\nu$ et s'il existe $g\in G(F)$ tel que $g{\cal F}={\cal F}_{1}$ et $g{\cal F}'={\cal F}'_{1}$.  Pour tout $({\cal F},{\cal F}', \nu)\in Trip(G)$, on va d\'efinir un nombre r\'eel $z({\cal F},{\cal F}',\nu)$.
On ne perd rien \`a supposer que ${\cal F}'$ (donc aussi ${\cal F}$) est contenue dans $App(A)$. 
Introduisons le sous-tore ${\bf A}$ de ${\bf G}_{{\cal F}}$ comme en 4(8) et (9). 
La facette ${\cal F}'$ d\'etermine un sous-espace parabolique ${\bf P}_{{\cal F}'}^{\nu}$ de ${\bf G}_{{\cal F}}^{\nu}$ et un sous-tore ${\bf A}_{{\cal F}'}^{\nu}$ de ${\bf A}$, dont le commutant  ${\bf M}_{{\cal F}'}^{\nu}$ dans ${\bf G}_{{\cal F}}^{\nu}$ est une composante de Levi de ${\bf P}_{{\cal F}'}^{\nu}$. On sait que  ${\bf P}_{{\cal F}'}^{\nu}$ d\'etermine une chambre ouverte $C_{{\cal F}'}^{{\cal F},\nu}$ dans l'espace $X_{*}({\bf A}_{{\cal F}'}^{\nu})\otimes_{{\mathbb Z}}{\mathbb R}$. Celui-ci s'identifie \`a ${\cal A}_{{\cal F}'}^{\nu}$. Or on a muni ${\cal A}$, donc aussi  ${\cal A}_{{\cal F}'}^{\nu}$, d'une structure d'espace euclidien. Notons $B$ la boule de centre $0$ et de rayon $1$ dans cet espace. On d\'efinit $z({\cal F},{\cal F}',\nu)$ comme le quotient de la mesure de $B\cap C_{{\cal F}'}^{{\cal F},\nu}$ par la mesure de $B$. Il est clair que

(1) pour deux \'el\'ements conjugu\'es $({\cal F},{\cal F}',\nu)$ et ${\cal F}_{1},{\cal F}_{1}',\nu)$, on a $z({\cal F},{\cal F}',\nu)=z({\cal F}_{1},{\cal F}_{1}',\nu')$.

Montrons que

(2) soit $({\cal F},\nu)\in Fac^*(G)$ et soit ${\bf M}^{\nu}$ un espace de Levi de ${\bf G}_{{\cal F}}^{\nu}$; d'apr\`es 5(6), \`a tout ${\bf P}^{\nu}\in {\cal P}({\bf M}^{\nu})$ est associ\'ee une facette, notons-la ${\cal F}_{{\bf P}^{\nu}}$, et on a $({\cal F},{\cal F}_{{\bf P}^{\nu}},\nu)\in Trip(G)$; alors 
$$\sum_{{\bf P}^{\nu}\in {\cal P}({\bf M}^{\nu})}z({\cal F},{\cal F}_{{\bf P}^{\nu}},\nu)=1.$$

  Pour le prouver, on peut encore supposer que ${\cal F}\subset App(A)$ et que ${\bf M}$ contient ${\bf A}$. On note ${\bf A}_{{\bf M}}^{\nu}$ le plus grand tore d\'eploy\'e de ${\bf G}_{{\cal F}}$ commutant \`a ${\bf M}^{\nu}$ et on pose ${\cal A}_{{\bf M}}^{\nu}=X_{*}({\bf A}_{{\bf M}}^{\nu})\otimes_{{\mathbb Z}}{\mathbb R}$. Les chambres dans ${\cal A}_{{\bf M}}^{\nu}$ associ\'ees aux diff\'erents ${\bf P}^{\nu}\in {\cal P}({\bf M}^{\nu})$ sont disjointes et leur r\'eunion est dense dans ${\cal A}_{{\bf M}}^{\nu}$. Donc la somme des mesures des intersections de ces chambres avec la boule de centre $0$ et de rayon $1$ dans ${\cal A}_{{\bf M}}^{\nu}$ est la mesure de cette boule. Cela prouve (2). $\square$

  Pour $({\cal F},\nu)\in Fac^*(G)$ et pour un entier $R>0$, notons $E_{R}({\cal F},\nu)$ l'ensemble des ${\cal F}_{1}\in Fac(G)$ telles que ${\cal F}_{1}\subset B_{R}$ et $({\cal F}_{1},{\cal F},\nu)\in Trip(G)$. Posons
  $$z_{R}({\cal F},\nu)=\sum_{{\cal F}_{1}\in E_{R}({\cal F},\nu)}(-1)^{dim({\cal F}_{1}^{\nu})}z({\cal F}_{1},{\cal F},\nu).$$

 Pour $({\cal F},\nu)\in Fac^*(G)$, on a d\'efini les fonctions $\phi_{{\cal F},\nu}$ et  $\phi_{{\cal F},\nu,cusp}$ sur ${\bf G}_{{\cal F}}^{\nu}({\mathbb F}_{q})$, que l'on peut consid\'erer comme des fonctions sur  $G(F)$.
   Montrons que, pour $f\in C_{c}(K_{\star}^{\dag}/H)$ et pour $R$ assez grand, on a l'\'egalit\'e
 
 $$(3) \qquad trace\, \pi(f)=\sum_{({\cal F},\nu)\in Fac^*(G) }  z_{R}({\cal F},\nu) \int_{K_{{\cal F}}^{\nu}} f(k) \phi_{{\cal F},\nu,cusp}(k)\,dk.$$

 Preuve.  Fixons $({\cal F},\nu)\in Fac^*(G)$ telle que ${\cal F}\subset B_{R}$. On d\'efinit une fonction  $f_{{\cal F},\nu}$ sur $K_{{\cal F}}^{\nu}$ par $f_{{\cal F},\nu}(g)=\int_{K_{{\cal F}^+}}f(gk)\,dk$. On peut la  descendre en une  fonction sur $ {\bf G}_{{\cal F}}^{\nu}({\mathbb F}_{q})$. Alors
 $$\int_{K_{{\cal F}}^{\nu}}f(k)trace\,\pi_{{\cal F}}(k)\,dk\,=\sum_{{\bf G}_{{\cal F}}^{\nu}({\mathbb F}_{q})}f_{{\cal F},\nu}(x)\phi_{{\cal F},\nu}(x).$$
 On calcule $\phi_{{\cal F},\nu}$ en appliquant la formule 2(5). On a expliqu\'e en 5 que les sous-espaces paraboliques de ${\bf G}_{{\cal F}}^{\nu}$ correspondaient bijectivement aux facettes ${\cal F}'$ dont l'adh\'erence contient ${\cal F}$ et telles que $\nu\in {\cal N}({\cal F}')$. Autrement dit aux facettes ${\cal F}'$ telles que $({\cal F},{\cal F}',\nu)\in Trip(G)$. D'apr\`es (2), pour une telle facette ${\cal F}'$ correspondant \`a un sous-espace ${\bf P}_{{\cal F}'}^{\nu}$, on peut choisir $z({\bf P}_{{\cal F}'}^{\nu})=z({\cal F},{\cal F}',\nu)$.  Un espace de Levi ${\bf M}_{{\cal F}'}^{\nu}$ s'identifie \`a ${\bf G}_{{\cal F}'}^{\nu}$ et   $res_{{\bf M}_{{\cal F}'}^{\nu}}^{{\bf G}_{{\cal F}}^{\nu}}(\phi_{{\cal F},\nu})$ s'identifie \`a $\phi_{{\cal F}',\nu}$ d'apr\`es 10(1). Il en r\'esulte que $proj_{cusp,{\bf M}_{{\cal F}'}^{\nu}}(\phi_{{\cal F},\nu})$ s'identifie \`a $\phi_{{\cal F}',\nu,cusp}$. A ce point, la formule 2(5) donne
 $$\sum_{{\bf G}_{{\cal F}}^{\nu}({\mathbb F}_{q})}f_{{\cal F},\nu}(x)\phi_{{\cal F},\nu}(x)=\sum_{x\in {\bf G}_{{\cal F}}^{\nu}({\mathbb F}_{q})}\sum_{{\cal F}'\in Fac(G); ({\cal F},{\cal F}',\nu)\in Trip(G)}z({\cal F},{\cal F}',\nu) f_{{\cal F},\nu}(x)\phi_{{\cal F}',\nu,cusp}[{\bf P}_{{\cal F}'}^{\nu}](x).$$
 Fixons ${\cal F}'$ intervenant dans le membre de droite. On voit que
 $$\sum_{x\in {\bf G}_{{\cal F}}^{\nu}({\mathbb F}_{q})}f_{{\cal F},\nu}(x)\phi_{{\cal F}',\nu,cus}[{\bf P}_{{\cal F}'}^{\nu}](x)=\sum_{x\in {\bf P}_{{\cal F}}^{\nu}({\mathbb F}_{q})}f_{{\cal F},\nu}(x)\phi_{{\cal F}',\nu,cusp}[{\bf P}_{{\cal F}'}^{\nu}](x)$$
 $$=\sum_{x\in {\bf G}_{{\cal F}'}^{\nu}({\mathbb F}_{q})}f_{{\cal F}',\nu}(x)\phi_{{\cal F}',\nu,cusp} (x)=\int_{K_{{\cal F}'}^{\nu}}f(k)\phi_{{\cal F}',\nu,cusp}(k)\,dk.$$
 D'o\`u
  $$\int_{K_{{\cal F}}^{\nu}}f(k)trace\,\pi_{{\cal F}}(k)\,dk\,=\sum_{{\cal F}'\in Fac(G); ({\cal F},{\cal F}',\nu)\in Trip(G)}z({\cal F},{\cal F}',\nu) \int_{K_{{\cal F}'}^{\nu}}f(k)\phi_{{\cal F}',\nu,cusp}(k)\,dk.$$
 Reportons cette \'egalit\'e dans la formule 13(10). On obtient
 $$trace\,\pi(f)=\sum_{({\cal F},\nu)\in Fac^*(G); {\cal F}\subset B_{R}}(-1)^{dim({\cal F}^{\nu})}$$
 $$\sum_{{\cal F}'\in Fac(G); ({\cal F},{\cal F}',\nu)\in Trip(G)}z({\cal F},{\cal F}',\nu) \int_{K_{{\cal F}'}^{\nu}}f(k)\phi_{{\cal F}',\nu,cusp}(k)\,dk.$$
 On change ${\cal F}$ en ${\cal F}_{1}$ et ${\cal F}'$ en ${\cal F}$ et on intervertit les sommations. Il suffit d'appliquer la d\'efinition de $z_{R}({\cal F},\nu)$ pour obtenir la formule (3). $\square$
 
On a fix\'e une facette ${\cal F}_{\star}\subset App(A)$. Fixons une facette ouverte ${\cal F}_{\star\star}\subset App(A)$ dont l'adh\'erence contient ${\cal F}_{\star}$. On pose $K_{\star\star}^0=K_{{\cal F}_{\star\star}}^0$.  Ce groupe est contenu dans $K_{\star}^0$.   
  Pour toute $f\in C_{c}(K_{\star}^{\dag}/H)$, notons $f_{\star\star}$ la fonction sur $G(F)$ d\'efinie par $f_{\star\star}(g)=\int_{K_{\star\star}^0}f(k^{-1}gk)\,dk$.  Elle appartient encore \`a $C_{c}(K_{\star}^{\dag}/H)$. Montrons que la formule (3) entra\^{\i}ne

  $$(4) \qquad trace\, \pi(f)=\sum_{({\cal F},\nu)\in Fac^*(G;A) } mes(K_{\star\star}^0\cap K_{{\cal F}}^0)^{-1} z_{R}({\cal F},\nu) \int_{K_{{\cal F}}^{\nu}} f_{\star\star}(k) \phi_{{\cal F},\nu,cusp}(k)\,dk.$$ 
  
  Preuve. A tout \'el\'ement $({\cal F},\nu)\in Fac^*(G;A)$, associons l'ensemble $X_{{\cal F},\nu}$ des \'el\'ements $(k({\cal F}),\nu)$ quand $k$ d\'ecrit $K_{\star\star}^0$. On a
  
  (5) quand $({\cal F},\nu)$ d\'ecrit $Fac^*(G;A)$, la r\'eunion des $X_{{\cal F},\nu}$ est $Fac^*(G)$ tout entier. 
  
  En effet, soit $({\cal F}_{1},\nu)\in Fac^*(G)$. On sait que l'on peut fixer $g\in G(F)$ tel que ${\cal F}'=g({\cal F}_{1})$ soit contenue dans l'adh\'erence de ${\cal F}_{\star\star}$. Le groupe $K_{\star\star}^0$ est un sous-groupe d'Iwahori associ\'e \`a une facette de $App(A)$. On a donc la d\'ecomposition $G(F)=K_{\star\star}^0Norm_{G(F)}(A)K_{\star\star}^0$. Ecrivons $g^{-1}=k_{1}nk_{2}$ conform\'ement \`a cette d\'ecomposition. On a $K_{\star\star}^0\subset K_{{\cal F}'}^0$ donc $k_{2}$ fixe ${\cal F}'$. Alors ${\cal F}_{1}=k_{1}n({\cal F}')$. Le couple $(n({\cal F}'),\nu)$ appartient \`a $Fac^*(G;A)$ et on a $({\cal F}_{1},\nu)=k_{1}(n({\cal F}'),\nu)$. Cela d\'emontre (5).
  
  Pour deux \'el\'ements $({\cal F},\nu)$ et $({\cal F}',\nu')$ de $Fac^*(G;A)$ les ensembles $X({\cal F},\nu)$ et $X({\cal F}',\nu')$ sont disjoints ou \'egaux. Ils sont \'egaux si et seulement si $\nu=\nu'$ et il existe $k\in K_{\star\star}^0$ tel que ${\cal F}'=k({\cal F})$. Les facettes ${\cal F}$ et ${\cal F}'$ \'etant toutes deux dans $App(A)$, on sait que cette derni\`ere condition ne peut \^etre v\'erifi\'ee que si ${\cal F}'={\cal F}$ (cela se voit comme en  13(6)). Autrement dit, les ensembles $X({\cal F},\nu)$ et $X({\cal F}',\nu')$ sont toujours disjoints si $({\cal F},\nu)\not=({\cal F}',\nu')$.    On peut alors r\'ecrire la somme (3) comme une somme sur les $({\cal F},\nu)\in Fac^*(G;A)$ suivie d'une somme sur les \'el\'ements de $X({\cal F},\nu)$. Cette derni\`ere est aussi une somme sur les $(k({\cal F}),\nu)$ pour $k\in K_{\star\star}^0/K_{\star\star}^0\cap K_{{\cal F}}^0$. Le terme que l'on somme est
    $$ z_{R}(k({\cal F}),\nu)\int_{K_{k({\cal F})}^{\nu}} f(x) \phi_{k({\cal F}),\nu,cusp}(x)\,dx,$$
  L'int\'egrale se r\'ecrit
    $$ \int_{K_{{\cal F}}^{\nu}} f(kxk^{-1}) \phi_{k({\cal F}),\nu,cusp}(kxk^{-1})\,dx.$$
    Il r\'esulte de 10(4) que $\phi_{k({\cal F}),\nu,cusp}(kxk^{-1})=\phi_{{\cal F},\nu,cusp}(x)$.  D'autre part, parce que $B_{R}$ est invariant par $K_{\star}^{\dag}$, on voit que $z_{R}(k({\cal F}),\nu)=z_{R}({\cal F},\nu)$. La contribution de $X({\cal F},\nu)$ est donc
    $$z_{R}({\cal F},\nu) \int_{K_{{\cal F}}^{\nu}}\phi_{{\cal F},\nu,cusp}(x) \sum_{k\in K_{\star\star}^0/K_{\star\star}^0\cap K_{{\cal F}}^0 }f(kxk^{-1})\,dx.$$
    On voit que la somme int\'erieure est \'egale \`a
    $$mes(K_{\star\star}^0\cap K_{{\cal F}}^0)^{-1}f_{\star\star}(x).$$
    En  rassemblant ces calculs, la formule (3) devient (4). $\square$

\bigskip

 \section{Descente aux sous-groupes de Levi} 
  
   Consid\'erons un entier $N>0$.  Soit $P\in {\cal F}(M_{min})$.   On note $X'_{N}(P)$ l'ensemble des $x\in App(A)$ tels que $\alpha(x)>N$ pour tout $\alpha\in \Sigma(U_{P})$. Si $Q\in {\cal F}(M_{min})$ et $Q\subset P$, on a $U_{P}\subset U_{Q}$ donc $X'_{N}(Q)\subset X'_{N}(P)$. On pose
   $$X_{N}(P)=X'_{N}(P)-\bigcup_{Q\in {\cal F}(M_{min}), Q\subsetneq P}X'_{N}(Q).$$
   Les ensembles $X'_{N}(P)$ et $X_{N}(P)$ sont des r\'eunions de facettes (parce que, pour tout $\alpha\in \Sigma$, l'ensemble $\Gamma_{\alpha}$ contient ${\mathbb Z}$). On a
   
   (1) pour toute facette ${\cal F}\in Fac(G;A)$, l'ensemble des $P\in {\cal F}(M_{min})$ tels que ${\cal F}\subset X'_{N}(P)$ poss\`ede un unique \'el\'ement minimal; celui-ci est 
    l'unique $P\in {\cal F}(M_{min})$ tel que ${\cal F}\subset X_{N}(P)$.
   
   Preuve. Fixons $x\in {\cal F}$. Notons $\Sigma_{0}(x)$, resp. $\Sigma_{+}(x)$, l'ensemble des \'el\'ements $\alpha\in \Sigma$ tels que $\alpha(x)=0$, resp. $\alpha(x)>0$. L'ensemble $\Sigma_{0}(x)\cup \Sigma_{+}(x)$ est clos, c'est-\`a-dire que pour deux \'el\'ements $\alpha,\beta$ de cet ensemble tels que $\alpha+\beta\in \Sigma$, $\alpha+\beta$ appartient encore \`a cet ensemble. La r\'eunion de $\Sigma_{0}(x)\cup \Sigma_{+}(x)$ et de son oppos\'e $\Sigma_{0}(x)\cup (-\Sigma_{+}(x))$ est $\Sigma$ tout entier. Il en r\'esulte    qu'il existe $P_{x}\in {\cal F}(M_{min})$ tel que $\Sigma(P_{x})=\Sigma_{0}(x)\cup \Sigma_{+}(x)$  et $\Sigma(U_{P_{x}})=\Sigma_{+}(x)$. On note $M_{x}$ la composante de Levi de $P_{x}$ contenant $M_{min}$. Choisissons un sous-groupe parabolique minimal $P_{min}\in {\cal P}(M_{min})$ tel que $P_{min}\subset P_{x}$. Notons $\Delta$ la base de $\Sigma$ associ\'ee \`a $P_{min}$, $\Delta^{M_{x}}$ le sous-ensemble des $\alpha\in \Delta$ qui interviennent dans $M_{x}$ et notons $\Delta_{N}$ le sous-ensemble des $\alpha\in \Delta$ tels que $\alpha(x)>N$. Il correspond \`a cet ensemble un sous-groupe parabolique $P$ contenant $P_{min}$: $\Delta-\Delta_{N}$ est l'ensemble des racines simples intervenant dans la composante de Levi $M$ de $P$ qui contient $M_{min}$. Pour $\alpha\in \Delta^{M_{x}}$, on a $\alpha(x)=0$ donc $\alpha\not\in \Delta_{N}$. Il en r\'esulte que $M_{x}\subset M$ donc $P_{x}\subset P$. Pour $\alpha\in \Sigma(U_{P})$, $\alpha$ est combinaison lin\'eaire \`a coefficients entiers positifs d'\'el\'ements de $\Delta$ et au moins un \'el\'ement de $\Delta_{N}$ intervient avec un coefficient $\geq1$. Donc $\alpha(x)>N$. Cela entra\^{\i}ne que $x\in X'_{N}(P)$.  Il reste \`a  prouver que, pour tout $Q\in {\cal F}(M_{min})$ tel que $x\in X'_{N}(Q)$, on a $P\subset Q$. Consid\'erons donc un tel sous-groupe parabolique $Q$. Pour $\alpha\in \Sigma(U_{Q})$, on a $\alpha(x)>N$ a fortiori $\alpha(x)>0$. Donc $\Sigma(U_{Q})\subset \Sigma_{+}(x)$, puis $U_{Q}\subset U_{P_{x}}$. Cela entra\^{\i}ne $P_{x}\subset Q$ d'o\`u aussi $P_{min}\subset Q$. A la composante de Levi $L$ de $Q$ contenant $M_{min}$ correspond un sous-ensemble $\Delta^L$ de $\Delta$. Pour $\alpha\in \Delta-\Delta_{L}$, on a $\alpha\in \Sigma(U_{Q})$ donc $\alpha(x)>N$, ce qui signifie que $\alpha\in \Delta_{N}$. Donc $L$ contient $M$ puis $Q$ contient $P$. Cela ach\`eve la preuve. $\square$
   
   Remarquons que
   
   (2) l'ensemble $X_{N}(G)$ est compact;
   
   Preuve.  C'est le compl\'ementaire de la r\'eunion des $X'_{N}(P)$ pour $P$ propre. Or chacun de ces ensembles est ouvert. Donc $X_{N}(G)$ est ferm\'e et il suffit de prouver qu'il est inclus dans un compact. 
   D'apr\`es la preuve pr\'ec\'edente, si $x$ est un point de $X_{N}(G)$, il existe $P_{min}\in {\cal P}(M_{min})$ tel qu'en notant $\Delta$ la base de $\Sigma$ associ\'ee \`a $P_{min}$, on ait $0\leq \alpha(x)\leq N$ pour tout $\alpha\in \Delta$. Pour chaque $P_{min}$, l'ensemble des $x$ v\'erifiant ces conditions est compact. Il n'y a qu'un nombre fini de $P_{min}$. $\square$
   
 Soit $P\in {\cal F}(M_{min})$, notons $M$ sa composante de Levi contenant $M_{min}$. L'application produit $U_{\bar{P}}(F)\times M(F)\times U_{P}(F)\to G(F)$ identifie l'ensemble de d\'epart \`a un ouvert dense de $G(F)$. On a d\'ej\`a fix\'e des mesures sur $G(F)$ et $M(F)$.  On munit les groupes $U_{\bar{P}}(F)$ et $U_{P}(F)$ de mesures de Haar de sorte que, pour toute fonction $f\in C_{c}^{\infty}(G(F))$, on ait l'\'egalit\'e
 $$\int_{G(F)}f(g)\,dg\,= \int_{U_{\bar{P}}(F)\times M(F)\times U_{P}(F)}f(\bar{u}mu)\delta_{P}(m)\,du\,dm\,d\bar{u}.$$
    
 \ass{Lemme}{Il existe un entier $N>0$  de sorte que, pour toute fonction $f\in C_{c}(K_{\star}^{\dag}/H)$ et pour tout $({\cal F},\nu)\in Fac^*(G;A)$ tel que ${\cal F}\subset X'_{N}(P)$, on ait les propri\'et\'es suivantes:
 
 (i)  si $M$ ne contient pas $ M_{{\cal F},\nu}$, 
 $$\int_{K_{{\cal F}}^{\nu}} f_{\star\star}(k) \phi_{{\cal F},\nu,cusp}(k)\,dk=0;$$
 
 (ii) si $M$ contient $M_{{\cal F},\nu}$, 
 $$\int_{K_{{\cal F}}^{\nu}} f_{\star\star}(k) \phi_{{\cal F},\nu,cusp}(k)\,dk=mes(K_{{\cal F}}^+\cap U_{\bar{P}}(F))\int_{K_{{\cal F}^M}^{\nu}} f_{\star\star,[P]}(k) \phi_{{\cal F}^M,\nu,cusp}(k)\,dk.$$}
 
 Rappelons que la fonction $\phi_{{\cal F}^M,\nu,cusp}$ a \'et\'e d\'efinie en 10.
 
 Preuve.  
 On fixe un entier $c\geq1$   de sorte que
 
 (3) pour tout $\alpha\in \Sigma$, 
 $\vert c_{\alpha,{\cal F}_{\star}}\vert <c$, $\vert c_{\alpha,{\cal F}_{\star}}^+\vert<c$ et $U_{\alpha,-c}\subset H$. 
 
  On  impose \`a $N$   de v\'erifier
 
 (4) $N\geq 3c$. 
 
 D'autre part, puisque $K_{\star}^{\dag}/Z(G)(F)$ est compact, on peut aussi imposer

(5) soit $m\in M(F)$ et, pour tout $\alpha\in \Sigma(U_{P})$, soit $u_{\alpha}\in U_{\alpha}(F)$; supposons $m\prod_{\alpha\in \Sigma(U_{P})}u_{\alpha}\in K_{\star}^{\dag}$ (pour un ordre fix\'e  quelconque sur $\Sigma(U_{P})$); alors $u_{\alpha}\in U_{\alpha,N-1}$ pour tout $\alpha\in \Sigma(U_{P})$.

Soit $({\cal F},\nu)\in Fac^*(G;A)$ tel que ${\cal F}\subset X'_{N}(P)$. Supposons d'abord que $M$ contient $M_{{\cal F},\nu}$. Gr\^ace au lemme 6, on a l'\'egalit\'e
$$\int_{K_{{\cal F}}^{\nu}} f_{\star\star}(k) \phi_{{\cal F},\nu,cusp}(k)\,dk=\int_{K_{{\cal F}}^+\cap U_{\bar{P}}(F)\times K_{{\cal F}^M}^{\nu}\times K_{{\cal F}}^+\cap U_{P}(F)}f_{\star\star}(\bar{u}mu)\phi_{{\cal F},\nu,cusp}(\bar{u}mu)\,du\,dm\,d\bar{u}.$$
Dans l'expression ci-dessus, on a $\phi_{{\cal F},\nu,cusp}(\bar{u}mu)=\phi_{{\cal F}^M,\nu,cusp}(m)$  d'apr\`es 10(2) et parce que $\phi_{{\cal F},\nu,cusp}$ est invariante par $K_{{\cal F}}^+$. Pour $\alpha\in \Sigma(U_{\bar{P}})$, on a $c_{\alpha,{\cal F}}<-N$ d'apr\`es l'hypoth\`ese ${\cal F}\subset X'_{N}(P)$. Les relations 4(7) et (3) et (4) ci-dessus entra\^{\i}nent $K_{{\cal F}}^+\cap U_{\alpha}(F)\subset H$. Dans l'expression ci-dessus, on a donc $\bar{u}\in H$, d'o\`u 
  $f_{\star\star}(\bar{u}mu)=f_{\star\star}(mu)$. La formule ci-dessus se r\'ecrit
$$(6) \qquad \int_{K_{{\cal F}}^{\nu}} f_{\star\star}(k) \phi_{{\cal F},\nu,cusp}(k)\,dk=mes(K_{{\cal F}}^+\cap U_{\bar{P}}(F))\int_{K_{{\cal F}^M}^{\nu}}\psi(m)\phi_{{\cal F}^M,\nu,cusp}(m)\,dm,$$
o\`u 
$$\psi(m)=\int_{K_{{\cal F}}^+\cap U_{P}(F)}f_{\star\star}(mu)\,du.$$
Mais $K_{{\cal F}}^+\cap U_{P}(F)$ est produit des $K_{{\cal F}}^+\cap U_{\alpha}(F)$ pour $\alpha\in \Sigma(U_{P})$. En utilisant 4(7), (5) ci-dessus et  l'hypoth\`ese ${\cal F}\subset X'_{N}(P)$, on voit que, pour tous $m\in M(F)$ et $u\in U_{P}(F)$, la condition $mu\in K_{\star}^{\dag}$ impose $u\in K_{{\cal F}}^+\cap U_{P}(F)$. Puisque $f_{\star\star}$ est \`a support dans $K_{\star}^{\dag}$, on a $f_{\star\star}(mu)=0$ si $u\not\in K_{{\cal F}}^+\cap U_{P}(F)$. On peut donc remplacer dans la d\'efinition de $\psi$ l'int\'egration sur $K_{{\cal F}}^+\cap U_{P}(F)$ par l'int\'egration sur $U_{P}(F)$ tout entier. Alors $\psi(m)=\delta_{P}(m)^{-1/2}f_{\star\star,[P]}(m)$. Mais  on int\`egre sur  $m\in K_{{\cal F}^M}^{\nu}$. L'image d'un tel $m$ dans $M_{ad}(F)$ appartient \`a un sous-groupe compact donc $\delta_{P}(m)=1$. On peut donc remplacer la fonction $\psi$ par $f_{\star\star,[P]}$ et  la formule (6) devient celle du (ii) de l'\'enonc\'e.

Supposons maintenant que $M$ ne contient pas $M_{{\cal F},\nu}$. Remarquons que, puisque $f_{\star\star}$ est \`a support dans $K_{\star}^{\dag}$, l'int\'egrale \`a calculer est nulle si $K_{\star}^{\dag}$ ne coupe pas $K_{{\cal F}}^{\nu}$. On peut donc supposer que  $K_{\star}^{\dag}$ coupe $K_{{\cal F}}^{\nu}$. Cela implique que $\nu\in {\cal N}({\cal F}_{\star})$ et que $K_{\star}^{\dag}\cap K_{{\cal F}}^{\nu}=K_{\star}^{\nu}\cap K_{{\cal F}}^{\nu}$. Notre int\'egrale vit en fait sur cet ensemble. Fixons un point $x_{\star}\in {\cal F}_{\star}^{\nu}$. Toute racine $\alpha\in \Sigma-\Sigma^{M_{{\cal F},\nu}}$ est non constante sur ${\cal F}^{\nu}$. On peut donc fixer un point $x\in {\cal F}^{\nu}$ tel que, en posant $v=x-x_{\star}$, on ait $\alpha(v)\not=0$ pour tout $\alpha\in \Sigma-\Sigma^{M_{{\cal F},\nu}}$. Consid\'erons l'ensemble des sous-groupes paraboliques $P'\in {\cal F}(M_{min})$ tels que $\alpha(v)>N-c$ pour tout $\alpha\in \Sigma(U_{P'})$.   La condition (3) et l'hypoth\`ese ${\cal F}\in X'_{N}(P)$ entra\^{\i}nent que $P$ appartient \`a cet ensemble. La m\^eme preuve qu'en (1) montre que cet ensemble poss\`ede un unique \'el\'ement minimal. On note $Q$ cet \'el\'ement et $L$ sa composante de Levi contenant $M_{min}$. On a donc $Q\subset P$ et $L\subset M$. A fortiori, $L$ ne contient pas $M_{{\cal F},\nu}$. On pose $Q'=M_{{\cal F},\nu}\cap Q$, $L'=M_{{\cal F},\nu}\cap L$. Alors $Q'$ est un sous-groupe parabolique propre de $M_{{\cal F},\nu}$ et $L'$ est sa composante de Levi contenant $M_{min}$.  Il existe une unique facette ${\cal F}'$ telle que ${\cal F}'$ coupe le segment $[x_{\star},x]$ selon un ouvert de ce segment et que $x$ appartienne \`a l'adh\'erence de cette intersection (il se peut que ${\cal F}'={\cal F}$). De m\^eme, il  existe une unique facette ${\cal F}_{\star}'$ telle que ${\cal F}_{\star}'$ coupe le segment $[x_{\star},x]$ selon un ouvert de ce segment et que $x_{\star}$ appartienne \`a l'adh\'erence de cette intersection. Montrons que

(7) $\nu\in{\cal N}({\cal F}'_{\star})$ et $\nu\in {\cal N}({\cal F}')$; $K_{\star}^{\nu}\cap K_{{\cal F}}^{\nu}\subset K_{{\cal F}_{\star}'}^{\nu}\subset K_{\star}^{\nu}$ et $K_{\star}^{\nu}\cap K_{{\cal F}}^{\nu}\subset K_{{\cal F}'}^{\nu}\subset K_{{\cal F}}^{\nu}$;   $M_{{\cal F}'_{\star},\nu}\subset L'$ et $M_{{\cal F}',\nu}\subset L'$.

On consid\`ere le cas de ${\cal F}'$, celui de ${\cal F}'_{\star}$ \'etant analogue. Soit $g\in K_{\star}^{\nu}\cap K_{{\cal F}}^{\nu}$. Alors l'action de $g$ fixe $x_{\star}$ et $x$. Elle fixe donc tout point du segment. En particulier, elle fixe un point de ${\cal F}'$ et $g$ appartient \`a $K_{{\cal F}'}^{\dag}$. Puisque $w_{G}(g)=\nu$, on a $\nu\in {\cal N}({\cal F}')$ et $g\in K_{{\cal F}'}^{\nu}$. L'inclusion $K_{{\cal F}'}^{\nu}\subset K_{{\cal F}}^{\nu}$ est \'evidente si ${\cal F}'={\cal F}$. Sinon elle r\'esulte de 5(6) car  ${\cal F}$ est incluse dans l'adh\'erence de ${\cal F}'$. L'intersection de ${\cal F}'$ et du segment $[x_{\star},x]$ est contenue dans l'ensemble des points fixes d'un \'el\'ement $g$ comme ci-dessus, donc incluse dans ${\cal F}^{_{'}\nu}$. Puisque cette intersection est ouverte dans le segment, une racine qui est constante sur ${\cal F}^{_{'}\nu}$ annule $v$. D'apr\`es le choix du point $x$, cela entra\^{\i}ne $\alpha\in \Sigma^{M_{{\cal F},\nu}}$. Cela entra\^{\i}ne aussi $\alpha\in \Sigma^L$ puisque, pour une racine $\beta\not\in \Sigma^L$, on a $\vert \beta(v)\vert >N-c$. Donc $\alpha\in \Sigma^{L'}$, d'o\`u l'inclusion $M_{{\cal F}',\nu}\subset L'$. Cela d\'emontre (7).

Une cons\'equence de ces propri\'et\'es est que
 $$K_{\star}^{\nu}\cap K_{{\cal F}}^{\nu} =K_{{\cal F}'_{\star}}^{\nu}\cap K_{{\cal F}'}^{\nu}.$$
 Les groupes qui conservent ces espaces sont donc aussi \'egaux, c'est-\`a-dire
   $$K_{\star}^0\cap K_{{\cal F}}^0=K_{{\cal F}'_{\star}}^{0}\cap K_{{\cal F}'}^{0}.$$

 On choisit un sous-groupe parabolique $P'$ de $G$ contenant $M_{min}$ et tel que $P'\cap M_{{\cal F},\nu}=Q'$. L'ensemble $\Sigma(U_{P'})$ se d\'ecompose en $\Sigma_{1}^+\cup \Sigma_{2}^+$, o\`u $\Sigma_{1}^+=\Sigma(U_{Q'})$ et $\Sigma_{2}^+$ est son compl\'ementaire. On ordonne $\Sigma(U_{P'})$ de sorte que les \'el\'ements de $\Sigma_{1}^+$ soient inf\'erieurs \`a ceux de $\Sigma_{2}^+$. On ordonne $\Sigma(U_{\bar{P}'})=-\Sigma(U_{P'})$ de sorte que les \'el\'ements de $-\Sigma_{2}^+$ soient inf\'erieurs \`a ceux de $-\Sigma_{1}^+$. 
 On applique le lemme 6 \`a $({\cal F}',\nu)$ et au couple $(P',L')$. C'est loisible puisque $M_{{\cal F}',\nu}\subset L'$. On a donc
$$K_{{\cal F}'}^{\nu}=\left(\prod_{\alpha\in \Sigma(U_{\bar{P}'})}K_{{\cal F}'}^0\cap U_{\alpha}(F)\right)(K_{{\cal F}'}^{\nu}\cap L'(F))\left(\prod_{\alpha\in \Sigma(U_{P'})}K_{{\cal F}'}^0\cap U_{-\alpha}(F)\right).$$
Une d\'ecomposition analogue vaut pour $K_{{\cal F}'_{\star}}^{\nu}$. Puisque les produits ci-dessus sont directs, on obtient une d\'ecomposition
$$ K_{\star}^{\nu}\cap K_{{\cal F}}^{\nu}=K_{{\cal F}'_{\star}}^{\nu}\cap K_{{\cal F}'}^{\nu}=\left(\prod_{\alpha\in \Sigma(U_{\bar{P}'})}K_{\alpha}\right)K_{L'}^{\nu}\left(\prod_{\alpha\in \Sigma(U_{P'})}K_{\alpha}\right),$$
o\`u $K_{L'}^{\nu}=K_{{\cal F}'}^{\nu}\cap K_{{\cal F}'_{\star}}^{\nu}\cap L'(F)$ et, pour tout $\alpha\in \Sigma(U_{P'})$, 
$$K_{\pm \alpha}=K_{{\cal F}'_{\star}}^0\cap K_{{\cal F}'}^0\cap U_{\pm \alpha}(F)=K_{\star}^0\cap K_{{\cal F}}^0\cap U_{\pm \alpha}(F).$$
  Cette d\'ecomposition se r\'ecrit
 $$ K_{\star}^{\nu}\cap K_{{\cal F}}^{\nu}=K_{2}^-K_{1}^-K_{L'}^{\nu}K_{1}^+K_{2}^+,$$
 o\`u, pour $i=1,2$ et $\epsilon=\pm$, $K_{i}^{\epsilon}$ est le produit des $K_{\epsilon \alpha}$ pour $\alpha\in \Sigma_{i}^+$.   On peut  remplacer notre int\'egration par une int\'egration sur le membre de droite ci-dessus (en choisissant convenablement une d\'ecomposition de la mesure). D'o\`u
  $$\int_{K_{{\cal F}}^{\nu}} f_{\star\star}(k) \phi_{{\cal F},\nu,cusp}(k)\,dk=\int_{K_{2}^-...K_{2}^+}f_{\star\star}(k_{2}^-k_{1}^-k_{L'} k_{1}^+k_{2}^+)\phi_{{\cal F},\nu,cusp}(k_{2}^-k_{1}^-k_{L'}k_{1}^+k_{2}^+)\,dk_{2}^-...dk_{2}^+.$$
  Pour $\alpha\in \Sigma_{2}^+$, $\alpha$ n'intervient pas dans $M_{{\cal F},\nu}$ et  il r\'esulte du lemme 6 que $K_{{\cal F}}^0\cap U_{\pm\alpha}(F)=K_{{\cal F}}^+\cap U_{\pm \alpha}(F)$. Donc $K_{\pm \alpha}\subset K_{{\cal F}}^+$. La fonction $\phi_{{\cal F},\nu,cusp}$ est invariante par ce groupe, donc, dans la formule ci-dessus, on a
 $$\phi_{{\cal F},\nu,cusp}(k_{2}^-k_{1}^- k_{L'}k_{1}^+k_{2}^+)=\phi_{{\cal F},\nu,cusp}(k_{1}^- k_{L'}k_{1}^+).$$ 
 Pour $\alpha\in \Sigma_{1}^+$, $\alpha$ intervient dans $U_{Q'}$ donc dans $U_{Q}$. Donc $\alpha(v)>N-c$. On a $\alpha(v)=\alpha(x)-\alpha(x_{\star})<c_{\alpha,{\cal F}}^+-c_{\alpha,{\cal F}_{\star}}$. Gr\^ace \`a (4), on en d\'eduit que $c_{\alpha,{\cal F}}^->c_{\alpha,{\cal F}_{\star}}$. Gr\^ace \`a 4(6) et 4(7), cela entra\^{\i}ne que
  $K_{\star}^0\cap U_{\alpha}(F)\subset K_{{\cal F}}^+$, d'o\`u aussi $K_{\alpha}\subset K_{{\cal F}}^+$. Pour la m\^eme raison que ci-dessus, on obtient
  $$\phi_{{\cal F},\nu,cusp}(k_{1}^- k_{L'}k_{1}^+)=\phi_{{\cal F},\nu,cusp}(k_{1}^- k_{L'}).$$
  Toujours pour $\alpha\in -\Sigma_{1}^+$, $\alpha$ intervient dans $U_{\bar{Q}}$ donc $\alpha(v)<-N+c$. Puisque $\alpha(v)=\alpha(x)-\alpha(x_{\star})>c_{\alpha,{\cal F}}-c_{\alpha,{\cal F}_{\star}}^+$, on d\'eduit de (4)  que $c_{\alpha,{\cal F}}<-c$. Gr\^ace \`a 4(6) et (3) ci-dessus, cela entra\^{\i}ne 
   que $K_{{\cal F}}^0\cap U_{\alpha}(F)\subset H$. Donc $K_{\alpha}=K_{{\cal F}}^0\cap U_{\alpha}(F)\subset H$. Puisque $H$ est distingu\'e dans $K_{\star}^{\dag}$ et que $f_{\star\star}$ est invariante par ce groupe, cela entra\^{\i}ne d'abord l'\'egalit\'e
   $$f_{\star\star}(k_{2}^-k_{1}^- k_{L'}k_{1}^+k_{2}^+)=f_{\star\star}(k_{2}^- k_{L'}k_{1}^+k_{2}^+).$$
    Cela entra\^{\i}ne aussi  aussi que $K_{1}^-=K_{{\cal F}}^0\cap   U_{\bar{Q}'}(F)$. 
   A ce point, on obtient la formule
 $$(8) \qquad \int_{K_{{\cal F}}^{\nu}} f_{\star\star}(k) \phi_{{\cal F},\nu,cusp}(k)\,dk=$$
 $$\int_{K_{2}^-K_{L'}K_{1}^+K_{2}^+}  f_{\star\star}(k_{2}^- k_{L'}k_{1}^+k_{2}^+)\phi_{{\cal F},\nu,cusp,\bar{Q}'}(k_{L'})\, dk_{2}^-dk_{L'}dk_{1}^+dk_{2}^+,$$
 o\`u
$$\phi_{{\cal F},\nu,cusp,\bar{Q}'}(k_{L'})=\int_{K_{{\cal F}}^0\cap   U_{\bar{Q}'}(F)}\phi_{{\cal F},\nu,cusp}(\bar{u}k_{L'})\,d\bar{u}.$$
L'image de $K_{{\cal F}}^0\cap \bar{Q}(F)$ dans ${\bf G}_{{\cal F}}({\mathbb F}_{q})$ est le groupe des points sur ${\mathbb F}_{q}$ d'un sous-groupe parabolique de ${\bf \bar{Q}}'$ de   ${\bf G}_{{\cal F}}$ et l'image de $K_{{\cal F}}^0\cap L(F)$ est le groupe des points sur ${\mathbb F}_{q}$ d'une composante de Levi ${\bf L}'$  de ${\bf \bar{Q}}'$. D'apr\`es le lemme 6 appliqu\'e \`a ${\cal F}'$ et au groupe $L'$, on peut fixer $n\in Norm_{L'(F)}(A)\cap K_{{\cal F}'}^{\nu}$. A fortiori $n\in K_{{\cal F}}^{\nu}$.  Notons ${\bf n}$  l'image de $n$  dans ${\bf G}_{{\cal F}}^{\nu}$. L'image $w$ de $n$ dans $W$ appartient \`a   $W^{L'}$. La compatibilit\'e des actions de $n$ et de ${\bf n}$ entra\^{\i}ne alors que la conjugaison par ${\bf n}$   conserve ${\bf \bar{Q}}'$ et ${\bf L}'$. Alors ${\bf \bar{Q}}^{_{'}\nu}={\bf n}{\bf \bar{Q}}'$ est un sous-espace parabolique de ${\bf G}_{{\cal F}}^{\nu}$ d'espace de Levi ${\bf L}^{_{'}\nu}={\bf n}{\bf L}'$. Introduisons le sous-tore ${\bf A}$ de ${\bf G}_{{\cal F}}$ comme en 4(8) et (9). Identifions ${\cal A}$ \`a $X_{*}({\bf A})\otimes_{{\mathbb Z}}{\mathbb R}$. Au groupe ${\bf L}'$ et \`a l'espace ${\bf L}^{_{'}\nu}$ correspondent des sous-espaces ${\cal A}_{{\bf L}'}$ et ${\cal A}_{{\bf L}'}^{\nu}$. Ce dernier n'est autre que le sous-espace des points fixes de $w$ agissant dans le premier. Or, puisque $L'=M_{{\cal F},\nu}\cap L$, ${\cal A}_{{\bf L}'}$ contient le sous-espace ${\cal A}_{L'}$. Puisque $w\in W^{L'}$, il agit trivialement sur cet espace donc ${\cal A}_{{\bf L}'}^{\nu}$ contient lui aussi ${\cal A}_{L'}$. Puisque $L'$ ne contient pas $M_{{\cal F},\nu}$, cela entra\^{\i}ne ${\cal A}_{{\bf L}'}^{\nu}\not={\cal A}_{{\cal F}}^{\nu}$. Donc  ${\bf L}^{_{'}\nu}$ est un sous-espace de Levi  propre de ${\bf G}_{{\cal F}}^{\nu}$. 
En consid\'erant les fonctions $\phi_{{\cal F},\nu,cusp}$ et $\phi_{{\cal F},\nu,cusp,\bar{Q}'}$ comme des fonctions sur ${\bf G}_{{\cal F}}^{\nu}({\mathbb F}_{q})$ et ${\bf L}^{_{'}\nu}({\mathbb F}_{q})$, on voit que $\phi_{{\cal F},\nu,cusp,\bar{Q}'}$  est \'egale \`a $res_{{\bf L}^{_{'}\nu}}^{{\bf G}_{{\cal F}}^{\nu}}(\phi_{{\cal F},\nu,cusp})$, \`a une constante non nulle pr\`es provenant des choix de mesures. Mais $\phi_{{\cal F},\nu,cusp}$ est cuspidale et le sous-espace  ${\bf L}^{_{'}\nu}$ est propre. Donc $ \phi_{{\cal F},\nu,cusp,\bar{Q}'}=0$. Alors l'\'egalit\'e (8) d\'emontre le (i) de l'\'enonc\'e. $\square$

Consid\'erons le cas o\`u   $M$ contient $M_{{\cal F},\nu}$. On projette l'appartement $App(A)$ sur $App^M(A)$. Comme dans la preuve du (ii), on peut construire un segment $[x_{\star},x]$ et des facettes ${\cal F}'_{\star}$ et ${\cal F}'$.   
  Comme dans cette preuve, on obtient  la d\'ecomposition
   $$K_{\star}^0\cap K_{{\cal F}}^0=(K_{{\cal F}}^+\cap U_{\bar{P}}(F))(K_{\star}^0\cap K_{{\cal F}}^0\cap M(F))(K_{\star}^0\cap U_{P}(F))$$
   et l'\'egalit\'e
$$K_{\star}^0\cap K_{{\cal F}}^0\cap M(F)=K_{{\cal F}'_{\star}}^0\cap K_{{\cal F}'}^0\cap M(F).$$
Les facettes ${\cal F}_{\star}$, ${\cal F}'_{\star}$, ${\cal F}' $ et ${\cal F}$ se projettent dans des facettes ${\cal F}_{\star}^M$, ${\cal F}_{\star}^{_{'}M}$, ${\cal F}^{_{'}M}$ et ${\cal F}^M$ de $App^M(A)$.  Notons $x_{\star}^M$ et $x^M$ les projections de $x_{\star}$ et $x$. Ces points et les facettes pr\'ec\'edentes sont dans la m\^eme situation que $x_{\star}$, $x$ et les facettes  ${\cal F}_{\star}$, ${\cal F}'_{\star}$, ${\cal F}' $ et ${\cal F}$. C'est-\`a-dire que $x^M_{\star}\in {\cal F}_{\star}^M$, $x^M\in {\cal F}^M$,  les facettes ${\cal F}^{_{'}M}_{\star}$ et ${\cal F}^{_{'}M}$ coupent le segment $[x_{\star}^M,x^M]$ selon des ensembles ouverts, $x_{\star}^M$ est adh\'erent \`a  ${\cal F}^{_{'}M}_{\star}$ et $x^M$ est adh\'erent \`a ${\cal F}^{_{'}M}$. Il en r\'esulte que 
$$K_{{\cal F}^{_{'}M}_{\star}}^0\cap K_{{\cal F}^{_{'}M}}^0=K_{\star}^{M,0}\cap K_{{\cal F}^M}^0,$$
o\`u on a pos\'e $K_{\star}^{M,0}=K_{{\cal F}_{\star}^M}^0$. Or $M$ contient $M_{{\cal F}_{\star}^{_{'}}}$ et 
  $M_{{\cal F}^{_{'}}}$.  D'apr\`es le lemme 6, on a donc $K_{{\cal F}^{_{'}M}_{\star}}^0=K_{{\cal F}'_{\star}}^0\cap M(F)$ et $K_{{\cal F}^{_{'}M}}^0=K_{{\cal F}'}^0\cap M(F)$. Toutes ces relations conduisent \`a l'\'egalit\'e
$$K_{\star}^0\cap K_{{\cal F}}^0\cap M(F)=  K_{\star}^{M,0}\cap K_{{\cal F}^M}^0$$
et la d\'ecomposition ci-dessus devient
$$ K_{\star}^0\cap K_{{\cal F}}^0=(K_{{\cal F}}^+\cap U_{\bar{P}}(F))( K_{\star}^{M,0}\cap K_{{\cal F}^M}^0)(K_{\star}^0\cap U_{P}(F)).$$
On peut appliquer les m\^emes constructions en rempla\c{c}ant la facette ${\cal F}_{\star}$ par ${\cal F}_{\star\star}$ (a priori, il pourrait \^etre n\'ecessaire d'accro\^{\i}tre $N$, ce qui serait sans cons\'equence, mais on v\'erifie que ce n'est pas m\^eme n\'ecessaire). D'o\`u
$$(9) \qquad  K_{\star\star}^0\cap K_{{\cal F}}^0=(K_{{\cal F}}^+\cap U_{\bar{P}}(F))( K_{\star\star}^{M,0}\cap K_{{\cal F}^M}^0)(K_{\star\star}^0\cap U_{P}(F)).$$

 \bigskip
 
 \section{Un r\'esultat d'annulation}
 Fixons $N$ tel que le lemme pr\'ec\'edent soit v\'erifi\'e pour tout $P\in {\cal F}(M_{min})$.  Pour tout tel $P$ et pour tout $f\in C_{c}(K_{\star}^{\dag}/H)$,  d\'efinissons
 $$B_{N,R}(P,f)=\sum_{({\cal F},\nu)\in Fac^*(G;A); {\cal F}\subset X_{N}(P)}  mes(K_{\star\star}^0\cap K_{{\cal F}}^0)^{-1}z_{R}({\cal F},\nu)\int_{K_{{\cal F}}^{\nu}}f_{\star\star}(k)\phi_{{\cal F},\nu,cusp}(k)\,dk.$$
 D'apr\`es 15(1), la formule 13(4) se r\'ecrit
 $$(1) \qquad trace\,\pi(f)=\sum_{P\in {\cal F}(M_{min})}B_{N,R}(P,f)$$
 pourvu que $R$ soit assez grand. 
 
 Fixons $P$ et notons $M$ sa composante de Levi contenant $M_{min}$. D'apr\`es le (ii) du lemme 15, ne contribuent \`a $B_{N,R}(P,f)$ que des \'el\'ements $({\cal F},\nu)$ tels que $M$ contient $M_{{\cal F},\nu}$. Notons $Y_{N}(P)$ l'ensemble des $({\cal F},\nu)\in Fac^*(G;A)$ tels que ${\cal F}\subset X_{N}(P)$ et $M$ contient $M_{{\cal F},\nu}$.   Pour $({\cal F},\nu)\in Y_{N}(P)$, le lemme 15(i) calcule l'int\'egrale intervenant dans $B_{N,R}(P,f)$. La relation 15(9) conduit \`a l'\'egalit\'e
 $$mes(K_{\star\star}^0\cap K_{{\cal F}}^0)^{-1}mes(K_{{\cal F}}^+\cap U_{\bar{P}}(F))=mes(K_{\star\star}^0\cap U_{P}(F))^{-1}mes(K_{\star\star}^{M,0}\cap K_{{\cal F}^M}^0)^{-1}.$$
  Ainsi
 $$B_{N,R}(P,f)=mes(K_{\star\star}^0\cap U_{P}(F))^{-1}\sum_{({\cal F},\nu)\in Y_{N}(P)}mes(K_{\star\star}^{M,0}\cap K_{{\cal F}^M}^0)^{-1}z_{R}({\cal F},\nu)$$
 $$\int_{K_{{\cal F}^M}^{\nu}}f_{\star\star,[P]}(k)\phi_{{\cal F}^M,\nu,cusp}(k)\,dk.$$
 Pour tout $({\cal F}_{M},\nu)\in Fac^*(M;A)_{G-comp}$, notons $Y_{N}(P,{\cal F}_{M},\nu)$ l'ensemble des ${\cal F}\in Fac(G;A)$ telles que $({\cal F},\nu)\in Y_{N}(P)$ et ${\cal F}^M={\cal F}_{M}$. La formule ci-dessus se r\'ecrit
 $$(2) \qquad B_{N,R}(P,f)=mes(K_{\star\star}^0\cap U_{P}(F))^{-1}\sum_{({\cal F}_{M},\nu)\in Fac^*(M;A)_{G-comp}}mes(K_{\star\star}^{M,0}\cap K_{{\cal F}_{M}}^0)^{-1}z_{N,R}(P,{\cal F}_{M},\nu)$$
 $$\int_{K_{{\cal F}^M}^{\nu}}f_{\star\star,[P]}(k)\phi_{{\cal F}_{M},\nu,cusp}(k)\,dk,$$
 o\`u
 $$(3) \qquad z_{N,R}(P,{\cal F}_{M},\nu)=\sum_{{\cal F}\in Y_{N}(P,{\cal F}_{M},\nu)}z_{R}({\cal F},\nu).$$
 
   De m\^eme que l'on a d\'efini les sous-ensembles $X'_{N}(Q)$ et $X_{N}(Q)$ de $App(A)$ pour $Q\in {\cal F}(M_{min})$, on d\'efinit les sous-ensembles $X^{_{'}M}(Q)$ et $X^M(Q)$ de $App^M(A)$ pour $Q\in {\cal F}^M(M_{min})$. Notons $B_{M}$ la boule dans ${\cal A}_{M}$ de centre $0$ et de rayon $1$. Le sous-groupe parabolique $P$ d\'efinit une chambre $C_{P}$ dans ${\cal A}_{M}$. 
   
   \ass{Proposition}{Si $R$ est assez grand relativement \`a $N$, les propri\'et\'es suivantes sont v\'erifi\'ees pour tout $({\cal F}_{M},\nu)\in Fac^*(M;A)_{G-comp}$:
   
   (i) si ${\cal F}_{M}$ est contenu dans  $X^M_{N}(M)$ et ${\cal F}_{M}^{\nu}$ est r\'eduit \`a un point, $z_{N,R}(P,{\cal F}_{M},\nu)=mes(B_{M}\cap C_{P})mes(B_{M})^{-1}$;
   
   (ii) sinon, $z_{N,R}(P,{\cal F}_{M},\nu)=0$.}
   
   Preuve. D'apr\`es la d\'efinition (3) ci-dessus et celle de $z_{R}({\cal F},\nu)$ donn\'ee en 13, on a
  $$ z_{N,R}(P,{\cal F}_{M},\nu)=\sum_{{\cal F}\in Y_{N}(P,{\cal F}_{M},\nu)}\sum_{{\cal F}_{1}\in E_{R}({\cal F},\nu)} (-1)^{dim({\cal F}_{1}^{\nu})}z({\cal F}_{1},{\cal F},\nu).$$
  Pour simplifier la notation, on change ${\cal F}$ en ${\cal F}'$ et ${\cal F}_{1}$ en ${\cal F}$. En se rappelant les d\'efinitions, on obtient
  $$(4) \qquad z_{N,R}(P,{\cal F}_{M},\nu)=\sum_{{\cal F}\subset B_{R}(A); \nu\in {\cal N}({\cal F})}(-1)^{dim({\cal F}^{\nu})}\sum_{{\cal F}'\in Y_{N}(P,{\cal F}_{M},{\cal F},\nu)}z({\cal F},{\cal F}',\nu),$$
  o\`u on rappelle que $B_{R}(A)=B_{R}\cap App(A)$ et o\`u on a not\'e $Y_{N}(P,{\cal F}_{M},{\cal F},\nu)$ l'ensemble des ${\cal F}'\in Y_{N}(P,{\cal F}_{M},\nu)$ telles que ${\cal F}\subset \bar{{\cal F}}'$.
  Plus explicitement, $Y_{N}(P,{\cal F}_{M},{\cal F},\nu)$ est l'ensemble des facettes ${\cal F}'\in Fac(G;A)$ telles que $\nu\in {\cal N}({\cal F}')$, ${\cal F}'\subset X_{N}(P)$, $M$ contient $M_{{\cal F}',\nu}$, ${\cal F}^{_{'}M}={\cal F}_{M}$ et ${\cal F}\subset \bar{{\cal F}}'$.
  
  Notons $\Psi$ l'ensemble des facettes ${\cal F}\in Fac(G;A)$ qui coupent l'ensemble $p_{M}^{-1}(\bar{{\cal F}}_{M}^{\nu})$. D'apr\`es le lemme 7(i), pour toute ${\cal F}\in \Psi$, on a $\nu\in {\cal N}({\cal F})$ et ${\cal F}^{\nu}={\cal F}\cap p_{M}^{-1}(\bar{{\cal F}}_{M}^{\nu})$.  L'ensemble $p_{M}^{-1}(\bar{{\cal F}}_{M}^{\nu})$ est r\'eunion disjointe des ${\cal F}^{\nu}$ pour ${\cal F}\in \Psi$ et cette d\'ecomposition est une d\'ecomposition en facettes qui sont toutes des polysimplexes comme on l'a vu au paragraphe 5. Notons $\Psi_{R}$ l'ensemble des ${\cal F}\in \Psi$ telles que ${\cal F}\subset B_{R}(A)$. Puisque  $B_{R}$ est une r\'eunion de facettes, pour ${\cal F}\in Fac(G;A)$, on a ${\cal F}\in \Psi_{R}$ si et seulement si ${\cal F}$ coupe l'ensemble $p_{M}^{-1}(\bar{{\cal F}}_{M}^{\nu})\cap B_{R}(A)$.  L'ensemble $p_{M}^{-1}(\bar{{\cal F}}_{M}^{\nu})\cap B_{R}(A)$ est r\'eunion disjointe des ${\cal F}^{\nu}$ pour ${\cal F}\in \Psi_{R}$.

  Consid\'erons une facette ${\cal F}\in Fac(G;A)$ telle que $\nu\in {\cal N}({\cal F})$ et que l'ensemble $Y_{N}(P,{\cal F}_{M},{\cal F},\nu)$  ne soit pas vide. Pour ${\cal F}'\in Y_{N}(P,{\cal F}_{M},{\cal F},\nu)$, on a ${\cal F}^{_{'}\nu}={\cal F}'\cap p_{M}^{-1}({\cal F}_{M}^{\nu})$ d'apr\`es le lemme 6, donc $\bar{{\cal F}}^{_{'}\nu}\subset p_{M}^{-1}(\bar{{\cal F}}_{M}^{\nu})$. Les conditions $\nu\in {\cal N}({\cal F})$ et ${\cal F}\subset \bar{{\cal F}}'$ entra\^{\i}nent ${\cal F}^{\nu}\subset \bar{{\cal F}}^{_{'}\nu}$ d'apr\`es 5(2). Donc ${\cal F}\in \Psi$.  Si de plus ${\cal F}\subset B_{R}(A)$, on a ${\cal F}\in \Psi_{R}$. Dans l'expression (4), on peut donc  remplacer  
   l'ensemble de sommation en ${\cal F}$ par   $\Psi_{R}$. 
  
   Consid\'erons ${\cal F}\in \Psi_{R}$ et ${\cal F}'\in  Y_{N}(P,{\cal F}_{M},{\cal F},\nu)$. Les conditions $\nu\in {\cal N}({\cal F}')$,   $M$ contient $M_{{\cal F}',\nu}$ et ${\cal F}^{_{'}M}={\cal F}_{M}$ entra\^{\i}nent ${\cal A}_{{\cal F}'}^{\nu}={\cal A}_{{\cal F}_{M}}^{\nu}$ d'apr\`es le lemme 6. Rappelons qu'en posant $L=M_{{\cal F}_{M},\nu}$,  cet espace n'est autre que ${\cal A}_{L}$ comme on l'a vu au paragraphe 5. Au paragraphe 13, on a d\'efini  une chambre $C_{{\cal F}'}^{{\cal F},\nu}$ dans cet espace. Elle est d\'efinie par des in\'egalit\'es $\alpha(x)>0$ pour un certain ensemble de $\alpha\in \Sigma-\Sigma^L$. Il est bien connu que les annulateurs dans ${\cal A}_{L}$ des racines dans $\Sigma-\Sigma^L$ d\'ecoupent cet espace en chambres (minimales) associ\'ees aux sous-groupes paraboliques de composantes de Levi $L$. Pour tout tel sous-groupe $Q$, notons $C_{Q}$ la chambre associ\'ee. Alors la chambre  $C_{{\cal F}'}^{{\cal F},\nu}$ est r\'eunion disjointe  d'un certain nombre de  chambres $C_{Q}$ et de sous-vari\'et\'es ferm\'ees de mesures nulles. Rappelons que $z({\cal F},{\cal F}',\nu)=mes(B_{L}\cap C_{{\cal F}'}^{{\cal F},\nu})mes(B_{L})^{-1}$ o\`u $B_{L}$ est la boule dans ${\cal A}_{L}$ de centre $0$ et de rayon $1$. On a donc
   $$z({\cal F},{\cal F}',\nu)=\sum_{Q}mes(B_{L}\cap C_{Q})mes(B_{L})^{-1},$$
   o\`u $Q$ parcourt les \'el\'ements de ${\cal P}(L)$ tels que $C_{Q}\subset  C_{{\cal F}'}^{{\cal F},\nu}$. Inversement, partons de ${\cal F}\in \Psi_{R}$ et d'un sous-groupe parabolique $Q\in {\cal P}(L)$. Il y a au plus une facette ${\cal F}'\in Y_{N}(P,{\cal F}_{M},{\cal F},\nu)$ telle que $C_{Q}$ soit contenue dans $C_{{\cal F}'}^{{\cal F},\nu}$. En effet, posons $\underline{L}=M_{{\cal F},\nu}$ et  identifions ${\cal A}_{{\cal F}}^{\nu}$ \`a ${\cal A}_{\underline{L}}$. La condition ${\cal F}^{\nu}\subset p_{M}^{-1}(\bar{{\cal F}}_{M}^{\nu})$ entra\^{\i}ne que $L\subset \underline{L}$. La chambre $C_{{\cal F}'}^{{\cal F},\nu}$ associ\'ee \`a une facette ${\cal F}'\in Y_{N}(P,{\cal F}_{M},{\cal F},\nu)$ est d\'etermin\'ee par un sous-groupe parabolique $Q^{\underline{L}}\in {\cal P}^{\underline{L}}(L)$. La condition  $C_{Q}\subset  C_{{\cal F}'}^{{\cal F},\nu}$ \'equivaut \`a ce que $Q\cap \underline{L}=Q^{\underline{L}}$. Plus concr\`etement, fixons un \'el\'ement $\rho_{Q}\in C_{Q}$ en position g\'en\'erale. Il existe une unique facette telle que, pour $x\in {\cal F}^{\nu}$, il existe $\epsilon>0$ tel que le segment $]x,x+\epsilon\rho_{Q}]$ soit  contenu dans cette facette. Notons-la  ${\cal F}[Q]$. Alors  l'ensemble des ${\cal F}' \in Y_{N}(P,{\cal F}_{M},{\cal F},\nu)$  telles que $C_{Q}$ soit contenue dans $C_{{\cal F}'}^{{\cal F},\nu}$ est soit vide, soit r\'eduit \`a cette unique facette ${\cal F}[Q]$. Notons $\Psi_{R}[Q]$ l'ensemble des ${\cal F}\in \Psi_{R}$ telles que ${\cal F}[Q]$ appartient \`a    $Y_{N}(P,{\cal F}_{M},{\cal F},\nu)$.
    L'expression (4) se r\'ecrit
   
$$(5) \qquad z_{N,R}(P,{\cal F}_{M},\nu)=\sum_{Q\in {\cal P}(L)}mes(B_{L}\cap C_{Q})mes(B_{L})^{-1}Z_{N,R}[Q],$$
o\`u
$$Z_{N,R}[Q]=\sum_{{\cal F}\in \Psi_{R}[Q]}(-1)^{dim({\cal F}^{\nu})}.$$

Fixons $Q\in {\cal P}(L)$ et ${\cal F}\in \Psi_{R}$. Montrons que

(6) ${\cal F}\in \Psi_{R}[Q]$ si et seulement si pour $x\in {\cal F}^{\nu}$, il existe $\epsilon>0$ de sorte que le segment $]x,x+\epsilon\rho_{Q}]$ soit contenu dans $p_{M}^{-1}(\bar{{\cal F}}_{M}^{\nu})$ et dans $X_{N}(P)$. 

En effet, ce segment est contenu dans ${\cal F}[Q]$. Si ${\cal F}\in \Psi_{R}[Q]$,  on a ${\cal F}[Q]\in Y_{N}(P,{\cal F}_{M},{\cal F},\nu)$. Donc ${\cal F}[Q]\subset X_{N}(P)$ et le segment lui-m\^eme est inclus dans $X_{N}(P)$. D'autre part, ${\cal F}[Q]\subset p_{M}^{-1}({\cal F}_{M})$, donc $]x,x+\epsilon\rho_{Q}]$ est inclus dans cet ensemble. Sa projection $]p_{M}(x),p_{M}(x)+\epsilon p_{M}(\rho_{Q})]$ est  incluse dans ${\cal F}_{M}$ mais aussi dans $p_{M}(x)+{\cal A}_{L}/{\cal A}_{M}$ puisque $\rho_{Q}\in {\cal A}_{L}$. Or, puisque ${\cal A}_{L}={\cal A}_{M_{{\cal F}_{M},\nu}}$ et que $p_{M}(x)$ appartient \`a $\bar{{\cal F}}_{M}^{\nu}$, l'intersection de ${\cal F}$ et de cet espace affine est \'egale \`a ${\cal F}_{M}^{\nu}$. Donc $]x,x+\epsilon\rho_{Q}]$ est contenu dans $p_{M}^{-1}(\bar{{\cal F}}_{M}^{\nu})$. Inversement, supposons la conclusion de (6) v\'erifi\'ee.  Alors la facette ${\cal F}[Q]$ contenant notre segment coupe $p_{M}^{-1}(\bar{{\cal F}}_{M}^{\nu})$.  Parce que l'on suppose $\rho_{Q}$ en position g\'en\'erale, la facette ${\cal F}[Q]$ coupe  cet ensemble  selon un ouvert. D'apr\`es le (ii) du lemme 7, cela entra\^{\i}ne que $\nu\in {\cal N}({\cal F}[Q])$, $M$ contient $M_{{\cal F}[Q],\nu}$ et ${\cal F}[Q]^M={\cal F}_{M}$. Il est clair par construction que $x$ appartient \`a $\bar{{\cal F}}[Q]$, donc ${\cal F}$ tout enti\`ere est contenue dans $\bar{{\cal F}}[Q]$. Enfin ${\cal F}[Q]$ coupe $X_{N}(P)$ et est donc incluse dans cet ensemble. Toutes les conditions sont v\'erifi\'ees pour que ${\cal F}[Q]$ appartienne \`a $Y_{N}(P,{\cal F}_{M},{\cal F},\nu)$, autrement dit pour que ${\cal F}$ appartienne \`a $\Psi_{R}[Q]$. Cela prouve (6).

Pour $Q$ et ${\cal F}$ comme ci-dessus, \'etudions la condition: $]x,x+\epsilon\rho_{Q}]$ est contenu dans $p_{M}^{-1}(\bar{{\cal F}}_{M}^{\nu})$, ce qui \'equivaut \`a ce que le segment ferm\'e $[x,x+\epsilon\rho_{Q}]$ soit contenu dans le m\^eme ensemble. En posant $x^M=p_{M}(x)$ et $\rho_{Q}^M=p_{M}(\rho_{Q})$, la condition \'equivaut \`a ce que  $[x^M,x^M+\epsilon\rho_{Q}^M]$ soit contenu dans $\bar{{\cal F}}_{M}^{\nu}$. Remarquons que $x^M$ appartient \`a la facette ${\cal F}^M$ qui contient $p_{M}({\cal F})$, plus pr\'ecis\'ement $x^M$ appartient au sous-ensemble ${\cal F}^{M,\nu}$ qui est une facette de $\bar{{\cal F}}_{M}^{\nu}$. 
Supposons d'abord que $\Sigma^M$ est un syst\`eme de racines irr\'eductible et utilisons la description de l'ensemble $\bar{{\cal F}}_{M}^{\nu}$ donn\'ee au paragraphe 5. On y supprime les exposants $\nu$ pour simplifier la notation. L'ensemble $\bar{{\cal F}}_{M}^{\nu}$ est l'enveloppe convexe d'une famille de points $(s_{i})_{i\in I}$.  Le point $x^M$ appartient \`a la facette  ${\cal F}^{M,\nu}$ de cette enveloppe, \`a laquelle est associ\'e  un sous-ensemble non vide $I_{{\cal F}^{M,\nu}}\subset I$. En fixant un point $s\in \bar{{\cal F}}_{M}^{\nu}$, on peut \'ecrire $x^M-s=\sum_{i\in I}a_{i}(s_{i}-s)$ avec des $a_{i}>0$  pour $i\in I_{{\cal F}^{M,\nu}}$, $a_{i}=0$ pour $i\in I-I_{{\cal F}^{M,\nu}}$ et  $\sum_{i\in I}a_{i}=1$.    Tout \'el\'ement de ${\cal A}_{L}/{\cal A}_{M}$ peut s'\'ecrire de fa\c{c}on unique $\sum_{i\in I}\lambda_{i}(s_{i}-s)$ avec des r\'eels $\lambda_{i}$ tels que $\sum_{i\in I}\lambda_{i}=0$. En particulier, on peut \'ecrire ainsi $\rho_{Q}^M=\sum_{i\in I}\lambda_{i}(s_{i}-s)$.   Supposons d'abord $\vert I\vert \geq 2$. 
Parce que $\rho_{Q}$ est en position g\'en\'erale, on a $\lambda_{i}\not=0$ pour tout $i\in I$. Notons $I_{Q}$ l'ensemble des $i\in I$ tels que $\lambda_{i}<0$. Cet ensemble est non vide puisque la somme des $\lambda_{i}$ est nulle. Remarquons que, pour la m\^eme raison, on a $I_{Q}\not=I$. On a $x^M+\epsilon\rho_{Q}^M-s=\sum_{i\in I}(a_{i}+\epsilon\lambda_{i})(s_{i}-s)$. La condition pour que $x^M+\epsilon\rho_{Q}$ appartienne \`a l'enveloppe convexe de la famille $(s_{i})_{i\in I}$ pour $\epsilon$ assez petit est donc que $a_{i}>0$ quand $\lambda_{i}<0$. Autrement dit $I_{Q}\subset I_{{\cal F}^{M,\nu}}$. Si maintenant $I$ n'a qu'un \'el\'ement, forc\'ement \'egal \`a $s$, on a $x^M=s$ et $\rho_{Q}^M=0$. Donc le segment $[x^M,x^M+\epsilon\rho_{Q}^M]$  est r\'eduit au point $x^M$ et est inclus dans $\bar{{\cal F}}_{M}^{\nu}$. En posant formellement $I_{Q}=I$, le r\'esultat est le m\^eme que dans le cas $\vert I\vert \geq2$. 
Il est facile d'\'etendre ces r\'esultats au cas g\'en\'eral o\`u $\Sigma^M$ n'est plus suppos\'e irr\'eductible. On d\'ecompose $\Sigma^M$ en $\Sigma^M_{1}\times...\times \Sigma^M_{k}$ o\`u les $\Sigma^M_{l}$ sont irr\'eductibles.  Soit $l\in \{1,...,k\}$. A ${\cal F}_{M}^{\nu}$, resp.  ${\cal F}$, resp. $\rho_{Q}$,  est associ\'e  une famille de points $(s_{i,l})_{i\in I_{l}}$, resp. un sous-ensemble non vide $I_{{\cal F}^{M,\nu},l}\subset I_{l}$, resp. un sous-ensemble non vide $I_{Q,l}\subset I_{l}$. On note $I$, resp. $I_{{\cal F}^{M,\nu}}$, $I_{Q}$, la r\'eunion disjointe des $I_{l}$, resp. $I_{{\cal F}^{M,\nu},l}$, $I_{Q,l}$. Alors 

(7) $]x,x+\epsilon\rho_{Q}]$ est contenu dans $p_{M}^{-1}(\bar{{\cal F}}_{M}^{\nu})$ si et seulement si $I_{Q}\subset I_{{\cal F}^{M,\nu}}$.

On conserve $Q$ et ${\cal F}$ et on suppose que le segment $]x,x+\epsilon\rho_{Q}]$ est contenu dans $p_{M}^{-1}(\bar{{\cal F}}_{M}^{\nu})$. On a d\'ej\`a remarqu\'e que, parce que l'on suppose $\rho_{Q}$ en position g\'en\'erale, cela entra\^{\i}ne qu'il est en fait contenu dans $p_{M}^{-1}({\cal F}_{M}^{\nu})$. Etudions  la seconde condition: $]x,x+\epsilon\rho_{Q}]$ est contenu dans $X_{N}(P)$.  Montrons d'abord

(8) si ${\cal F}_{M}$ n'est pas contenue dans $X_{N}^M(M)$, $]x,x+\epsilon\rho_{Q}]$ n'est pas  contenu dans $X_{N}(P)$.

Si ${\cal F}_{M}$ n'est pas contenue dans $X_{N}^M(M)$, on peut fixer un sous-groupe parabolique $P_{1}^M\in {\cal F}^M(M_{min})$ tel que $P_{1}^M\not=M$ et ${\cal F}_{M}$ est contenue dans $X_{N}^{_{'}M}(P_{1}^M)$. Notons $P_{1}$ le sous-groupe parabolique contenu dans $P$ et tel que $P_{1}\cap M=P_{1}^M$. Supposons que notre segment soit inclus dans $X_{N}(P)$. Soit $y$ un point  du segment et soit $\alpha\in \Sigma(U_{P_{1}})=\Sigma(U_{P})\cup \Sigma^M(U_{P_{1}^M})$. Si $\alpha\in \Sigma(U_{P})$, on a $\alpha(y)>N$ par d\'efinition de $X_{N}(P)$.  Si $\alpha\in \Sigma(U_{P_{1}^M})$, on a aussi $\alpha(y)>N$ parce que $p_{M}(y)\in {\cal F}$ et que ${\cal F}$ est contenu dans $X_{N}^{_{'}N}(P_{1}^M)$. Cela prouve que $y$ appartient \`a $X'_{N}(P_{1})$. Parce que $P_{1}\subsetneq P$, cela contredit l'hypoth\`ese que $y$ appartient \`a $X_{N}(P)$. On a prouv\'e (8).

Supposons que  ${\cal F}_{M}$ est  contenue dans $X_{N}^M(M)$. 
On peut fixer un sous-groupe parabolique minimal $P_{min}^M\in {\cal P}^M(M_{min})$ de $M$ de sorte que $\bar{{\cal F}}_{M}$ soit contenue dans l'adh\'erence de la chambre positive associ\'ee \`a ce parabolique. Autrement dit, en notant $\Delta^M$ la base de $\Sigma^M$ associ\'ee \`a $P_{min}$, on a $\alpha(y)\geq0$ pour tout $\alpha\in \Delta^M$ et tout $y\in \bar{{\cal F}}_{M}$. On note $P_{min}$ l'unique sous-groupe parabolique contenu dans $P$ et tel que $P_{min}\cap M=P_{min}^M$. On note $\Delta$ la base de $\Sigma$ associ\'ee \`a $P_{min}$. Soit $y\in p_{M}^{-1}({\cal F}_{M})$. Alors 

(9) la condition $y\in X_{N}(P)$ \'equivaut \`a $\alpha(y)>N$ pour tout $\alpha\in \Delta-\Delta^M$. 

En effet, la condition $y\in X_{N}(P)$ implique ces relations parce que $\Delta-\Delta^M$ est contenu dans $\Sigma(U_{P})$. Inversement, supposons ces relations v\'erifi\'ees. Toute racine $\alpha\in \Sigma(U_{P})$ est somme \`a coefficients positifs ou nuls d'\'el\'ements de $\Delta$ et au moins un \'el\'ement de $\Delta-\Delta^M$ intervient avec un coefficient strictement positif. Puisque $\beta(y)\geq0$ pour tout $\beta\in \Delta$, cela entra\^{\i}ne $\alpha(y)>N$. Donc $y\in X'_{N}(P)$. Si $y$ n'appartenait pas \`a $X_{N}(P)$, il y aurait un sous-groupe parabolique $P_{1}$ strictement contenu dans $P$ tel que $y\in X_{N}(P_{1})$. Alors $p_{M}(y)$ appartiendrait \`a $X_{N}^M(P_{1}\cap M)$. Puisque $P_{1}\cap M$ est un sous-groupe parabolique propre de $M$, on aurait $y\not\in X_{N}^M(M)$, contrairement \`a l'hypoth\`ese sur ${\cal F}_{M}$. Cela prouve (9). 

Remarquons qu'une racine $\alpha\in \Delta-\Delta^M$ n'est pas constante sur $p_{M}^{-1}(\bar{{\cal F}}_{M}^{\nu})$. Puisque $\rho_{Q}$ est en position g\'en\'erale, on a donc $\alpha(\rho_{Q})\not=0$. Notons $\Delta_{P}[Q]$  l'ensemble des $\alpha\in \Delta-\Delta^M$ tels que $\alpha(\rho_{Q})<0$. On d\'eduit imm\'ediatement de (9) que  le segment $]x,x+\epsilon\rho_{Q}]$ est contenu dans $p_{M}^{-1}(\bar{{\cal F}}_{M}^{\nu})$ pour $\epsilon>0$ assez petit si et seulement si les deux conditions suivantes sont v\'erifi\'ees

 $\alpha(x)\geq N$ pour tout $\alpha\in \Delta-\Delta^M$;
 
 $\alpha(x)>N$ pour tout $\alpha\in \Delta_{P}[Q]$.
 
 Notons $\Delta_{P,N}({\cal F})$ l'ensemble des $\alpha\in \Delta-\Delta^M$ telles que $\alpha(x)=N$ (il ne d\'epend que de ${\cal F}$,  pas du choix de $x$). Si la premi\`ere condition ci-dessus est v\'erifi\'ee, la seconde \'equivaut \`a $\Delta_{P,N}({\cal F})\cap \Delta_{P}[Q]=\emptyset$. Autrement dit
 
 (10)  le segment $]x,x+\epsilon\rho_{Q}]$ est contenu dans $p_{M}^{-1}(\bar{{\cal F}}_{M}^{\nu})$ pour $\epsilon>0$ assez petit si et seulement si les deux conditions suivantes sont v\'erifi\'ees

 $\alpha(x)\geq N$ pour tout $\alpha\in \Delta-\Delta^M$;
 
 $\Delta_{P,N}({\cal F})\cap \Delta_{P}[Q]=\emptyset$. 
 
 \bigskip

Il r\'esulte de (6) et (8) que, si ${\cal F}_{M}$ n'est pas contenue dans $X_{N}^M(M)$, alors $\Psi_{R}[Q]$ est vide. Donc $Z_{N,R}[Q]=0$. Cela \'etant vrai pour tout $Q$, on a $z_{N,R}(P,{\cal F}_{M},\nu)=0$, ce qui d\'emontre une partie du (ii) de l'\'enonc\'e. On suppose d\'esormais que ${\cal F}_{M}$ est contenue dans $X_{N}^M(M)$. Notons $\Psi_{R,N}$ l'ensemble des facettes ${\cal F}\in \Psi_{R}$ telles que, pour $x\in {\cal F}$, on a $\alpha(x)\geq N$ pour tout $\alpha\in \Delta-\Delta^M$ (avec la d\'efinition ci-dessus de $\Delta$, qui ne d\'epend pas de $Q$). 
Il r\'esulte de (6), (7) et (10) que $\Psi_{R}[Q]$ est l'ensemble des ${\cal F}\in \Psi_{R,N}$ telles que

$I_{Q}\subset I_{{\cal F}^{M,\nu}}$;

 $\Delta_{P,N}({\cal F})\cap \Delta_{P}[Q]=\emptyset$. 
 
 Chacune de ces conditions est ouverte car, en passant de ${\cal F}$ \`a une facette de son adh\'erence, $I_{{\cal F}^{M,\nu}}$ diminue tandis que $\Delta_{P,N}({\cal F})$ augmente. On va les convertir en conditions ferm\'ees. On remarque que, pour ${\cal F}\in \Psi_{R,N}$, on a
 
 $$ \sum_{J\subset I_{Q}; J\cap I_{{\cal F}^{M,\nu}}=\emptyset}(-1)^{\vert J\vert }=\left\lbrace \begin{array}{cc}1,& \text{si }\,\, I_{Q}\subset I_{{\cal F}^{M,\nu}},\\ 0,&\text{sinon;}\\ \end{array}\right.$$

  $$\sum_{D\subset \Delta_{P,N}({\cal F})\cap \Delta_{P}[Q]}(-1)^{\vert D\vert }=\left\lbrace \begin{array}{cc}1,& \text{si }\,\,\Delta_{P,N}({\cal F})\cap \Delta_{P}[Q]=\emptyset ,\\ 0,&\text{sinon;}\\ \end{array}\right.$$
  
  Pour des sous-ensembles $J\subset I_{Q}$ et $D\subset \Delta_{P}[Q]$, notons $\Psi_{R,N}[J,D]$ l'ensemble des ${\cal F}\in \Psi_{R,N}$ telles que $J\cap I_{{\cal F}^{M,\nu}}=\emptyset$ et $D\subset \Delta_{P,N}({\cal F})$. Alors la fonction caract\'eristique de $\Psi_{R}[Q]$ est la somme des fonctions caract\'eristiques des $\Psi_{R,N}[J,D]$ affect\'ees du coefficient $(-1)^{\vert J\vert +\vert D\vert }$. D'o\`u
  $$(11) \qquad Z_{N,R}[Q]=\sum_{J\subset I_{Q}}\sum_{D\subset \Delta_{P}[Q]}(-1)^{\vert J\vert +\vert D\vert }Z_{N,R}[J,D],$$
  o\`u
  $$Z_{N,R}[J,D]=\sum_{{\cal F}\in \Psi_{R,N}[J,D]}(-1)^{dim({\cal F}^{\nu})}.$$
  
Fixons $J$ et $D$. Il r\'esulte des d\'efinitions que $\Psi_{R,N}[J,D]$ est l'ensemble des ${\cal F}\in \Psi$ telles que ${\cal F}^{\nu}$ soit contenue dans un certain sous-ensemble $E[J,D]$ de $p_{M}^{-1}(\bar{{\cal F}}_{M}^{\nu})$. Pr\'ecis\'ement, c'est l'ensemble des $x\in App(A)$ tels que

(12)(a) $x\in p_{M}^{-1}(\bar{F}_{M}^{\nu})$;

(12)(b)  $x\in B_{R}(A)$;

 (12)(c) $p_{M}(x)$ appartient \`a une facette ${\cal F}_{M}^{_{'}\nu}$ de $\bar{{\cal F}}_{M}^{\nu}$ dont l'ensemble de points associ\'es $I_{{\cal F}^{_{'}\nu}}$ est disjoint de $J$;
 
  (12)(d) $\alpha(x)=N$ pour $\alpha\in D$ tandis que $\alpha(x)\geq N$ pour $\alpha\in \Delta-\Delta^M$ et $\alpha\not\in D$. 
  
Cet ensemble est ferm\'e et est r\'eunion de facettes ${\cal F}^{\nu}$. Le nombre $Z_{N,R}[J,D]$ appara\^{\i}t comme la caract\'eristique d'Euler-Poincar\'e de $E[J,D]$. 
Avec les notations introduites plus haut,  posons $J_{l}=I_{l}\cap J$ pour tout $l=1,...,k$. On a $J=\bigcup_{l=1,...,k}J_{l}$. Pour toute facette ${\cal F}_{M}^{_{'}\nu}$ de $\bar{{\cal F}}_{M}^{\nu}$, on a de m\^eme $I_{{\cal F}_{M}^{_{'}\nu}}=\bigcup_{l=1,...,k}I_{{\cal F}_{M}^{_{'}\nu},l}$ et chaque ensemble $I_{{\cal F}_{M}^{_{'}\nu},l}$ est non vide. Un tel ensemble ne peut \^etre disjoint de $J_{l}$ que si $J_{l}\not=I_{l}$. S'il existe $l$ tel que $J_{l}=I_{l}$, l'ensemble $ E[J,D]$ est donc vide et $Z_{R,N}[J,D]=0$. Notons ${\cal J}$ l'ensemble des sous-ensembles $J\subset I_{Q}$ tels que $J_{l}\not=I_{l}$ pour tout $l=1,...,k$.  Supposons maintenant  $J\in {\cal J}$. L'ensemble $I-J$ d\'efinit alors une facette de $\bar{{\cal F}}_{M}^{\nu}$. Fixons un point $y^M[J]$ de cette facette. Puisque les \'el\'ements de $\Delta-\Delta^M$ se restreignent en  des formes lin\'eaires sur ${\cal A}_{M}$ qui sont lin\'eairement ind\'ependantes, on peut fixer un point $y[J,D]\in App(A)$ tel que $p_{M}(y[J,D])=y^M[J]$ et  $\alpha(y[J,D])=N$ pour   $\alpha\in \Delta-\Delta^M$.    Ce point $y[J,D]$ v\'erifie toutes les conditions pour appartenir \`a $E[J,D]$, sauf peut-\^etre la condition $y[J,D]\in B_{R}(A)$. Mais, $N$ \'etant fix\'e, les donn\'ees sous-jacentes intervenant sont en nombre fini. C'est clair pour les donn\'ees $P$ et $Q$. Puisqu'on a suppos\'e ${\cal F}_{M}\subset X^M_{N}(M)$, il n'y a qu'un nombre fini de telles facettes. L'ensemble des $\nu$ peut \^etre infini mais ${\cal F}_{M}^{\nu}$ ne d\'epend en fait que de l'image de $\nu$ modulo le groupe $w_{G}(Z(G)(F))$ et ces images sont en nombre fini. Enfin, $J$ et $D$ sont en nombre fini. On n'a donc qu'un nombre fini de points $y[J,D]$. Si $R$ est assez grand, tous ces points appartiennnent \`a $B_{R}(A)$. Supposons qu'il en soit ainsi. Alors $y[J,D]\in E[J,D]$ et cet ensemble n'est pas vide.   Montrons que

(13)  $E[J,D]$ est convexe. 

  Les conditions (12)(a) et (12)(d) d\'efinissent \'evidemment des convexes. Il en est de m\^eme de (12)(b) car on a impos\'e \`a $B_{R}(A)$ d'\^etre convexe. On v\'erifie que, pour deux facettes ${\cal F}_{M}^{_{'}\nu}$ et ${\cal F}_{M}^{_{''}\nu}$ de $\bar{{\cal F}}_{M}^{\nu}$, un point appartenant \`a un segment joignant un point de la premi\`ere \`a un point de la seconde appartient \`a une facette ${\cal F}_{M,1}^{\nu}$ telle que $I_{{\cal F}_{M,1}^{\nu}}\subset I_{{\cal F}_{M}^{_{'}\nu}}\cup I_{{\cal F}_{M}^{_{''}\nu}}$. Donc (12)(c) d\'efinit aussi un convexe. Cela d\'emontre (13).
  
  Puisque $E[J,D]$ est un convexe non vide, sa caract\'eristique d'Euler-Poincar\'e vaut $1$, c'est-\`a-dire $Z_{N,R}[J,D]=1$. On obtient
 $$Z_{N,R}[Q]=\sum_{J\in{\cal J}}\sum_{D\subset \Delta_{P}[Q]}(-1)^{\vert J\vert +\vert D\vert }.$$
 La somme en $D$ vaut $1$ si $\Delta_{P}[Q]=\emptyset$  et $0$ sinon. La somme en $J$ se d\'ecompose en produit sur $l=1,...,k$ de sommes
 $$\sum_{J_{l}\subset I_{Q,l}, J_{l}\not=I_{l}} (-1)^{\vert J_{l}\vert }.$$
 Fixons $l$. Si $ I_{l} $ n'a qu'un \'el\'ement, alors $I_{Q,l}=I_{l}$ et le seul terme $J_{l}$ intervenant est l'ensemble vide. La somme ci-dessus vaut $1$. Si $I_{l}$ a au moins $2$ \'el\'ements, $I_{Q,l}$ est un sous-ensemble non vide de $I_{l}$ et on a remarqu\'e plus haut qu'il n'\'etait pas \'egal \`a $I_{l}$ tout entier. La condition $J_{l}\not=I_{l}$ est donc inutile dans la d\'efinition ci-dessus et la somme est nulle. On obtient que $Z_{R,N}[Q]$ vaut $1$ si $\Delta_{P}[Q]=\emptyset$ et $\vert I_{l}\vert =1$ pour tout $l$ et vaut $0$ sinon. La deuxi\`eme condition \'equivaut \`a ce que ${\cal F}_{M}^{\nu}$ soit r\'eduit \`a un point. Si elle est v\'erifi\'ee, on a $M_{{\cal F}_{M},\nu}=M$, c'est-\`a-dire $L=M$. Donc $Q$ et $P$ sont deux sous-groupes paraboliques de m\^eme composante de Levi. Puisque $\rho_{Q}$ est un point g\'en\'eral de la chambre associ\'ee \`a $Q$,  la condition $\Delta_{P}[Q]=\emptyset$ \'equivaut alors \`a $Q=P$. Donc $Z_{R,N}[Q]=1$ si ${\cal F}_{M}^{\nu}$ est r\'eduit \`a un point et $Q=P$ tandis que $Z_{R,N}[Q]=0$ sinon. La relation (5) entra\^{\i}ne alors que $z_{N,R}(P,{\cal F}_{M},\nu)$ vaut $mes(B_{M}\cap C_{P})mes(B_{M})^{-1}$ si ${\cal F}_{M}^{\nu}$ est r\'eduit \`a un point et $0$ sinon. Cela ach\`eve la d\'emonstration. $\square$
 
   \bigskip
   
   \section{D\'emonstration du th\'eor\`eme 12}
   Pour $f\in C_{c}(K_{\star}^{\dag}/H)$ et $P\in {\cal F}(M_{min})$, on a d\'efini en 16 le terme $B_{N,R}(P,f)$. 
   La proposition de ce paragraphe convertit la formule 16(2) en
   $$B_{N,R}(P,f)=mes(B_{M}\cap C_{P})mes(B_{M})^{-1}mes(K_{\star\star}^0\cap U_{P}(F))^{-1}$$
   $$\sum_{({\cal F}_{M},\nu)\in Fac^*_{max}(M;A)_{G-comp}, {\cal F}_{M}\subset X_{N}^M(M)}mes(K_{\star\star}^{M,0}\cap K_{{\cal F}_{M}}^0)^{-1}\int_{K_{{\cal F}_{M}}^{\nu}}f_{\star\star,[P]}(k)\phi_{{\cal F}_{M},\nu,cusp}(k)\,dk.$$
   Cela vaut pourvu que $R$ soit assez grand relativement \`a $N$. Le lemme 11 s'applique \`a $M$. Il montre que, si $N$ est assez grand, l'int\'egrale ci-dessus est nulle pour un \'el\'ement $({\cal F}_{M},\nu)\in Fac^*_{max}(M;A)_{G-comp}$ tel que ${\cal F}_{M}$ n'est pas inclus dans  $X_{N}^M(M)$. En choisissant $N$ assez grand (puis $R$ assez grand), on peut supprimer la condition ${\cal F}_{M}\subset X_{N}^M(M)$ dans la formule ci-dessus. On obtient une expression qui ne d\'epend plus de $N$ et que l'on note simplement $B(P,f)$. Posons $h=f_{\star\star,[P]}$. En imitant la d\'efinition de 14, on d\'efinit une fonction $h_{\star\star}$ sur $M(F)$ par
   $$h_{\star\star}(m)=\int_{K_{\star\star}^{M,0}}h(x^{-1}mx)\,dx.$$
   Mais il est clair que $h$ est par construction invariante par conjugaison par $K_{\star\star}^{M,0}$. Donc $h_{\star\star}=mes(K_{\star\star}^{M,0})h$. Alors
   $$B(P,f)=mes(B_{M}\cap C_{P})mes(B_{M})^{-1}mes(K_{\star\star}^0\cap U_{P}(F))^{-1}mes(K_{\star\star}^{M,0})^{-1}$$
   $$\sum_{({\cal F}_{M},\nu)\in Fac^*_{max}(M;A)_{G-comp} }mes(K_{\star\star}^{M,0}\cap K_{{\cal F}_{M}}^0)^{-1}\int_{K_{{\cal F}_{M}}^{\nu}}h_{\star\star}(k)\phi_{{\cal F}_{M},\nu,cusp}(k)\,dk$$
 $$=mes(B_{M}\cap C_{P})mes(B_{M})^{-1}mes(K_{\star\star}^0\cap P(F))^{-1} $$
   $$\sum_{({\cal F}_{M},\nu)\in Fac^*_{max}(M;A)_{G-comp} }mes(K_{\star\star}^{M,0}\cap K_{{\cal F}_{M}}^0)^{-1}\int_{K_{{\cal F}_{M}}^{\nu}}h_{\star\star}(k)\phi_{{\cal F}_{M},\nu,cusp}(k)\,dk$$ 
   o\`u la mesure sur $P(F)=M(F)U_{P}(F)$ est le produit des mesures sur $M(F)$ et sur $U_{P}(F)$.  
   Le m\^eme calcul qui a permis de convertir la formule 14(3) en (4), mais en sens inverse,  conduit \`a l'\'egalit\'e
  $$B(P,f)=mes(B_{M}\cap C_{P})mes(B_{M})^{-1}mes(K_{\star\star}^0\cap P(F))^{-1} $$
  $$\sum_{({\cal F}_{M},\nu)\in Fac^*_{max}(M)_{G-comp} } \int_{K_{{\cal F}_{M}}^{\nu}}h(k)\phi_{{\cal F}_{M},\nu,cusp}(k)\,dk.$$  
  Autrement dit, en utilisant les d\'efinitions,
  $$B(P,f)=mes(B_{M}\cap C_{P})mes(B_{M})^{-1}mes(K_{\star\star}^0\cap P(F))^{-1}mes(K_{\star\star}^{M,0})^{-1}  \Theta_{\pi,cusp}^M(f_{\star\star,[P]}).$$
  Remarquons que cette expression est ind\'ependante de  la mesure fix\'ee sur $P(F)$ puisque le dernier terme est proportionnel \`a cette mesure tandis que l'avant-dernier lui est inversement proportionnel. D'apr\`es 16(1), on a
  $$trace\,\pi(f)=\sum_{P\in {\cal F}(M_{min})}B(P,f).$$
  Pour tout $w\in W$, on peut remplacer ci-dessus $P$ par $w^{-1}(P)$ puis moyenner en $w$. D'o\`u
  $$trace\,\pi(f)=\vert W\vert ^{-1}\sum_{P\in {\cal F}(M_{min})}\sum_{w\in W}B(w^{-1}(P),f).$$
  Fixons $P$ et, pour tout $w\in W$, relevons $w$ en un \'el\'ement $n_{w}\in Norm_{G(F)}(A)$. Il est clair que $mes(B_{w^{-1}(M)}\cap C_{w^{-1}(P)})=mes(B_{M}\cap C_{P})$ et que $mes(K_{\star\star}^0\cap w^{-1}(P)(F))=mes(n_{w}K_{\star\star}^0n_{w}^{-1}\cap P(F))$. En d\'evissant les d\'efinitions, on v\'erifie que
  $$ \Theta^{w^{-1}(M)}_{\pi,cusp}(f_{\star\star,[w^{-1}(P)]})= \Theta^M_{\pi,cusp}((^{n_{w}}(f_{\star\star}))_{[P]}).$$ 
     D'o\`u
  $$B(w^{-1}(P),f)=mes(B_{M}\cap C_{P})mes(B_{M})^{-1}  \Theta^M_{\pi,cusp}( f[w]_{[P]}),$$
  o\`u $f[w]=mes(n_{w}K_{\star\star}^0n_{w}^{-1}\cap P(F))^{-1}{^{n_{w}}(f_{\star\star})}$. 
  Par d\'efinition, pour $g\in G(F)$, on a  
  $$ f[w](g)=mes(n_{w}K_{\star\star}^0n_{w}^{-1}\cap P(F))^{-1}\int_{K_{\star\star}^0}f(k^{-1}n_{w}^{-1}gn_{w}k)\,dk.$$
 Des mesures sur $G(F)$ et sur $P(F)$ se d\'eduit une mesure sur $P(F)\backslash G(F)$, cf. 11. On v\'erifie que, pour toute fonction $\varphi$ sur $G(F)$ v\'erifiant $\varphi(mug)=\delta_{P}(m)\varphi(g)$ pour tous $m\in M(F)$, $u\in U_{P}(F)$ et $g\in G(F)$, on a l'\'egalit\'e
 $$\int_{P(F)\backslash P(F)n_{w}K_{\star\star}^0}\varphi(g)\,dg =mes(n_{w}K_{\star\star}^0n_{w}^{-1}\cap P(F))^{-1}\int_{K_{\star\star}^0}\varphi(n_{w}k)\,dk.$$
 En cons\'equence, 
 $$\Theta^M_{\pi,cusp}( f[w]_{[P]})=\int_{P(F)\backslash P(F)n_{w}K_{\star\star}^0}\Theta^M_{\pi,cusp}((^gf)_{[P]})\,dg.$$
 La r\'eunion des ensembles $P(F)\backslash P(F)n_{w}K_{\star\star}^0$  quand $w$ d\'ecrit $W$ est $P(F)\backslash G(F)$ tout entier. Les ensembles associ\'es \`a $w$ et $w'$ sont disjoints ou confondus. Ils sont confondus si et seulement si $w\in W^Mw'$. On obtient que
 $$\sum_{w\in W}\Theta^M_{\pi,cusp}( f[w]_{[P]})=\vert W^M\vert  \int_{P(F)\backslash G(F)}\Theta^M_{\pi,cusp}((^gf)_{[P]})\,dg=\vert W^M\vert ind_{M}^G(\Theta^M_{\pi,cusp})(f).$$
 Puis
 $$\sum_{w\in W}B(w^{-1}(P),f)=\vert W^M\vert mes(B_{M}\cap C_{P})mes(B_{M})^{-1}ind_{M}^G(\Theta^M_{\pi,cusp})(f),$$
 et
 $$trace\,\pi(f)=\vert W\vert ^{-1}\sum_{P\in {\cal F}(M_{min})}\vert W^M\vert mes(B_{M}\cap C_{P})mes(B_{M})^{-1} ind_{M}^G(\Theta^M_{\pi,cusp})(f),$$
 o\`u, comme ci-dessus, $M$ d\'esigne la composante de Levi de $P$ contenant $M_{min}$. 
 On peut d\'ecomposer la somme en $P$ en une somme sur $M\in {\cal L}(M_{min})$ puis une somme sur $P\in {\cal P}(M)$. Pour $M$ fix\'e, la somme des $mes(B_{M}\cap C_{P})mes(B_{M})^{-1}$ quand $P$ d\'ecrit ${\cal P}(M)$ vaut $1$ puisque les chambres $C_{P}$ sont disjointes et que l'adh\'erence de leur r\'eunion est ${\cal A}_{M}$ tout entier. D'o\`u
 $$ trace\,\pi(f)=\sum_{M\in {\cal L}(M_{min})}\vert W^M\vert \vert W\vert ^{-1}ind_{M}^G(\Theta^M_{\pi,cusp})(f).$$
 Cela d\'emontre le th\'eor\`eme 12 pour une fonction $f\in C_{c}(K_{\star}^{\dag}/H)$.  Puisqu'on peut choisir des groupes $H$ v\'erifiant les conditions impos\'ees en 14 et aussi petits que l'on veut, le r\'esultat vaut pour toute fonction $f\in C_{c}^{\infty}(G(F))$ \`a support dans $K_{\star}^{\dag}$. Puisque ${\cal F}_{\star}$ \'etait une facette  quelconque dans $App(A)$ et que toutes les distributions ci-dessus sont invariantes par conjugaison, le r\'esultat vaut pour toute fonction $f\in C_{c}^{\infty}(G(F))$ telle qu'il existe une facette ${\cal F}\in Fac(G)$ de sorte que le support de $f$ soit contenu dans $K_{{\cal F}}^{\dag}$.  Consid\'erons maintenant une fonction $f\in C_{c}^{\infty}(G(F))$ dont le support est form\'e d'\'el\'ements compacts modulo $Z(G)$. Chacun de ces \'el\'ements est contenu dans $K_{{\cal F}}^{\dag}$ pour une certaine facette ${\cal F}$. Puisque le support de $f$ est compact et que les ensembles $K_{{\cal F}}^{\dag}$ sont ouverts, une partition de l'unit\'e permet d'\'ecrire $f$ comme une combinaison lin\'eaire finie de fonctions  $f'$ pour lesquelles il existe une facette ${\cal F}'$ de sorte que $f'$ soit \`a support dans $K_{{\cal F}'}^{\dag}$. L'\'egalit\'e ci-dessus vaut pour ces fonctions $f'$, donc aussi pour $f$. Cela ach\`eve la preuve du th\'eor\`eme 12.
 
 \bigskip
 
 \section{Expression du caract\`ere \`a l'aide d'int\'egrales orbitales pond\'er\'ees}
 On sait que la distribution $\Theta_{\pi}$ est localement int\'egrable, plus pr\'ecis\'ement, elle est associ\'ee \`a une fonction  $\theta_{\pi}$ sur $G(F)$ qui est localement int\'egrable sur $G(F)$ et localement constante sur le sous-ensemble des \'el\'ements semi-simples fortement r\'eguliers. La formule de Weyl permet  d'\'ecrire
$$(1) \qquad \Theta_{\pi}(f)=\sum_{M\in {\cal L}(M_{min})}\vert W^M\vert \vert W^G\vert ^{-1} \int_{M(F)_{ell}}f_{P_{M}}(m)\theta_{\pi}(m)D^G(m)^{1/2}dm$$
pour toute $f\in C_{c}^{\infty}(G(F))$, o\`u, pour tout $M$, $P_{M}$ est un \'el\'ement quelconque de ${\cal P}(M)$. 

En utilisant les formules  9(1) et 11(1), on a l'\'egalit\'e
$$\Theta_{\pi,cusp}(f)=\sum_{({\cal F},\nu)\in \underline{Fac}^*_{max}(G)}mes(A_{G}(F)\backslash K_{{\cal F}}^{\dag})^{-1}\sum_{M\in {\cal L}(M_{min})}\vert W^M\vert \vert W^G\vert ^{-1} (-1)^{dim(A_{M})-dim(A_{G})}$$
$$\int_{M(F)_{ell}}f_{P_{M}}(m) J_{M}^G(m,\phi_{{\cal F},\nu,cusp})\,dm.$$
Pour $L\in {\cal L}(M_{min})$, on a une formule similaire pour la distribution $\Theta^L_{\pi,cusp}$.  
  Pour $M\subset L$ deux sous-groupes de Levi contenant $M_{min}$, pour $P\subset Q$ deux sous-groupes paraboliques de composantes de Levi $M$, resp. $L$, on a l'\'egalit\'e $(f_{Q})_{P\cap L}=f_{P}$ pour tout $f\in C_{c}^{\infty}(G(F))$. On en d\'eduit
$$ind_{L}^G(\Theta^L_{\pi,cusp})(f)=\sum_{({\cal F}_{L},\nu)\in \underline{Fac}^*_{max}(L)_{G-comp}}mes(A_{L}(F)\backslash K_{{\cal F}_{L}}^{\dag})^{-1}\sum_{M\in {\cal L}^L(M_{min})}\vert W^M\vert \vert W^L\vert ^{-1}$$
$$(-1)^{dim(A_{M})-dim(A_{L})} \int_{M(F)_{ell}}f_{P_{M}}(m) J_{M}^L(m,\phi_{{\cal F}_{L},\nu,cusp})\,dm$$
pour toute $f\in C_{c}^{\infty}(G(F))$. 
Donc
$$(2) \qquad \sum_{L\in {\cal L}(M_{min})}\vert W^L\vert \vert W^G\vert ^{-1} ind_{L}^G(\Theta^L_{\pi,cusp})(f)=\sum_{M\in {\cal L}^L(M_{min})}\vert W^M\vert \vert W^L\vert ^{-1} $$
$$\int_{M(F)_{ell}}f_{P_{M}}(m)\theta^M_{\pi, comp}(m)D^G(m)^{1/2}\,dm,$$
o\`u, pour $M\in {\cal L}(M_{min})$ et $m\in M(F)_{ell}$, on a pos\'e
$$\theta^M_{\pi,comp}(m)=D^G(m)^{-1/2}\sum_{L\in {\cal L}(M)} (-1)^{dim(A_{M})-dim(A_{L})}$$
$$\sum_{({\cal F}_{L},\nu)\in \underline{Fac}^*_{max}(L)_{G-comp}}  mes(A_{L}(F)\backslash K_{{\cal F}_{L}}^{\dag})^{-1}J_{M}^L(m,\phi_{{\cal F}_{L},\nu,cusp}).$$
Le th\'eor\`eme 12 dit que les deux expressions (1) et (2) sont \'egales si $f$ est \`a support compact modulo $Z(G)$. On en d\'eduit le calcul suivant pour le caract\`ere $\theta_{\pi}$ restreint aux \'el\'ements semi-simples fortement r\'eguliers et compacts modulo $Z(G)$. Soit $x$ un tel \'el\'ement. Quitte \`a conjuguer $x$, on peut supposer qu'il existe $M\in {\cal L}(M_{min})$ de sorte que $x$ appartienne \`a $M(F)$ et soit elliptique dans $M$. On a alors  l'\'egalit\'e 

$$(3) \qquad \theta_{\pi}(x)=\theta^M_{\pi,comp}(x).$$ 

On peut interpr\'eter la d\'efinition de la fonction $\theta^M_{\pi,comp}$ de la fa\c{c}on suivante.   Fixons une fonction $\underline{P}:{\cal L}(M_{min})\to {\cal F}(M_{min})$ telle que $\underline{P}(L)\in {\cal P}(L)$ pour tout $L\in {\cal L}(M_{min})$.   Pour tout $M\in {\cal L}(M_{min})$, on d\'efinit une fonction $\theta^M_{\pi, \underline{P}}$  sur l'ensemble des \'el\'ements semi-simples fortement r\'eguliers de $M(F)$ par
$$\theta^M_{\pi,\underline{P}}(m)=D^G(m)^{-1/2}\sum_{L\in {\cal L}(M)} (-1)^{dim(A_{M})-dim(A_{L})}$$
$$\sum_{({\cal F}_{L},\nu_{L})\in \underline{Fac}^*_{max}(L)}  mes(A_{L}(F)\backslash K_{{\cal F}_{L}}^{\dag})^{-1}J_{M}^L(m,\phi_{\pi_{\underline{P}(L)},{\cal F}_{L},\nu_{L},cusp}).$$
On a:

(4) si $m$ est compact modulo $Z(G)$, alors $\theta^M_{\pi,comp}(m)=\theta^M_{\pi,\underline{P}}(m)$.

 La preuve est la m\^eme qu'en 11(2).

Soit $g\in G(F)$ un \'el\'ement semi-simple et fortement r\'egulier. D'apr\`es Casselman, cf. \cite{Ca} paragraphe 2, on lui associe un sous-groupe parabolique $Q_{g}$ de $G$ et une composante de Levi $L_{g}
$ de $P_{g}$ de la fa\c{c}on suivante. Notons $T$ le sous-tore maximal de $G$ contenant $g$. Pla\c{c}ons-nous pour un instant sur une cl\^oture alg\'ebrique $\bar{F}$ de $F$. Prolongeons la valeur absolue $\vert .\vert _{F}$  \`a $\bar{F}$. Notons $\Sigma_{T}$ l'ensemble des racines de $T$ dans  $G$ sur $\bar{F}$. On note $\Sigma_{T}^1$, resp. $\Sigma_{T}^<$, l'ensemble des $\alpha\in \Sigma_{T}$ telles que $\vert \alpha(g)\vert _{F}=1$, resp. $\vert \alpha(g)\vert < 1$. Alors il existe un unique  sous-groupe parabolique $Q_{g}$ et une unique composante de Levi $L_{g}$ de $Q_{g}$ contenant $T$, d\'efinis tous deux sur $\bar{F}$, de sorte que $\Sigma_{T}^1$, resp. $\Sigma_{T}^<$, soit l'ensemble des racines de $T$ dans $L_{g}$, resp. dans $U_{Q_{g}}$. On v\'erifie que ces deux groupes sont en fait d\'efinis sur $F$. L'\'el\'ement $g\in M_{g}(F)$ est compact modulo $Z(M_{g})$. Quitte \`a conjuguer $g$, on peut supposer qu'il existe $M\in {\cal L}(M_{min})$ de sorte que $g\in M(F)_{ell}$. Un tel \'el\'ement est compact modulo $Z(M)$. Il en r\'esulte ais\'ement que $M\subset L_{g}$, a fortiori $L_{g}\in {\cal L}(M_{min})$. 
Casselman a prouv\'e que $\theta_{\pi}(g)=\delta_{P_{g}}(g)^{1/2} \theta_{\pi_{P_{g}}}(g)$, cf. \cite{Ca} th\'eor\`eme 5.2. Ce dernier terme est calcul\'e par la formule (3) o\`u l'on remplace $G$ et $\pi$ par $L_{g}$ et $\pi_{P_{g}}$. Cela d\'emontre le th\'eor\`eme suivant.
 
\ass{Th\'eor\`eme}{Soit $M\in {\cal L}(M_{min})$ et soit $m\in M(F)_{ell}$. Alors on a l'\'egalit\'e
$$\theta_{\pi}(m)=\delta_{P_{m}}(m)^{1/2}\theta^M_{\pi_{Q_{m}},comp}(m).$$}

Pour \'eviter les confusions, \'ecrivons explicitement le terme apparaissant ci-dessus:
$$\theta^M_{\pi_{Q_{m}},comp}(m)=D^{L_{m}}(m)^{-1/2}\sum_{L\in {\cal L}^{L_{m}}(M)} (-1)^{dim(A_{M})-dim(A_{L})}$$
$$\sum_{({\cal F}_{L},\nu)\in \underline{Fac}^*_{max}(L)_{L_{m}-comp}}  mes(A_{L}(F)\backslash K_{{\cal F}_{L}}^{\dag})^{-1}J_{M}^L(m,\phi_{\pi_{Q_{m}},{\cal F}_{L},\nu,cusp}).$$

De m\^eme qu'en (4), on peut reformuler la d\'efinition de $\theta^M_{\pi_{Q_{m}},comp}(m)$ de la fa\c{c}on suivante. Fixons une fonction $\underline{P}^{Q_{m}}:{\cal L}^{L_{m}}(M_{min})\to {\cal F}(M_{min})$ telle que, pour tout $L\in {\cal L}^{L_{m}}(M_{min})$, $\underline{P}^{Q_{m}}(L)\in {\cal P}(L)$ et $\underline{P}^{Q_{m}}(L)\subset Q_{m}$.   On a alors
$$\theta^M_{\pi_{Q_{m}},comp}(m)=D^{L_{m}}(m)^{-1/2}\sum_{L\in {\cal L}^{M_{m}}(M)} (-1)^{dim(A_{M})-dim(A_{L})}$$
$$\sum_{({\cal F}_{L},\nu_{L})\in \underline{Fac}^*_{max}(L)}  mes(A_{L}(F)\backslash K_{{\cal F}_{L}}^{\dag})^{-1}J_{M}^L(m,\phi_{\pi_{\underline{P}^{Q_{m}}(L)},{\cal F}_{L},\nu_{L},cusp}).$$

\bigskip

CNRS-Institut de Math\'ematiques de Jussieu

4 place Jussieu, 75005 Paris

jean-loup.waldspurger@imj-prg.fr
\end{document}